\documentclass[11pt]{article}
\usepackage{amssymb,latexsym,times}
\usepackage{mathptmx}
\usepackage{amscd}
\usepackage[all]{xy}

\newtheorem{Thm}{Theorem}[section]
\newtheorem{Lem}[Thm]{Lemma}

\newtheorem{Prop}[Thm]{Proposition}
\newtheorem{Conj}[Thm]{Conjecture}
\newtheorem{example}[Thm]{Example}
\newtheorem{remark}[Thm]{Remark}
\newtheorem{Def}[Thm]{Definition}

\newcommand{\Z}{\mathbb{Z}}
\newcommand{\Q}{\mathbb{Q}}

\newcommand{\N}{\mathbb{N}}
\newcommand{\C}{\mathbb{C}}

\newcommand{\Hom}{\operatorname{Hom}}

\newcommand{\dimv}{\underline{\dim}}

\newcommand{\bsm}{\begin{smallmatrix}}
\newcommand{\esm}{\end{smallmatrix}}

\newcommand{\bg}{\mathbf{g}}

\def\proof{\medskip\noindent {\it Proof --- \ }}
\def\cqfd{\hfill $\Box$ \bigskip}
\def\CC{{\mathcal C}}
\def\resp{{\em resp.\ }}

\def\<{\langle\,}
\def\>{\,\rangle}

\def\eg{{\em e.g. }}

\def\g{\mathfrak g}

\def\<{\langle}
\def\>{\rangle}

\def\1{\mathbf 1}
\def\Gr{{\rm Gr}}

\def\a{\alpha} 
\def\Y{\mathcal{Y}}

\def\L{\Lambda}

\def\Hom{{\rm Hom}}

\def\de{\delta}
\def\De{\Delta}

\def\AA{\mathcal{A}}

\def\YY{\mathcal{Y}}

\def\1{{\mathbf 1}}
\def\ie{{\em i.e. }}
\def\hg{\widehat{\mathfrak{g}}}

\def\G{\Gamma}
\def\tG{\widetilde{\Gamma}}

\def\b{\beta}
\def\ga{\gamma}

\def\bz{\mathbf z}

\def\by{\mathbf{y}}
\def\btrois{3}
\def\bdeux{2}
\def\bun{1}

\def\bquat{4}

\def\bY{\mathbf{Y}}
\def\Si{\Sigma}
\def\SS{\mathcal{S}}
\def\dh{h\,\check{}\,}
\def\ds{\displaystyle}
\def\bW{\mathbf{W}}
\def\bV{\mathbf{V}}
\def\Ker{\mathrm{Ker\,}}

\begin{document}

\title{\bf A cluster algebra approach to $q$-characters of Kirillov-Reshetikhin modules}
\author{D. Hernandez, B. Leclerc}

\date{}

\maketitle

\begin{abstract}
We describe a cluster algebra algorithm for calculating $q$-characters of 
Kirillov-Resheti\-khin modules for any untwisted quantum affine algebra $U_q(\hg)$.
This yields a geometric $q$-character formula for tensor products 
of Kirillov-Reshetikhin modules. When $\g$ is of type $A, D, E$,
this formula extends Nakajima's formula for $q$-characters of 
standard modules in terms of homology of graded quiver varieties. 
\end{abstract}

\setcounter{tocdepth}{1}
{\footnotesize \tableofcontents}

\section{Introduction}

Let $\g$ be a simple Lie algebra over $\C$, and let $U_q(\hg)$ be the corresponding
untwisted quantum affine algebra with quantum parameter $q\in\C^*$ not a root of unity.
The finite-dimensional complex representations of $U_q(\hg)$ have been studied by many
authors during the past twenty years. We refer the reader to \cite{CP} for a 
classical introduction, and to \cite{CH,L} for recent surveys on this topic.

In \cite{HL}, we started to explore some new connections between this rich representation
theory and the cluster algebras of Fomin and Zelevinsky.
The main result, proved in \cite{HL} in type $A_n$ and $D_4$, and extended to any
$A$-$D$-$E$ type by Nakajima \cite{N3}, shows the existence of a tensor category $\CC_1$ of
finite-dimensional $U_q(\hg)$-modules whose Grothendieck ring is a cluster algebra of the 
same finite Dynkin type, such that the classes of simple modules coincide with the set of
cluster monomials. As a consequence, the $q$-characters of the simple objects of $\CC_1$ 
can be computed algorithmically using the combinatorics
of cluster algebras. Moreover, the Caldero-Chapoton formula for cluster expansions leads to
some new geometric formulae for these characters, in terms of Euler characteristics of
quiver Grassmannians.

Unfortunately the category $\CC_1$ is quite small. For instance it contains only three
Kirillov-Resheti\-khin modules for each node of the Dynkin diagram of $\g$.
Another limitation of the papers \cite{HL} and \cite{N3} is that $\g$ is assumed to
be of simply laced type. In fact, the general proof of Nakajima uses in a crucial way
his geometric construction of the standard $U_q(\hg)$-modules \cite{N1}, which is only 
available in the simply laced case.

In this paper we drop the assumption of being simply laced, and we consider a much larger tensor 
subcategory $\CC^-$ which contains, up to spectral shifts, all the irreducible finite-dimensional 
representations of $U_q(\hg)$. Our first main result (Theorem~\ref{thm1}) is an algorithm
which calculates the $q$-character of an arbitrary Kirillov-Reshetikhin module in $\CC^-$
as the result of a sequence of cluster mutations. The only input for this calculation is
the initial seed of our cluster algebra $\AA$, which is encoded in a quiver obtained 
from the Cartan matrix of $\g$ by a simple and uniform recipe. (It may be worth noting that
$\AA$ is always a skew-symmetric cluster algebra, even when $\g$ is not simply laced.)

The proof of this theorem is based on the fact that the $q$-characters of the Kirillov-Reshetikhin 
modules are solutions of the corresponding $T$-system of Kuniba, Nakanishi and Suzuki \cite{KNS1,N2,H}.
This will come as no surprise, given the many papers already devoted to the relationships between 
cluster algebras and $T$-systems (see in particular \cite{IIKNS}, \cite{IIKKN1}, \cite{IIKKN2};
in fact our algorithm is inspired from \cite[\S13]{GLS2}, where similar $T$-system formulas
are obtained for generalized minors of symmetric Kac-Moody groups).
We find it nevertheless remarkable that, by interpreting the $T$-system equations as appropriate
cluster transformations, one is able to obtain the Kirillov-Reshetikhin $q$-characters starting from
their highest weight monomials via a procedure of successive approximations.
To the best of our knowledge this simple ``bootstrap'' algorithm  had not been noticed before,
although, in retrospect, it could certainly 
have been formulated and proved without knowledge of the cluster algebra theory.

At this stage, we should recall that Frenkel and Mukhin \cite{FM} have described long ago a completely
different algorithm, which can be used for computing the $q$-characters of the Kirillov-Reshetikhin modules \cite{N2,H}.
The advantage of our approach is that we are now in a position to apply  deep results of the theory
of cluster algebras and obtain new formulas for the Kirillov-Reshetikhin $q$-characters.
In \cite{DWZ,DWZ2}, Derksen, Weyman and Zelevinsky have constructed a categorical model for a large class
of cluster algebras using quivers with potentials. In particular they have proved a far-reaching
generalization of the Caldero-Chapoton formula, expressing any cluster variable in terms of the $F$-polynomial
of an associated quiver representation (see also \cite{P1} for a different proof of this generalized formula). 
Applying this formula in our context, we get a geometric
character formula for arbitrary Kirillov-Reshetikhin modules, and also for their tensor products
(Theorem~\ref{th_geom_form}). 

When $\g$ is simply laced, and we restrict our attention to the simplest
Kirillov-Reshetikhin modules and their tensor products, namely the fundamental modules and 
the standard modules, the quiver Grassmannians involved in our formula are homeomorphic to the projective
varieties $\mathfrak{L}^\bullet(V,W)$ used by Nakajima \cite[\S4]{NAnnals} in his geometric construction of the standard modules. 
This suggests that the quiver Grassmannians we introduce, in connection with general Kirillov-Reshetikhin modules of  
not necessarily simply laced type, might be interesting new varieties. 

When $\g$ is a classical Lie algebra of type $A$, $B$, $C$, $D$, there exist tableau sum formulas for 
the $q$-characters of certain Kirillov-Reshetikhin modules (see \cite[\S7]{KNS2} and references therein). 
From the geometric point of view of Theorem~\ref{th_geom_form},
these formulas can be explained by the fact that the corresponding quiver representations
have a nice and regular ``grid structure'', and in many cases their quiver Grassmannians
are reduced to points (see \eg \S\ref{append-B2}, \S\ref{append-B3}, \S\ref{append-C3}).

The cluster algebra approach also suggests that our results should extend far beyond
the Kirillov-Reshetikhin modules.
Indeed, we show (Theorem~\ref{thm_Grothendieck_ring}) that the cluster algebra $\AA$
is isomorphic to the Grothendieck ring of $\CC^-$. It is then natural to conjecture
that this isomorphism maps all cluster monomials of $\AA$ to the classes of certain simple
objects of $\CC^-$ (Conjecture~\ref{conj1}), and to extend the above geometric character
formula to all these simple objects (Conjecture~\ref{Conj2}).
The results of \cite{HL,HL2} and \cite{N3} provide some evidence supporting
these conjectures in the simply laced case.

Here is a more precise outline of the paper. 
In Section~\ref{section1} we associate with every simple Lie algebra $\g$ some quivers (\S\ref{subsect-quiver}), 
from which we define a cluster algebra $\AA$ (\S\ref{subsect-cluster}). We also introduce the untwisted
quantum affine algebra $U_q(\hg)$ (\S\ref{subsect-qaff}).
In Section~\ref{section2} we state and prove our algorithm for computing the 
Kirillov-Reshetikhin $q$-characters as special cluster variables of $\AA$. 
The proof uses $T$-systems (\S\ref{ssect-T-syst}) 
and the notion of truncated $q$-characters (\S\ref{ssect-truncated}).
In Section~\ref{section3}, we consider an algebra $A$ defined by a quiver
with potential, coming from our initial seed for~$\AA$~(\S\ref{potential}).
We introduce certain distinguished $A$-modules $K^{(i)}_{k,m}$ (\S\ref{ssect-A-mod}), 
and we state our geometric formula for the Kirillov-Reshetikhin $q$-characters
in terms of Grassmannians of submodules of the $K^{(i)}_{k,m}$
(Theorem~\ref{th_geom_form}).
To prove it, we calculate the $g$-vectors of these $q$-characters, regarded
as cluster variables of $\AA$, and we apply a result of Plamondon \cite{P} which allows
to reconstruct the $A$-module corresponding to a given cluster variable
from the knowledge of its $g$-vector. To be in a position 
to apply this result, we show that the defining potential of
$A$ is rigid, and that appropriate truncations of $A$
are finite-dimensional (Proposition~\ref{finite+rigid}).
In Section~\ref{section4}, we prove Theorem~\ref{thm_Grothendieck_ring} and 
we formulate Conjecture~\ref{conj1} and Conjecture~\ref{Conj2}.  
The paper closes with an appendix 
illustrating our results with many examples.

\section{Definitions and notation}\label{section1}
\subsection{Quivers}\label{subsect-quiver}
\subsubsection{Cartan matrix}
Let $C=(c_{ij})_{i,j\in I}$ be an indecomposable $n\times n$ Cartan matrix of finite type \cite[\S4.3]{K}.
There is a diagonal matrix $D = \mbox{diag}(d_i\mid i\in I)$ with entries in $\Z_{>0}$ 
such that the product
\[
B=D\,C = (b_{ij})_{i,j\in I}
\]
is symmetric. We normalize $D$ so that $\min\{d_i \mid i\in I\} = 1$, and
we put $t:=\max\{d_i \mid i\in I\}$.
Thus 
\[
t =
\left\{
\begin{array}{cl}
1 & \mbox{if $C$ is of type $A_n$, $D_n$, $E_6$, $E_7$ or $E_8$}, \cr
         2 & \mbox{if $C$ is of type $B_n$, $C_n$ or $F_4$}, \cr
         3 & \mbox{if $C$ is of type $G_2$}.
\end{array}
\right.
\]
It is easy to check by inspection that
\begin{equation}\label{easy-prop}
(d_i > 1 \mbox{ and } c_{ij} < 0) \Longrightarrow (c_{ij} = -1). 
\end{equation}
One attaches to $C$ a Dynkin diagram $\de$ with vertex set $I$ \cite[\S 4.7]{K}.
Since $C$ is assumed to be indecomposable and of finite type, $\de$ is a tree. 

All the objects that we consider below depend on $C$,
but we shall not always repeat it, neither record it explicitly in our notation. 

\begin{example}\label{example_Cartan}
{\rm
The Cartan matrix $C$ of type $B_3$ in the Cartan-Killing classification is
defined by
\[
C = \pmatrix{
2 & -1 & 0 \cr
-1 & 2 & -1\cr
0 & -2 & 2}
\]
We have $D=\mbox{diag}(2,2,1)$ and the symmetric matrix $B$ is given by
\[
B = \pmatrix{
4 & -2 & 0 \cr
-2 & 4 & -2\cr
0 & -2 & 2}
\]
}
\end{example}

\subsubsection{Infinite quiver} \label{subsubsect_doubly}

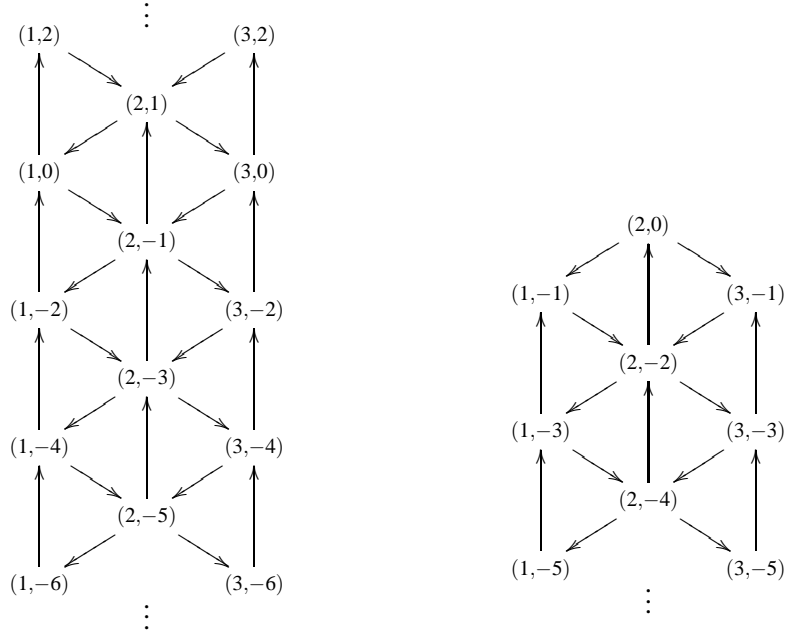
\begin{figure}[t]
\[
\def\objectstyle{\scriptstyle}
\def\lablestyle{\scriptstyle}
\xymatrix@-1.0pc{
&&&&\\
&{(1,2)}\ar[rd]&{}\save[]+<0cm,2ex>*{\vdots}\restore
&\ar[ld] (3,2) 
\\
&&\ar[ld] (2,1) \ar[rd]&&
\\
&\ar[uu]{(1,0)}\ar[rd]&
&\ar[ld] (3,0) \ar[uu]
\\
&&\ar[uu]\ar[ld] (2,-1) \ar[rd]&&
\\
&\ar[uu](1,-2) \ar[rd] &&\ar[ld] (3,-2)\ar[uu]
\\
&&\ar[ld] \ar[uu](2,-3) \ar[rd]&&
\\
&\ar[uu](1,-4) \ar[rd] &&\ar[ld] (3,-4)\ar[uu]
\\
&&\ar[ld] \ar[uu](2,-5) \ar[rd]&&
\\
&\ar[uu](1,-6) &{}\save[]+<0cm,-2ex>*{\vdots}\restore&\ar[uu] (3,-6) 
\\
}
\qquad\qquad
\xymatrix@-1.0pc{
&&&&\\
&&{}& 
\\
&&&&
\\
&&&
\\
&&&
\\
&&\ar[ld] (2,0) \ar[rd]&&
\\
&(1,-1) \ar[rd] &&\ar[ld] (3,-1)
\\
&&\ar[ld] \ar[uu](2,-2) \ar[rd]&&
\\
&\ar[uu](1,-3) \ar[rd] &&\ar[ld] (3,-3)\ar[uu]
\\
&&\ar[ld] \ar[uu](2,-4) \ar[rd]&&
\\
&\ar[uu](1,-5) &{}\save[]+<0cm,-2ex>*{\vdots}\restore&\ar[uu] (3,-5) 
}
\]
\caption{\label{Fig0} {\it The quivers $\G$ and $G^-$ in type $A_3$.}}
\end{figure}
\begin{figure}[t]
\[
\def\objectstyle{\scriptstyle}
\def\lablestyle{\scriptstyle}
\xymatrix@-1.0pc{
&&&&\\
&&{}\save[]+<0cm,0ex>*{\vdots}\restore
& {}\save[]+<0cm,0ex>*{\vdots}\restore
\\
&&\ar[ld] (\bdeux,5)& \ar[ld](\bun,5)&
\\
&{(\bun,3)}\ar[rd]\ar[uu] &\ar[u] (\bdeux,3) \ar[rd]&&
\\
&&\ar[u]\ar[ld] (\bdeux,1) &\ar[ld] (\bun,1) \ar[uu]  &
\\
&\ar[uu](\bun,-1)\ar[rd]   & \ar[u](\bdeux,-1) \ar[rd]&&
\\
&&\ar[u]\ar[ld] (\bdeux,-3) &\ar[ld] (\bun,-3) \ar[uu] &
\\
&\ar[uu]\ar[rd](\bun,-5) & \ar[u]\ar[rd](\bdeux,-5)&&
\\
&&\ar[u]\ar[ld] (\bdeux,-7) &\ar[ld] (\bun,-7) \ar[uu]  &
\\
&\ar[uu](\bun,-9)   & \ar[u](\bdeux,-9) &&
\\
&{}\save[]+<0cm,0ex>*{\vdots}\restore  &{}\save[]+<0cm,0ex>*{\vdots}\restore&{}\save[]+<0cm,0ex>*{\vdots}\ar[uu]\restore 
\\
}
\qquad\qquad
\xymatrix@-1.0pc{
&&&&\\
&&& \\
&&&&\\
&&&&\\
&&&&\\
&&&&\\
&  & (\bdeux,0) \ar[rd]&&
\\
&&\ar[u]\ar[ld] (\bdeux,-2) &\ar[ld] (\bun,-1) &
\\
&\ar[rd](\bun,-3) & \ar[u]\ar[rd](\bdeux,-4)&&
\\
&&\ar[u]\ar[ld] (\bdeux,-6) &\ar[ld] (\bun,-5) \ar[uu]  &
\\
&\ar[uu](\bun,-7)   & \ar[u](\bdeux,-8) &&
\\
&{}\save[]+<0cm,0ex>*{\vdots}\restore  &{}\save[]+<0cm,0ex>*{\vdots}\restore&{}\save[]+<0cm,0ex>*{\vdots}\ar[uu]\restore 
\\
}
\]
\caption{\label{Fig1} {\it The quivers $\G$ and $G^-$ in type $B_2$.}}
\end{figure}
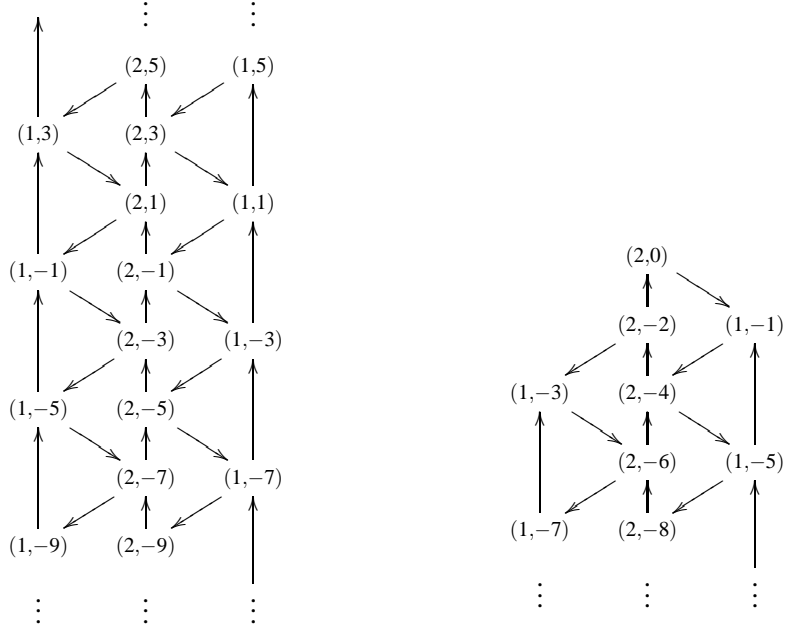

Put $\widetilde{V} = I \times \Z$.
We introduce a quiver $\tG$ with vertex set $\widetilde{V}$.
The arrows of $\tG$ are given by
\[
((i,r) \to (j,s)) 
\quad \Longleftrightarrow \quad
(b_{ij}\not = 0 
\quad \mbox{and} \quad
s=r+b_{ij}).
\]

\begin{Lem}
The quiver $\tG$ has two isomorphic connected components.
\end{Lem}
\proof
Let $i\in I$ be such that $d_i=1$. For every $r\in\Z$ we have an arrow
$(i,r) \to (i,r+2)$. Since the Dynkin diagram $\de$ is connected, every vertex $(j,s)\in\widetilde{V}$
is connected to a vertex of the form $(i,r)$, so $\tG$ has at most 
two connected components. Moreover, since $\de$ is a tree, any path 
from $(i,r)$ to $(i,s)$ in $\tG$ contains as many arrows of the form
$(j,p)\to (k,p+b_{jk})$ with $j\not =k$, as it contains arrows of the form $(k,t) \to (j,t+b_{kj})$.
Since $b_{jk}=b_{kj}$, and since $b_{jj} \in 2\Z$ for every $j\in I$, 
it follows that if there is a path from
$(i,r)$ to $(i,s)$ then $s-r\in 2\Z$. Therefore $\tG$ has exactly two connected components.
These two components are isomorphic via the map $(j,r) \mapsto (j,r+1) \ ((j,r) \in \widetilde{V}\times \Z))$.
\cqfd

We pick one of the two isomorphic connected components of $\tG$ and call it
$\G$. The vertex set of $\G$ is denoted by $V$.

\subsubsection{Semi-infinite quiver} \label{subsubsect_simply}

We will have to use a second labelling of the vertices of $\G$.
It is deduced from the first one by means of the function $\psi$ defined by
\begin{equation}\label{def_psi}
 \psi(i,r) = (i, r+d_i),\qquad ((i,r)\in V).
\end{equation}
Let $W\subset I \times \Z$ be the image of $V$ under $\psi$.
We shall denote by $G$ the same quiver as $\G$ 
but with vertices labelled by $W$.
Put $W^- := W \cap (I \times \Z_{\le 0})$.
Let $G^-$ be the full subquiver of $G$ with vertex set $W^-$.

\begin{example}\label{example_quivers}
{\rm
The definitions of \S\ref{subsubsect_doubly} and \S\ref{subsubsect_simply} are illustrated 
in Figure~\ref{Fig0} and Figure~\ref{Fig1}.
We find it convenient to always display the quivers $\G$ in the following way.
We decide to draw all arrows of the form $(i,r) \to (i,r+b_{ii})$ vertically, going upwards. 
Moreover, if $(i,r)$ and $(i,s)$ 
are two vertices with $r-s\not \in b_{ii}\Z$, we draw them in different \emph{columns}.
Hence, the quivers attached to $C$ always have $\sum_{i\in I}d_i$ columns.
Finally, the integer $r$ determines the \emph{altitude} of the vertex $(i,r)$ in $\G$.
Therefore, since for $i\not=j$ we have $b_{ij}\le 0$, the arrows $(i,r) \to (j,r+b_{ij})$
are represented as oblique arrows going down.

Figure~\ref{Fig0} displays the quivers $\G$ and $G^-$ for $C$ of type $A_3$.
Figure~\ref{Fig1} shows
$\G$ and $G^-$ for $C$ of type $B_2$.
In both cases  
we have chosen to call $\G$ the connected component of~$\widetilde{\G}$
containing the vertex $(2,1)$.
For another illustration, with $C$ of type $G_2$, see Figure~\ref{Fig2}.
More examples can be found in the Appendix, \S\ref{append-B3}, \S\ref{append-C3}, \S\ref{append-F4}. 
}
\end{example}

\subsection{Cluster algebras}\label{subsect-cluster}

We refer the reader to \cite{FZsurvey} and \cite{GSV} for an introduction
to cluster algebras, and for any un\-defined terminology.

\subsubsection{Cluster algebra attached to $G^-$} \label{subsect_ca}
Consider an infinite set of indeterminates 
$
\bz^-=\{z_{i,r}\mid (i,r)\in W^-\}
$
over $\Q$.
Let $\AA$ be the cluster algebra defined by the initial seed 
$(\bz^-, G^-)$.
Thus,  $\AA$ is the $\Q$-subalgebra of the field of rational functions
$\Q(\bz^-)$ generated by all the elements obtained from
some element of $\bz^-$ via a finite sequence of seed mutations, see \cite[Definition 3.1]{GG}.
Note that there are no frozen variables.

Cluster algebras of infinite rank have not received much attention
up to now. (In fact we are not aware of another paper than \cite{GG};
in \cite{GG}, a specific example of type $A_\infty$ is developed, 
in connection with a triangulated category studied by Holm and Jorgensen \cite{HoJo}.)

For our purposes in this paper, it is always possible to 
work with sufficiently large finite subseeds 
of the seed $(\bz^-,G^-)$, and replace $\AA$ by the genuine cluster subalgebras
attached to them.
On the other hand, statements become nicer if we allow ourselves to
formulate them in terms of the infinite rank cluster algebra $\AA$.

\subsubsection{Monomial change of variables}\label{subsect_chvar}

Let $\bY=\{Y_{i,r}\mid (i,r)\in W\}$ be a new set of indeterminates over $\Q$.
Let $\bY^-=\{Y_{i,r}\in \bY\mid (i,r)\in W^-\}$.
For $(i,r)\in W^-$, we perform the substitution
\begin{equation}\label{chvar}
z_{i,r} = \prod_{k\ge 0,\ r+kb_{ii}\le 0} Y_{i,\,r+kb_{ii}}. 
\end{equation}
Note that all variables in the right-hand side of (\ref{chvar}) belong to $\bY^-$.

\begin{example}
{\rm
If $G^-$ is as in Figure~\ref{Fig1}, we have
\[
\begin{array}{llll}
z_{\bdeux,0}= Y_{\bdeux,0},
&z_{\bdeux,-2}= Y_{\bdeux,-2}Y_{\bdeux,0},
&z_{\bdeux,-4}= Y_{\bdeux,-4}Y_{\bdeux,-2}Y_{\bdeux,0},
&z_{\bdeux,-6}= Y_{\bdeux,-6}Y_{\bdeux,-4}Y_{\bdeux,-2}Y_{\bdeux,0},
\\[3mm]
z_{\bun,-1} = Y_{\bun,-1},
&z_{\bun,-5} = Y_{\bun,-5}Y_{\bun,-1},
&z_{\bun,-9} = Y_{\bun,-9}Y_{\bun,-5}Y_{\bun,-1},
&z_{\bun,-13} = Y_{\bun,-13}Y_{\bun,-9}Y_{\bun,-5}Y_{\bun,-1},
\\[3mm]
z_{\bun,-3} = Y_{\bun,-3},
&z_{\bun,-7} = Y_{\bun,-7}Y_{\bun,-3},
&z_{\bun,-11} = Y_{\bun,-11}Y_{\bun,-7}Y_{\bun,-3},
&\emph{etc.}
\end{array}
\]
}
\end{example}

\subsubsection{Sequence of vertices}\label{subsect_seq_vert}

As explained in Example~\ref{example_quivers},
the arrows of $G^-$ of the form $(i, r) \longrightarrow (i, r+b_{ii})$
are called vertical and displayed in columns.
To each column we attach an initial \emph{label} given by the index of its
top vertex $(i, r)$, for which $r$ is maximal among the vertices of the
column. 

We now form a sequence of $tn$ columns by induction as follows.
At each step we pick a column whose label $(i, r)$ has maximal $r$
among labels of all columns.
After picking a column with label $(i,r)$, we change its
label to $(i, r-b_{ii})$. 
Finally, reading column after column in this ordering,
from top to bottom, we get an infinite sequence $\SS$ of vertices of $G^-$.

\begin{example}
{\rm
If $G^-$ is as in Figure~\ref{Fig0}, then $t=1$,
the sequence of columns consists of 3 columns, and we obtain the following
sequence of vertices:  
\[
\SS=((2,0), (2,-2), (2,-4), \ldots, (1, -1), (1,-5), (1, -9), \ldots, 
(3, -1), (3,-3), (3, -5), \ldots ).
\]
(Here, the columns labelled (1,-1) and (3,-1) could be interchanged.) 

If $G^-$ is as in Figure~\ref{Fig1}, then $t=2$,
the sequence of columns consists of 4 columns and gives the following
sequence of vertices:  
\[
\begin{array}{rcl}
\SS&=&((\bdeux,0), (\bdeux,-2), (\bdeux,-4), \ldots, (\bun, -1), (\bun,-5), (\bun, -9), \ldots, 
\\[3mm]
&&(\bdeux,0), (\bdeux,-2), (\bdeux,-4), \ldots, (\bun, -3), (\bun,-7), (\bun, -11), \ldots )
\end{array}
\]
Note that the column with vertices $(\bdeux,r)$ appears two times. It appears first because
its initial label is $(\bdeux,0)$. After picking it, its label is changed to $(\bdeux,-2)$, 
so it appears again between the columns labelled $(\bun,-1)$ and $(\bun,-3)$.  
}
\end{example}

Finally, for $(i,r)\in G^-$, we define $k_{i,r}$ to be the unique positive integer $k$ satisfying
\begin{equation}\label{def_k}
0< kb_{ii} - |r| \le b_{ii}. 
\end{equation}
In other words, $(i,r)$ is the $k$th vertex in its column, counting from the top.

\begin{example}
{\rm
If $G^-$ is as in Figure~\ref{Fig1}, then 
\[
k_{\bdeux,-2} = 2,\quad k_{\bun,-9} = 3. 
\]
}
\end{example}

\subsection{Quantum affine algebras}\label{subsect-qaff}

\subsubsection{The algebra $U_q(\hg)$}
Let $\g$ be the simple Lie algebra over $\C$ with Cartan matrix $C$.
We denote by $\a_i\ (i\in I)$ the simple roots of $\g$, and by 
$\varpi_i\ (i\in I)$ the fundamental weights.
They are related by
\begin{equation}
\a_i = \sum_{j\in I} c_{ji} \varpi_j. 
\end{equation}
Let $\dh$ be the dual Coxeter number of $\g$, see \cite[\S 6.1]{K}.
The values of $\dh$ are recalled in Table~\ref{table1}.
\begin{table}[t]%
\begin{center}
\begin{tabular}{c|c|c|c|c|c|c|c|c|c}
$\g$ & $A_n$ & $B_n$ & $C_n$ & $D_n$ & $E_6$ & $E_7$ & $E_8$ & $F_4$ & $G_2$\\
\hline
$t$ & 1 & 2 & 2 & 1 & 1 & 1 & 1 & 2 & 3 \\
$\dh$ & $n+1$ & $2n-1$ & $n+1$ & $2n-2$ & 12 & 18 & 30 & 9 & 4 
\end{tabular}
\end{center}
\bigskip
\caption{Dual Coxeter numbers}
\label{table1}
\end{table}

Let $\hg$ be the corresponding untwisted affine Lie algebra.
Thus if $\g$ has type $X_n$ in the Cartan-Killing classification,
$\hg$ has type $X_n^{(1)}$ in the Kac classification \cite[\S4.8]{K}.
Let $U_q(\hg)$ be the Drinfeld-Jimbo quantum enveloping algebra of $\hg$,
see \eg \cite{CP}.
We regard $U_q(\hg)$ as a $\C$-algebra with quantum parameter $q\in \C^*$
not a root of unity.

\subsubsection{$q$-characters}

Frenkel and Reshetikin \cite{FR} have attached to every complex finite-dimensional representation
of $U_q(\hg)$ a $q$-character $\chi_q(M)$. If $M$ is irreducible, it is determined up to isomorphism
by its $q$-character. The irreducible finite-dimensional representations of
$U_q(\hg)$ have been classified by Chari and Pressley in terms of Drinfeld polynomials,
see \cite[Theorem 12.2.6]{CP}.
Equivalently, irreducible finite-dimensional representations of $U_q(\hg)$
can be parametrized by the highest dominant monomial of their $q$-character \cite{FR},
and this is the parametrization we shall use.

By definition, the $q$-character $\chi_q(M)$ is a Laurent polynomial with positive integer
coefficients in the infinite set of variables 
$\YY=\{Y_{i,a}\mid i\in I,\ a\in \C^*\}$, which should be seen as a quantum affine analogue of 
$\{e^{\varpi_i}\mid i\in I\}$.
In this paper we will be concerned only with polynomials involving the
subset of variables
\[
 Y_{i, q^r},\qquad ((i,r) \in W).
\]
For simplicity of notation, we shall therefore write $Y_{i, r}$ instead of
$Y_{i, q^r}$.
Thus our $q$-characters will be Laurent polynomials in the variables
of the set $\bY$ introduced in \S\ref{subsect_chvar}.

Let $m$ be a \emph{dominant} monomial in the variables $Y_{i, r} \in \bY$,
that is, a monomial with nonnegative exponents. We denote by 
$L(m)$ the corresponding irreducible representation of $U_q(\hg)$,
and by $\chi_q(m)=\chi_q(L(m))$ its $q$-character.
For example, if $m$ is of the form
\[
m = \prod_{j=0}^{k-1}Y_{i,\,r+jb_{ii}},\qquad (i\in I,\ r\in\Z,\ k\ge 1),  
\]
$L(m)$ is called a \emph{Kirillov-Reshetikhin module}, and usually denoted
by $W^{(i)}_{k,r}$.
In particular, if $k=1$ we get a \emph{fundamental module}
$W^{(i)}_{1,r} = L(Y_{i,r})$.
By convention, if $k=0$ the module $W^{(i)}_{0,r}$ is the trivial
one-dimensional module for every $(i,r)$, and its $q$-character is
equal to 1.

Finally, following \cite{FR}, for $(i,r) \in V$ we introduce
the following quantum affine analogue of~$e^{\a_i}$:
\begin{equation}\label{def_A_var}
 A_{i,r} := Y_{i,r-d_i}Y_{i,r + d_i} 
\left(\prod_{j:\ c_{ji} = -1} Y_{j,r}\right)^{-1}
\left(\prod_{j:\ c_{ji} = -2} Y_{j,r-1}Y_{j,r+1}\right)^{-1}
\left(\prod_{j:\ c_{ji} = -3} Y_{j,r-2}Y_{j,r}Y_{j,r+2}\right)^{-1}
\end{equation}
Note that since $(i,r)\in V$, we have $(i,r\pm d_i)\in W$. If $c_{ji}<0$, we also have, because 
of (\ref{easy-prop}),
\[
(j,r+ c_{ji}+1) = (j,r+d_j(c_{ji}+1)) = (j,r+b_{ij}+d_j)\in W.
\]
It follows that $A_{i,r}$ is a Laurent monomial in the variables $Y_{j,s}$ with $(j,s)\in W$.

\section{An algorithm for the $q$-characters of Kirillov-Reshetikhin modules}\label{section2}

\subsection{Statement and examples}
Let $\AA$ be the cluster algebra defined in \S\ref{subsect_ca}, with initial seed
$\Si=(\bz^-, G^-)$, and
let 
\[
\SS = ((i_1,r_1), (i_2,r_2), (i_3,r_3), \ldots )
\] 
be the sequence of vertices of the quiver of $\AA$ defined in \S\ref{subsect_seq_vert}.
We denote by $\mu_\SS(\Si)$ the new seed obtained after performing the sequence of mutations
indexed by $\SS$, that is,
by mutating first at vertex $(i_1,r_1)$, then at vertex $(i_2,r_2)$, \emph{etc.} 
More generally, for $m\ge 1$, let $\Si_m = \mu_\SS^m(\Si)$ be the seed obtained from $\Si$ after $m$
repetitions of the mutation sequence $\mu_\SS$. 
Let $z^{(m)}_{i,r}$ be the cluster variable of 
$\Si_m$ sitting at vertex $(i,r)\in W^-$; this is a Laurent polynomial in the 
initial variables $z_{j,s},\ (j,s)\in W^-$. 
Let $y^{(m)}_{i,r}$ be the Laurent
polynomial obtained from $z^{(m)}_{i,r}$ by performing the change of variables (\ref{chvar})
of \S\ref{subsect_chvar}; this is a Laurent polynomial in the variables $Y_{j,s},\ (j,s)\in W^-$.

\begin{Thm}\label{thm1}
\begin{itemize}
 \item[(a)] The quiver of $\mu_\SS(\Si)$ is equal to the quiver of $\Si$, that is, to $G^-$.
 \item[(b)] 
Suppose that $m\ge \dh/2$. Then, the $y^{(m)}_{i,r}$ are the $q$-characters of the Kirillov-Reshetikhin modules.
More precisely, for $m\ge \dh/2$ there holds
\[
y^{(m)}_{i,r} = \chi_q\left(W^{(i)}_{k,\ r-2tm}\right). 
\]
where $k=k_{i,r}$ is defined as in \S\ref{subsect_seq_vert}.
\end{itemize}
 
\end{Thm}

\begin{remark}
{\rm 
It is well known that, for $p\in\Z$, the $q$-character $\chi_q(W^{(i)}_{k,r+p})$ is deduced
from $\chi_q(W^{(i)}_{k,r})$ by applying the ring automorphism mapping $Y_{j,s}$
to $Y_{j,s+p}$ for every $(j,s)\in I\times \Z$. Therefore, modulo these straightforward 
automorphisms, Theorem~\ref{thm1} describes the $q$-characters of \emph{all} Kirillov-Reshetikhin
modules. 
} 
\end{remark}

\begin{remark}
{\rm
Although the statement of Theorem~\ref{thm1} involves an infinite seed $\Si$
and an infinite sequence of mutations $\SS$, the calculation of 
the $q$-character of a given Kirillov-Reshetikhin module requires only
a finite number of mutations on a finite initial
segment of the semi-infinite quiver. More precisely,
the proof of Theorem~\ref{thm1} will show that all the
$q$-characters $\chi_q(W^{(i)}_{k,s})$ with $k=1,\ldots,l$ can be calculated
using $(h'+2l-1)h'n/2$ mutations, where $h' = \lceil\dh/2\rceil$.
} 
\end{remark}

\begin{example}
{\rm
Let $\g$ be of type $A_3$. The quiver $G^-$ of the initial
seed is displayed in Figure~\ref{Fig0}. The initial cluster variables are
\[
\begin{array}{llll}
z_{2,0}= Y_{2,0},
&z_{2,-2}= Y_{2,-2}Y_{2,0},
&z_{2,-4}= Y_{2,-4}Y_{2,-2}Y_{2,0},
&\emph{etc.}
\\[3mm]
z_{1,-1} = Y_{1,-1},
&z_{1,-3} = Y_{1,-3}Y_{1,-1},
&z_{1,-5} = Y_{1,-5}Y_{1,-3}Y_{1,-1},
&\emph{etc.}
\\[3mm]
z_{3,-1} = Y_{3,-1},
&z_{3,-3} = Y_{3,-3}Y_{3,-1},
&z_{3,-5} = Y_{3,-5}Y_{3,-3}Y_{3,-1},
&\emph{etc.}
\end{array}
\]
After the mutation sequence $\mu_\SS$, the first new cluster variables are
\[
\begin{array}{rcl}
y^{(1)}_{2,0}&=& Y_{2,-2}+Y_{1,-1}Y_{3,-1}Y_{2,0}^{-1},\\[1mm]
y^{(1)}_{2,-2}&=& Y_{2,-4}Y_{2,-2}+Y_{1,-1}Y_{3,-1}Y_{2,-4}Y_{2,0}^{-1}
+Y_{1,-3}Y_{1,-1}Y_{2,-2}^{-1}Y_{2,0}^{-1}Y_{3,-3}Y_{3,-1},\\[1mm]
y^{(1)}_{2,-4}&=& Y_{2,-6}Y_{2,-4}Y_{2,-2} + Y_{1,-1}Y_{3,-1}Y_{2,-6}Y_{2,-4}Y_{2,0}^{-1}
+Y_{1,-3}Y_{1,-1}Y_{2,-6}Y_{2,-2}^{-1}Y_{2,0}^{-1}Y_{3,-3}Y_{3,-1},\\
&&+\ Y_{1,-5}Y_{1,-3}Y_{1,-1}Y_{2,-4}^{-1}Y_{2,-2}^{-1}Y_{2,0}^{-1}Y_{3,-5}Y_{3,-3}Y_{3,-1},
\\[1mm]
y^{(1)}_{1,-1} &=& Y_{1,-3}+Y_{1,-1}^{-1}Y_{2,-2}+Y_{2,0}^{-1}Y_{3,-1},\\[1mm]
y^{(1)}_{1,-3} &=& Y_{1,-5}Y_{1,-3}+Y_{1,-5}Y_{1,-1}^{-1}Y_{2,-2}+Y_{1,-5}Y_{2,0}^{-1}Y_{3,-1}
+Y_{1,-3}^{-1}Y_{1,-1}^{-1}Y_{2,-4}Y_{2,-2}\\
&&+\ Y_{1,-3}^{-1}Y_{2,-4}Y_{2,0}^{-1}Y_{3,-1} + Y_{2,-2}^{-1}Y_{2,0}^{-1}Y_{3,-3}Y_{3,-1},\\[1mm]
y^{(1)}_{1,-5} &=& Y_{1,-7}Y_{1,-5}Y_{1,-3} + Y_{1,-7}Y_{1,-5}Y_{1,-1}^{-1}Y_{2,-2}
+Y_{1,-7}Y_{1,-5}Y_{2,0}^{-1}Y_{3,-1} + Y_{1,-7}Y_{1,-3}^{-1}Y_{1,-1}^{-1}Y_{2,-4}Y_{2,-2}\\
&&+\ Y_{1,-7}Y_{1,-3}^{-1}Y_{2,-4}Y_{2,0}^{-1}Y_{3,-1} + Y_{1,-7}Y_{2,-2}^{-1}Y_{2,0}^{-1}Y_{3,-3}Y_{3,-1}
+ Y_{1,-5}^{-1}Y_{1,-3}^{-1}Y_{1,-1}^{-1}Y_{2,-6}Y_{2,-4}Y_{2,-2}\\
&&+\ Y_{1,-5}^{-1}Y_{1,-3}^{-1}Y_{2,-6}Y_{2,-4}Y_{2,0}^{-1}Y_{3,-1} + 
Y_{1,-5}^{-1}Y_{2,-6}Y_{2,-2}^{-1}Y_{2,0}^{-1}Y_{3,-3}Y_{3,-1}\\
&&+\ Y_{2,-4}^{-1}Y_{2,-2}^{-1}Y_{2,0}^{-1}Y_{3,-5}Y_{3,-3}Y_{3,-1}, 
\end{array}
\]
(We omit the variables $y^{(1)}_{3,-1}$, $y^{(1)}_{3,-5}$, $y^{(1)}_{3,-5}$, since they are
readily obtained from $y^{(1)}_{1,-1}$, $y^{(1)}_{1,-5}$, $y^{(1)}_{1,-5}$ via the symmetry 
($1\leftrightarrow 3$).)
After a second application of the mutation sequence $\mu_\SS$, the first new cluster variables are
\[
\begin{array}{rcl}
y^{(2)}_{2,0}&=& Y_{2,-4}+Y_{1,-3}Y_{3,-3}Y_{2,-2}^{-1}+Y_{1,-3}Y_{3,-1}^{-1}+Y_{1,-1}^{-1}Y_{3,-3}
+Y_{1,-1}^{-1}Y_{2,-2}Y_{3,-1}^{-1}+Y_{2,0}^{-1},\\[1mm]
y^{(2)}_{2,-2}&=& Y_{2,-6}Y_{2,-4}+Y_{1,-3}Y_{3,-3}Y_{2,-6}Y_{2,-2}^{-1}
+Y_{1,-5}Y_{1,-3}Y_{2,-4}^{-1}Y_{2,-2}^{-1}Y_{3,-5}Y_{3,-3}+Y_{1,-5}Y_{2,0}^{-1}Y_{3,-3}^{-1}\\
&&+\ Y_{1,-3}^{-1}Y_{2,0}^{-1}Y_{3,-5} + Y_{1,-3}^{-1}Y_{2,-4}Y_{2,0}^{-1}Y_{3,-3}^{-1}
+Y_{2,-6}Y_{2,0}^{-1} + Y_{1,-5}Y_{2,-4}^{-1}Y_{2,0}^{-1}Y_{3,-5} + Y_{1,-3}Y_{2,-6}Y_{3,-1}^{-1}\\
&&+\ Y_{1,-5}Y_{1,-3}Y_{3,-3}^{-1}Y_{3,-1}^{-1} + Y_{1,-5}Y_{1,-3}Y_{2,-4}^{-1}Y_{3,-5}Y_{3,-1}^{-1}
+ Y_{1,-1}^{-1}Y_{2,-6}Y_{3,-3} + Y_{1,-3}^{-1}Y_{1,-1}^{-1}Y_{3,-5}Y_{3,-3} \\
&&+\ Y_{1,-5}Y_{1,-1}^{-1}Y_{2,-4}^{-1}Y_{3,-5}Y_{3,-3} 
+ Y_{1,-3}^{-1}Y_{1,-1}^{-1}Y_{2,-4}Y_{2,-2}Y_{3,-3}^{-1}Y_{3,-1}^{-1}
+ Y_{1,-1}^{-1}Y_{2,-6}Y_{2,-2}Y_{3,-1}^{-1}\\
&&+\ Y_{1,-3}^{-1}Y_{1,-1}^{-1}Y_{2,-2}Y_{3,-5}Y_{3,-1}^{-1} 
+ Y_{1,-5}Y_{1,-1}^{-1}Y_{2,-2}Y_{3,-3}^{-1}Y_{3,-1}^{-1} 
+ Y_{1,-5}Y_{1,-1}^{-1}Y_{2,-4}^{-1}Y_{2,-2}Y_{3,-5}Y_{3,-1}^{-1}\\
&&\ + Y_{2,-2}^{-1}Y_{2,0}^{-1}, 
\\[1mm]
y^{(2)}_{1,-1} &=& Y_{1,-5}+Y_{1,-3}^{-1}Y_{2,-4}+Y_{2,-2}^{-1}Y_{3,-3}+Y_{3,-1}^{-1},\\[1mm]
y^{(2)}_{1,-3} &=& Y_{1,-7}Y_{1,-5}+Y_{1,-7}Y_{1,-3}^{-1}Y_{2,-4}+Y_{1,-7}Y_{2,-2}^{-1}Y_{3,-3}
+Y_{1,-5}^{-1}Y_{1,-3}^{-1}Y_{2,-6}Y_{2,-4}\\
&&+\ Y_{1,-5}^{-1}Y_{2,-6}Y_{2,-2}^{-1}Y_{3,-3} + Y_{2,-4}^{-1}Y_{2,-2}^{-1}Y_{3,-5}Y_{3,-3}
+Y_{1,-5}^{-1}Y_{2,-6}Y_{3,-1}^{-1} + Y_{1,-7}Y_{3,-1}^{-1}\\
&&+\ Y_{2,-4}^{-1}Y_{3,-5}Y_{3,-1}^{-1} + Y_{3,-3}^{-1}Y_{3,-1}^{-1},
\end{array}
\]
Here $\dh/2=2$, so we can observe that the cluster variables obtained after performing 
2 times the mutation sequence $\mu_\SS$ are indeed $q$-characters of Kirillov-Reshetikhin
modules, namely,
\[
\begin{array}{lll}
y^{(2)}_{2,0} = \chi_q(Y_{2,-4}),\quad&
y^{(2)}_{2,-2} = \chi_q(Y_{2,-6}Y_{2,-4}),\quad&
\emph{etc.}\\[1mm]
y^{(2)}_{1,-1} = \chi_q(Y_{1,-5}),\quad&
y^{(2)}_{1,-3} = \chi_q(Y_{1,-7}Y_{1,-5}),\quad&
\emph{etc.}\\[1mm]
y^{(2)}_{3,-1} = \chi_q(Y_{3,-5}),\quad& 
y^{(2)}_{3,-3} = \chi_q(Y_{3,-7}Y_{3,-5}),\quad&
\emph{etc.} 
\end{array}
\]
} 
\end{example}

\begin{example}
{\rm
Let $\g$ be of type $B_2$. The quiver $G^-$ of the initial
seed is displayed in Figure~\ref{Fig1}. The initial cluster variables are
\[
\begin{array}{llll}
z_{\bdeux,0}= Y_{\bdeux,0},
&z_{\bdeux,-2}= Y_{\bdeux,-2}Y_{\bdeux,0},
&z_{\bdeux,-4}= Y_{\bdeux,-4}Y_{\bdeux,-2}Y_{\bdeux,0},
&\emph{etc.}
\\[3mm]
z_{\bun,-1} = Y_{\bun,-1},
&z_{\bun,-5} = Y_{\bun,-5}Y_{\bun,-1},
&z_{\bun,-9} = Y_{\bun,-9}Y_{\bun,-5}Y_{\bun,-1},
&\emph{etc.}
\\[3mm]
z_{\bun,-3} = Y_{\bun,-3},
&z_{\bun,-7} = Y_{\bun,-7}Y_{\bun,-3},
&z_{\bun,-11} = Y_{\bun,-11}Y_{\bun,-7}Y_{\bun,-3},
&\emph{etc.}
\end{array}
\]
After the mutation sequence $\mu_\SS$, the first new cluster variables are
\[
\begin{array}{rcl}
y^{(1)}_{\bdeux,0}&= &Y_{\bdeux,-4} + Y_{\bun,-3}Y_{\bdeux,-2}^{-1}, \\[1mm]
y^{(1)}_{\bdeux,-2}&= &Y_{\bdeux,-6}Y_{\bdeux,-4} + Y_{\bun,-3}Y_{\bdeux,-6}Y_{\bdeux,-2}^{-1}
+ Y_{\bun,-5}Y_{\bun,-3}Y_{\bdeux,-4}^{-1}Y_{\bdeux,-2}^{-1}
+ Y_{\bun,-3}Y_{\bun,-1}^{-1}, \\[1mm]
y^{(1)}_{\bdeux,-4}&= &Y_{\bdeux,-8}Y_{\bdeux,-6}Y_{\bdeux,-4} + Y_{\bun,-3}Y_{\bdeux,-8}Y_{\bdeux,-6}Y_{\bdeux,-2}^{-1}
+ Y_{\bun,-5}Y_{\bun,-3}Y_{\bdeux,-8}Y_{\bdeux,-4}^{-1}Y_{\bdeux,-2}^{-1}\\
&&+\ Y_{\bun,-7}Y_{\bun,-5}Y_{\bun,-3}Y_{\bdeux,-6}^{-1}Y_{\bdeux,-4}^{-1}Y_{\bdeux,-2}^{-1}
+ Y_{\bun,-7}Y_{\bun,-3}Y_{\bun,-1}^{-1}Y_{\bdeux,-6}^{-1}
+Y_{\bun,-1}^{-1}Y_{\bun,-3}Y_{\bdeux,-8}
,\\[1mm]
y^{(1)}_{\bun,-1}&= &Y_{\bun,-5} + Y_{\bun,-1}^{-1}Y_{\bdeux,-4}Y_{\bdeux,-2}
+Y_{\bdeux,-4}Y_{\bdeux,0}^{-1}+Y_{\bun,-3}Y_{\bdeux,-2}^{-1}Y_{\bdeux,0}^{-1},\\[1mm]
y^{(1)}_{\bun,-5}&= & Y_{\bun,-5}Y_{\bun,-9} + Y_{\bun,-9}Y_{\bun,-1}^{-1}Y_{\bdeux,-4}Y_{\bdeux,-2}
+ Y_{\bun,-9}Y_{\bdeux,-4}Y_{\bdeux,0}^{-1} + Y_{\bun,-9}Y_{\bun,-3}Y_{\bdeux,-2}^{-1}Y_{\bdeux,0}^{-1}\\
&&+\ Y_{\bun,-5}^{-1}Y_{\bun,-1}^{-1}Y_{\bdeux,-8}Y_{\bdeux,-6}Y_{\bdeux,-4}Y_{\bdeux,-2}
+ Y_{\bun,-5}^{-1}Y_{\bdeux,-8}Y_{\bdeux,-6}Y_{\bdeux,-4}Y_{\bdeux,0}^{-1}\\
&&+\ Y_{\bun,-5}^{-1}Y_{\bun,-3}Y_{\bdeux,-8}Y_{\bdeux,-6}Y_{\bdeux,-2}^{-1}Y_{\bdeux,0}^{-1}
+ Y_{\bun,-3}Y_{\bdeux,-8}Y_{\bdeux,-4}^{-1}Y_{\bdeux,-2}^{-1}Y_{\bdeux,0}^{-1}\\
&&+\ Y_{\bun,-7}Y_{\bun,-3}Y_{\bdeux,-6}^{-1}Y_{\bdeux,-4}^{-1}Y_{\bdeux,-2}^{-1}Y_{\bdeux,0}^{-1},
\\[1mm]
y^{(1)}_{\bun,-3}&= &Y_{\bun,-7} + Y_{\bun,-3}^{-1}Y_{\bdeux,-6}Y_{\bdeux,-4}
+Y_{\bdeux,-6}Y_{\bdeux,-2}^{-1}+Y_{\bun,-5}Y_{\bdeux,-4}^{-1}Y_{\bdeux,-2}^{-1}
+Y_{\bun,-1}^{-1},\\[1mm]
y^{(1)}_{\bun,-7}&= & Y_{\bun,-7}Y_{\bun,-11} + Y_{\bun,-11}Y_{\bun,-3}^{-1}Y_{\bdeux,-6}Y_{\bdeux,-4}
+ Y_{\bun,-11}Y_{\bdeux,-6}Y_{\bdeux,-2}^{-1} + Y_{\bun,-11}Y_{\bun,-5}Y_{\bdeux,-4}^{-1}Y_{\bdeux,-2}^{-1}\\
&&+\ Y_{\bun,-7}^{-1}Y_{\bun,-3}^{-1}Y_{\bdeux,-10}Y_{\bdeux,-8}Y_{\bdeux,-6}Y_{\bdeux,-4}
+ Y_{\bun,-7}^{-1}Y_{\bdeux,-10}Y_{\bdeux,-8}Y_{\bdeux,-6}Y_{\bdeux,-2}^{-1}\\
&&+\ Y_{\bun,-7}^{-1}Y_{\bun,-5}Y_{\bdeux,-10}Y_{\bdeux,-8}Y_{\bdeux,-4}^{-1}Y_{\bdeux,-2}^{-1}
+ Y_{\bun,-5}Y_{\bdeux,-10}Y_{\bdeux,-6}^{-1}Y_{\bdeux,-4}^{-1}Y_{\bdeux,-2}^{-1}\\
&&+\ Y_{\bun,-9}Y_{\bun,-5}Y_{\bdeux,-8}^{-1}Y_{\bdeux,-6}^{-1}Y_{\bdeux,-4}^{-1}Y_{\bdeux,-2}^{-1}
+ Y_{\bun,-9}Y_{\bun,-1}^{-1}Y_{\bdeux,-8}^{-1}Y_{\bdeux,-6}^{-1}+ Y_{\bun,-11}Y_{\bun,-1}^{-1}\\
&&+\ Y_{\bun,-1}^{-1}Y_{\bdeux,-10}Y_{\bdeux,-6}^{-1} + 
Y_{\bun,-7}^{-1}Y_{\bun,-1}^{-1}Y_{\bdeux,-10}Y_{\bdeux,-8}+Y_{\bun,-5}^{-1}Y_{\bun,-1}^{-1}.
\end{array}
\]
Here $\dh/2=3/2$, and we can observe that certain cluster variables are not yet
$q$-characters of Kirillov-Reshetikhin modules.
But some already are, namely
\[
\begin{array}{lll}
y^{(1)}_{\bun,-3}= \chi_q(Y_{\bun,-7}),\quad&
y^{(1)}_{\bun,-7}= \chi_q(Y_{\bun,-7}Y_{\bun,-11}),\quad&
\emph{etc.}
\end{array}
\]
After a second application of the mutation sequence $\mu_\SS$, since $2>3/2$,
all the new cluster variables are $q$-characters of Kirillov-Reshetikhin modules.
For example
\[
\begin{array}{rcl}
y^{(2)}_{\bdeux,0}&= &Y_{\bdeux,-8} + Y_{\bun,-7}Y_{\bdeux,-6}^{-1} + Y_{\bun,-3}^{-1}Y_{\bdeux,-4} + Y_{\bdeux,-2}^{-1}
= \chi_q(Y_{\bdeux,-8}).
\end{array}
\]
}
\end{example}

\subsection{Proof of Theorem~\ref{thm1}}

The proof relies on two main ingredients which we shall first review, namely, 
the $T$-systems, and the truncated $q$-characters.

\subsubsection{$T$-systems}\label{ssect-T-syst} 

With the quantum affine algebra $U_q(\hg)$ is associated a system
of difference equations called a $T$-system \cite{KNS1}. 
Its unknowns are denoted by
\[
T^{(i)}_{k, r},\qquad (i \in I,\ k\in \N,\ r\in\Z). 
\]
We fix the initial boundary condition 
\begin{equation}
T^{(i)}_{0, r} = 1,\qquad (i\in I,\ r\in\Z).
\end{equation}
If $\g$ is of type $A_n, D_n, E_n$, the $T$-system
equations are
\begin{equation}
T^{(i)}_{k,r+1}T^{(i)}_{k,r-1} =  T^{(i)}_{k-1,r+1}T^{(i)}_{k+1,r-1} + 
\prod_{j :\ c_{ij}=-1} T^{(j)}_{k,r},\qquad (i \in I,\ k \ge 1,\ r\in\Z). 
\end{equation}
If $\g$ is not of simply laced type, the $T$-system equations are more complicated. 
They can be written in the form
\begin{equation}\label{Tsystem_general}
T^{(i)}_{k,r+d_i}T^{(i)}_{k,r-d_i} =  T^{(i)}_{k-1,r+d_i}T^{(i)}_{k+1,r-d_i} + 
S^{(i)}_{k,r},\qquad (i \in I,\ k \ge 1,\ r\in\Z), 
\end{equation}
where $S^{(i)}_{k,r}$ is defined as follows. 
If $d_i\ge 2$ then
\begin{equation}
 S^{(i)}_{k,r} = \prod_{j:\ c_{ji}=-1} T^{(j)}_{k,r} 
\prod_{j:\ c_{ji}\le -2} T^{(j)}_{d_ik,\ r-d_i+1}. 
\end{equation}
If $d_i=1$ and $t=2$, then
\begin{equation}
S^{(i)}_{k,r} = 
\left\{
\begin{array}{ll}
\ds\prod_{j:\ c_{ij}=-1} T^{(j)}_{k,r}
\prod_{j:\ c_{ij}= -2} T^{(j)}_{l,r} T^{(j)}_{l,r+2},
& \mbox{if $k=2l$,}\\[5mm]
\ds\prod_{j:\ c_{ij}=-1} T^{(j)}_{k,r}
\prod_{j:\ c_{ij}= -2} T^{(j)}_{l+1,r} T^{(j)}_{l,r+2}
& \mbox{if $k=2l+1$.}
\end{array}
\right.
\end{equation}
Finally, if $d_i=1$ and $t=3$, that is, if $\g$ is of type $G_2$,
denoting by $j$ the other vertex of $\de$ we have $d_j=3$ and
\begin{equation}
S^{(i)}_{k,r} = 
\left\{
\begin{array}{ll}
T^{(j)}_{l,r} T^{(j)}_{l,r+2} T^{(j)}_{l,r+4} & \mbox{if $k=3l$,}\\[2mm]
T^{(j)}_{l+1,r} T^{(j)}_{l,r+2} T^{(j)}_{l,r+4} & \mbox{if $k=3l+1$,}\\[2mm]
T^{(j)}_{l+1,r} T^{(j)}_{l+1,r+2} T^{(j)}_{l,r+4} & \mbox{if $k=3l+2$.}
\end{array}
\right.
\end{equation}

\begin{example}\label{example_Tsystem}
{\rm
Let $\g$ be of type $B_2$. 
The Cartan matrix is 
\[
 C =
\pmatrix{2 & -1\cr
-2 & 2}
\]
and we have $d_1=2$ and $d_2=1$.
The $T$-system reads:
\[
\begin{array}{lcll}
T^{(\bun)}_{k,r+2}T^{(\bun)}_{k,r-2} &=&  T^{(\bun)}_{k-1,r+2}T^{(\bun)}_{k+1,r-2} + 
T^{(\bdeux)}_{2k,r-1},& (k \ge 1,\ r\in\Z),\\[2mm] 
T^{(\bdeux)}_{2l,r+1}T^{(\bdeux)}_{2l,r-1} &=&  T^{(\bdeux)}_{2l-1,r+1}T^{(\bdeux)}_{2l+1,r-1} + 
T^{(\bun)}_{l,r} T^{(\bun)}_{l,r+2}, &(l \ge 1,\ r\in\Z), \\[2mm]
T^{(\bdeux)}_{2l+1,r+1}T^{(\bdeux)}_{2l+1,r-1} &=&  T^{(\bdeux)}_{2l,r+1}T^{(\bdeux)}_{2l+2,r-1} + 
T^{(\bun)}_{l+1,r} T^{(\bun)}_{l,r+2}, &(l \ge 0,\ r\in\Z). 
\end{array}
\]

} 
\end{example}

It was conjectured in \cite{KNS1}, and proved in \cite{N2} (for $\g$ of type $A,D,E$) and \cite{H} (general case),
that the $q$-characters of the Kirillov-Reshetikhin modules of $U_q(\hg)$ satisfy the 
corresponding $T$-system. More precisely, we have
\begin{Thm}[\cite{N2}\cite{H}]\label{thm_Tsystem}
For $i \in I,\ k\in \N,\ r\in\Z$,
\[
T^{(i)}_{k,r} = \chi_q\left(W^{(i)}_{k,r}\right),
\]
is a solution of the $T$-system in the ring $\Z\left[Y_{i,r}^{\pm1}\mid (i,r) \in I\times\Z\right]$.
\end{Thm}

\subsubsection{Truncated $q$-characters}\label{ssect-truncated}

Let $\CC^-$ be the full subcategory of the category of finite-dimensional
$U_q(\hg)$-modules whose objects have all their composition factors of
the form $L(m)$ where $m$ is a dominant monomial in the variables of $\bY^-$.

\begin{Lem} \label{lemCC^-}
The $q$-character of an object in $\CC^-$ belongs to $\Z\left[Y_{i,r}^{\pm1}\mid Y_{i,r}\in \bY\right]$.
\end{Lem}

\proof
A simple object of $\CC^-$ is a quotient of a tensor product of fundamental representations of $\CC^-$.
But the $q$-character of a fundamental representation can be calculated by means of the Frenkel-Mukhin 
algorithm \cite{FM}. At each step the algorithm produces monomials which involve only variables 
$Y_{i,r}\in \bY$. Hence the result.
\cqfd

Note that for a dominant monomial $m$ in the variables of $\bY^-$, 
the $q$-character $\chi_q(m)$
may contain Laurent monomials $m'$ involving variables $Y_{i,r}\in\bY\setminus\bY^-$. 
Following \cite{HL},
we define the \emph{truncated $q$-character}
$\chi_q^-(m)$ to be the Laurent polynomial obtained from 
$\chi_q(m)$ by discarding all these monomials $m'$. 
So, by definition, $\chi_q^-(m) \in \Z\left[Y_{i,r}^{\pm1}\mid Y_{i,r} \in \bY^-\right]$.

\begin{example}
{\rm 
Let $\g$ be of type $B_2$. We keep the notation of Example~\ref{example_Tsystem}.
The fundamental modules $L(Y_{\bun,-3})$ and $L(Y_{\bdeux,-4})$ have $q$-characters 
equal to
\[
\begin{array}{lcl}
\chi_q(Y_{\bun,-3}) &= & Y_{\bun,-3} + Y_{\bun,1}^{-1}Y_{\bdeux,-2}Y_{\bdeux,0}
+Y_{\bdeux,-2}Y_{\bdeux,2}^{-1}+Y_{\bun,-1}Y_{\bdeux,0}^{-1}Y_{\bdeux,2}^{-1}
+Y_{\bun,3}^{-1}, \\[1mm]
\chi_q(Y_{\bdeux,-4}) &= &Y_{\bdeux,-4} + Y_{\bun,-3}Y_{\bdeux,-2}^{-1}
 + Y_{\bun,1}^{-1}Y_{\bdeux,0} 
+ Y_{\bdeux,2}^{-1}.
\end{array}
\]
The corresponding truncated $q$-characters are
\[
\begin{array}{lcl}
\chi^-_q(Y_{\bun,-3}) &= & Y_{\bun,-3}, \\[1mm]
\chi^-_q(Y_{\bdeux,-4}) &= &Y_{\bdeux,-4} + Y_{\bun,-3}Y_{\bdeux,-2}^{-1}.
\end{array}
\]
} 
\end{example}

\begin{Prop}\label{prop_truncated}
\begin{itemize}
\item[(i)]$\CC^-$ is a tensor category. 
\item[(ii)] The assignment $[L(m)]\mapsto \chi_q^-(m)$ extends to an injective
ring homomorphism from the Grothendieck ring $K_0(\CC^-)$ to  $\Z\left[Y_{i,r}^{\pm1}\mid Y_{i,r}\in \bY^-\right]$.
\end{itemize}
\end{Prop}
\proof
The argument follows the same lines as \cite[\S5.2.4, \S6.2]{HL}.
Recall the Laurent monomials $A_{i,r}$ introduced in (\ref{def_A_var}).
By \cite{FR}, a Laurent monomial $m'$ of the $q$-character
of a simple object of $\CC^-$ can always be written in the form $m' = mM$ where $m$
is a dominant monomial in the variables of $\bY^-$,
and $M$ is a monomial in the variables $A_{i,k}^{-1}$ with $(i, k+d_i)\in W$.
Note that the $Y$-variable appearing in $A_{i,r}$ with the highest
spectral parameter is $Y_{i,r+d_i}$. It follows that $A_{i,r}^{-1}$
is a \emph{right-negative} monomial in the sense of \cite{FM}, that is,
the $Y$-variable with highest spectral parameter occuring in $A_{i,r}^{-1}$ has 
a negative exponent. 

Let $L(m)$ and $L(m')$ be simple objects of $\CC^-$, that is, $m$ and $m'$ are 
dominant monomials in the variables of $\bY^-$. If $L(m'')$ is a composition 
factor of $L(m)\otimes L(m')$, then $m''$ is a product of monomials of $\chi_q(m)$
and $\chi_q(m')$. So, we have $m''=mm'M$ where $M$ is a monomial
in the variables $A_{i,r}^{-1}$.
We claim that, since $m''$ is dominant, the spectral parameters $r$ have to 
satisfy $r+d_i\le 0$. Indeed otherwise $m''$ would be right-negative. 
Therefore, using Lemma~\ref{lemCC^-}, the monomial $m''$ contains only variables of~$\bY^-$, 
hence $L(m'')$ is in $\CC^-$, and $\CC^-$ is a monoidal category.
Moreover, by \cite[Prop. 5.1]{CP2}, the category
$\CC^-$ is stable by duals, so it is a tensor category. This proves (i).

To prove (ii) consider now an arbitrary Laurent monomial $m'$ of the $q$-character
of an object of $\CC^-$. As above, it can be written in the form $m' = mM$ where $m$
is a dominant monomial in the variables of $\bY^-$,
and $M$ is a monomial in the variables $A_{i,k}^{-1}$ with $(i, k+d_i)\in W$.
Now $m'$ contains a variable $Y_{j,s}\not\in\bY^-$ 
if and only if $M$ contains a negative power of $A_{i,r}$
for some pair $(i,r)$ such that $(i,r+d_i)\not\in W^-$.
So, if $R$ denotes the subring of $\Z\left[Y_{i,r}^{\pm1}\mid Y_{i,r}\in \bY\right]$
generated by all the monomials of the $q$-characters of the objects of $\CC^-$,
and if $I$ denotes the linear span of those monomials containing 
a variable $Y_{j,s}\in\bY\setminus\bY^-$, we see that $I$
is an ideal of $R$. Hence, if $\pi\colon R \to R/I$ is the natural projection,
we can realize the truncated $q$-character map $\chi_q^-$
as 
\[
\chi_q^- = \pi\circ\chi_q,
\] 
which shows that $\chi_q^-$ is 
a ring homomorphism $K_0(\CC^-) \to \Z\left[Y_{i,r}^{\pm1}\mid Y_{i,r}\in \bY^-\right]$.
Finally, the fact that $\chi_q^-$ is injective follows from the
fact that $I$ contains only non-dominant monomials, and that
two $q$-characters having the same dominant monomials with the same coefficients
are equal.
\cqfd

\subsubsection{Proof of the theorem}

We first notice that the initial cluster variables $z_{i,r}$ are
equal, after the change of variables (\ref{chvar}), to the truncated
$q$-characters of certain Kirillov-Reshetikhin modules, namely,
\[
z_{i,r} = \prod_{k\ge 0,\ r+kb_{ii}\le 0} Y_{i,r+kb_{ii}} = \chi_q^-\left(W^{(i)}_{k_{i,r},r}\right),
\]
where $k_{i,r}$ is defined as in (\ref{def_k}).
Indeed, the level of truncation is chosen so that
after truncation only the highest dominant monomial of these $q$-characters 
survives.

Now, the main idea of the proof is that the quiver $G^-$ and the mutation sequence 
$\mu_\SS$ are designed in such a way
that, at every step of the mutation sequence, the exchange relation is nothing else
than a $T$-system equation. 
Let us first check this when $\g$ is of rank two.

For $\g$ of type $A_2$,  
the sequence of mutated quivers obtained at each step of $\mu_\SS$
is shown in Appendix \S\ref{appendixA3}.
The mutations take place at the boxed vertices. 
Reading the second quiver of \S\ref{appendixA3}, we see
that the new cluster variable obtained after the first mutation is equal to
\[
\frac{z_{2,-2} + z_{1,-1}}{z_{2,0}}
=\frac{\chi_q^-\left(W^{(2)}_{2,-2}\right)
+\chi_q^-\left(W^{(1)}_{1,-1}\right)}{\chi_q^-\left(W^{(2)}_{1,0}\right)} 
=
\chi_q^-\left(W^{(2)}_{1,-2}\right).
\]
Here we have used Theorem~\ref{thm_Tsystem} and Proposition~\ref{prop_truncated}.
Similarly, reading the third quiver of \S\ref{appendixA3}, 
the new cluster variable obtained after the second mutation is equal to
\[
\frac{\chi_q^-\left(W^{(2)}_{3,-4}\right)\chi_q^-\left(W^{(2)}_{1,-2}\right)
+\chi_q^-\left(W^{(1)}_{2,-3}\right)}{\chi_q^-\left(W^{(2)}_{2,-2}\right)} 
=
\chi_q^-\left(W^{(2)}_{2,-4}\right).
\] 
An easy induction shows that, after every vertex of the second column
has been mutated, each cluster variable of the form $\chi_q^-\left(W^{(2)}_{k,-2k+2}\right)$
has been replaced by the new cluster variable $\chi_q^-\left(W^{(2)}_{k,-2k}\right)$.
We now continue by mutating vertices of the first column. We first get, at the top vertex
\[
\frac{\chi_q^-\left(W^{(1)}_{2,-3}\right)+\chi_q^-\left(W^{(2)}_{1,-2}\right)}{\chi_q^-\left(W^{(1)}_{1,-1}\right)} 
=
\chi_q^-\left(W^{(1)}_{1,-3}\right).
\] 
Then, mutating at the next vertex gives
\[
\frac{\chi_q^-\left(W^{(1)}_{3,-5}\right)\chi_q^-\left(W^{(1)}_{1,-3}\right)
+\chi_q^-\left(W^{(2)}_{2,-4}\right)}{\chi_q^-\left(W^{(1)}_{2,-3}\right)} 
=
\chi_q^-\left(W^{(1)}_{2,-5}\right).
\] 
By induction one sees that, after every vertex of the first column
has been mutated, each cluster variable of the form $\chi_q^-\left(W^{(1)}_{k,-2k+1}\right)$
has been replaced by a new cluster variable $\chi_q^-\left(W^{(1)}_{k,-2k-1}\right)$.
Moreover, one sees that the new quiver obtained after $\mu_\SS$ is nothing
else than $G^-$. Hence we conclude that one application of $\mu_\SS$
produces a seed with the same quiver, and in which every cluster variable 
$\chi_q^-\left(W^{(i)}_{k,r}\right)$ has been replaced
by $\chi_q^-\left(W^{(i)}_{k,r-2}\right)$. In other words, the effect of
$\mu_\SS$ is merely a uniform shift of the spectral parameters $r$ by $-2$.

The argument is similar for $\g$ of type $B_2$.  
The sequence of mutated quivers obtained at each step of $\mu_\SS$
is displayed in Appendix \S\ref{appendixB2}. 
Reading the second quiver of \S\ref{appendixB2}, we see
that the new cluster variable obtained after the first mutation is equal to
\[
\frac{z_{\bdeux,-2} + z_{\bun,-1}}{z_{\bdeux,0}}
=\frac{\chi_q^-\left(W^{(\bdeux)}_{2,-2}\right)
+\chi_q^-\left(W^{(\bun)}_{1,-1}\right)}{\chi_q^-\left(W^{(\bdeux)}_{1,0}\right)} 
=
\chi_q^-\left(W^{(\bdeux)}_{1,-2}\right).
\]
Similarly, reading the third quiver of \S\ref{appendixB2}, the new cluster variable 
obtained after the second mutation is equal to
\[
 \frac{\chi_q^-\left(W^{(\bdeux)}_{3,-4}\right)\chi_q^-\left(W^{(\bdeux)}_{1,-2}\right)
+\chi_q^-\left(W^{(\bun)}_{1,-1}\right)\chi_q^-\left(W^{(\bun)}_{1,-3}\right)}
{\chi_q^-\left(W^{(\bdeux)}_{2,-2}\right)} 
=
\chi_q^-\left(W^{(\bdeux)}_{2,-4}\right).
\]
By induction, after every vertex of the second column
has been mutated, each cluster variable of the form $\chi_q^-\left(W^{(\bdeux)}_{k,-2k+2}\right)$
has been replaced by the new cluster variable $\chi_q^-\left(W^{(\bdeux)}_{k,-2k}\right)$.
We now continue by mutating vertices of the third column. We first get, at the top vertex
\[
\frac{\chi_q^-\left(W^{(\bun)}_{2,-5}\right)+\chi_q^-\left(W^{(\bdeux)}_{2,-4}\right)}
{\chi_q^-\left(W^{(\bun)}_{1,-1}\right)} 
=
\chi_q^-\left(W^{(\bun)}_{1,-5}\right).
\] 
Then, mutating at the next vertex gives
\[
\frac{\chi_q^-\left(W^{(\bun)}_{3,-9}\right)\chi_q^-\left(W^{(\bun)}_{1,-5}\right)
+\chi_q^-\left(W^{(\bdeux)}_{4,-8}\right)}
{\chi_q^-\left(W^{(\bun)}_{2,-5}\right)} 
=
\chi_q^-\left(W^{(\bun)}_{2,-9}\right).
\] 
By induction one sees that, after every vertex of the third column
has been mutated, each cluster variable of the form $\chi_q^-\left(W^{(\bun)}_{k,-4k+3}\right)$
has been replaced by the new cluster variable $\chi_q^-\left(W^{(\bun)}_{k,-4k-1}\right)$.
For the third part of $\mu_\SS$, we mutate again along the second column. 
One checks that after that, each cluster variable of the form $\chi_q^-\left(W^{(\bdeux)}_{k,-2k}\right)$
produced after the first part of $\mu_\SS$ has been replaced by $\chi_q^-\left(W^{(\bdeux)}_{k,-2k-2}\right)$.
Finally, the fourth part of $\mu_\SS$ along the first column replaces  
each cluster variable of the form $\chi_q^-\left(W^{(\bun)}_{k,-4k+1}\right)$
by the new cluster variable $\chi_q^-\left(W^{(\bun)}_{k,-4k-3}\right)$.
Moreover, one sees that the new quiver obtained after $\mu_\SS$ is nothing
else than $G^-$. Hence we conclude that one application of $\mu_\SS$
produces a seed with the same quiver, and in which every cluster variable 
$\chi_q^-\left(W^{(i)}_{k,r}\right)$ has been replaced
by $\chi_q^-\left(W^{(i)}_{k,r-4}\right)$. In other words, the effect of
$\mu_\SS$ is merely a uniform shift of the spectral parameters $r$ by $-4$.

\begin{figure}[t]
\[
\def\objectstyle{\scriptstyle}
\def\labelstyle{\scriptstyle}
\xymatrix@-1.0pc{
&&&&\\
&&(\bdeux,0)\ar[rrddd]& 
\\
&&&
\\
&&\ar[llddd]\ar[uu] (\bdeux,-2)& &
\\
&&&& (\bun,-1)\ar[llddd] \\
&&\ar[uu] (\bdeux,-4) \ar[rrrddd]&&
\\
(\bun,-3)\ar[rrddd] &&
\\
&&\ar[uu] \ar[rrddd](\bdeux,-6) &
\\
&&&&&(\bun,-5)\ar[lllddd]  
\\
&& \ar[uu](\bdeux,-8)\ar[llddd] 
\\
&&&& (\bun,-7)\ar[uuuuuu]\ar[llddd] &
\\
&&(\bdeux, -10)\ar[uu]\ar[rrrddd] & &
\\
(\bun,-9)\ar[uuuuuu]&  
\\
&&(\bdeux, -12)\ar[uu] & &
\\
{}\ar[uu]&&{}\ar[u]&&{}\ar[uuuu]&(\bun,-11)\ar[uuuuuu]
\\
{}\save[]+<0cm,2ex>*{\vdots}\restore&&{}\save[]+<0cm,2ex>*{\vdots}\restore
&&{}\save[]+<0cm,2ex>*{\vdots}\restore&{}\save[]+<0cm,2ex>*{\vdots}\restore
}
\]
\caption{\label{Fig2} {\it The quiver $G^-$ for $\g$ of type $G_2$.}}
\end{figure}
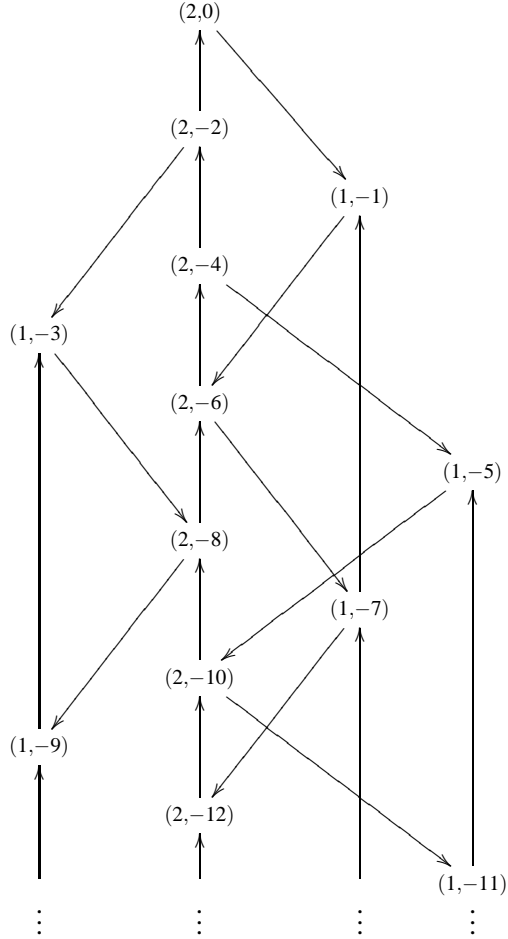
The argument is similar for $\g$ of type $G_2$. The quiver $G^-$
for this case is dispayed in Figure~\ref{Fig2}, and the mutation sequence is  
\[
\begin{array}{c}
(\bdeux,0), (\bdeux,-2), (\bdeux,-4), \ldots, (\bun,-1), (\bun,-7), (\bun,-13), \ldots, \\[1mm]
(\bdeux,0), (\bdeux,-2), (\bdeux,-4), \ldots, (\bun,-3), (\bun,-9), (\bun,-15), \ldots, \\[1mm]
(\bdeux,0), (\bdeux,-2), (\bdeux,-4), \ldots, (\bun,-5), (\bun,-11), (\bun,-17), \ldots. 
\end{array}
\]
The sequence of mutated quivers obtained at each step of $\mu_\SS$
is displayed in Appendix \S\ref{appendixG2}. 

For a general $\g$, we use a reduction to rank two.
Namely, we show that mutation sequences and $T$-systems equations 
are compatible with rank two reductions. 

First, by construction, the sequence of vertices $\mathcal{S}$ is a union of $tn$ columns:
\[
\mathcal{S} = (\mathcal{S}_1,\mathcal{S}_2,\cdots, \mathcal{S}_{tn}),
\]
where each column $\mathcal{S}_k$ is a subset of $i_k\times \mathbb{Z}_{\leq 0}$
for a certain $i_k\in {\bf I}$. As above we use $\mu_{\mathcal{S}_k}$ to denote the sequence of mutations
indexed by $\mathcal{S}_k$. So we have
\[
\mu_{\mathcal{S}} = \mu_{\mathcal{S}_{tn}}\circ \mu_{\mathcal{S}_{tn - 1}}\circ\cdots \circ \mu_{\mathcal{S}_1}.
\]
For $0\leq k\leq tn$, we get the mutated quiver
\[
\Sigma_k = (\mu_{\mathcal{S}_k}\circ \mu_{\mathcal{S}_{k - 1}}\circ\cdots \circ \mu_{\mathcal{S}_1})(\Sigma).
\]
For a subset ${J}\subset I$, let us denote by $(\Sigma_k)_{J}$ the subquiver of $\Sigma_k$ obtained 
by deleting the vertices $(i,r)$ such that $i\notin {J}$, and 
the edges whose tail or head is such a vertex. 
For any ${i} \in {I}$, the mutation sequence $\mu_{\mathcal{S}_k}$ modifies $(\Sigma_k)_{i}$ to itself.
Consequently, $(\Sigma_k)_{i} = (\Sigma)_{i}$ does not depend on $k$ 
(it is a disjoint union of $d_i$ semi-infinite linear quivers). 
Besides, the mutation sequence $\mu_{\mathcal{S}_k}$ modifies only the edges whose tail (resp. head) 
is in ${i_k}\times \mathbb{Z}$ 
and head (resp. tail) is in ${j}\times \mathbb{Z}$ where $c_{{i_k}{j}} < 0$.
This is because each mutation of the sequence takes place at a vertex $(i_k,r)$
having two incoming arrows from vertices $(i_k,r\pm d_i)$ and outgoing arrows to
vertices of the form $(j,s)$ with $c_{i_kj} < 0$.
Consequently, for each ${i}\neq {j}$ in ${I}$, the effect of the mutation sequence $\mu_{\mathcal{S}}$ 
on $(\Sigma)_{\{{i}, {j}\}}$ is the same as the effect of an iteration of the mutation sequence 
corresponding to the rank two Lie subalgebra of $\g$ attached to $\{{i}, {j}\}\subset {I}$. 
But we have established the result for rank two Lie algebras, so this implies
$$(\mu_{\mathcal{S}}(\Sigma))_{\{{i}, {j}\}}
= (\Sigma)_{\{{i}, {j}\}}.$$ 
As this is true for any ${i}\neq {j}$ in ${I}$, we get
$\mu_{\mathcal{S}}(\Sigma) = \Sigma$. 

Secondly, a $T$-system equation involves only a certain index ${i}\in {I}$ and the 
indices ${j} \in {I}$ with
$c_{{i}{j}} < 0$. 
The $T$-system equations do not change by reduction, in the sense that for such a $j$, 
the powers of the factors $T_{l,s}^{({j})}$ in the second term $S_{k,r}^{({i})}$ of the right-hand side
of~(\ref{Tsystem_general}) are the same as for the $T$-system equation associated with the 
rank two Lie subalgebra of $\g$ attached to $\{{i}, {j}\}$. 
Combining with our results above for the subquivers $(\Sigma_k)_{\{{i}, {j}\}}$, 
we have proved that, for a general $\g$, all exchange relations of cluster variables of 
our mutation sequence are in fact $T$-system equations. 
Moreover, the mutation sequence $\mu_\SS$ replaces
the initial seed $\Si$ by a seed with the same quiver; the cluster variables,
expressed in terms of the $Y_{i,r}$ via (\ref{chvar}), are truncated $q$-characters
of the same Kirillov-Reshetikhin modules, the only difference being that their spectral parameters 
are uniformly shifted by $-2t$. 

Hence, after $m$ applications of $\mu_\SS$ we
will get the truncated $q$-characters
\[
y^{(m)}_{i,r} = \chi_q^-\left(W^{(i)}_{k_{i,r},\ r-2tm}\right).  
\]
Now taking into account \cite[Corollary 6.14]{FM}, we see that if $2tm \ge t\dh$, then
all the monomials of the $q$-character of $W^{(i)}_{k_{i,r},\ r-2tm}$ are
lower than the level of truncation, that is,
\[
 \chi_q^-\left(W^{(i)}_{k_{i,r},\ r-2tm}\right) = 
 \chi_q\left(W^{(i)}_{k_{i,r},\ r-2tm}\right). 
\]
This finishes the proof of Theorem~\ref{thm1}.


\section{A geometric character formula for Kirillov-Reshetikhin modules}\label{section3}

\subsection{Semi-infinite quivers with potentials}\label{potential}

Recall the map $\psi \colon V \to W$ of \S\ref{subsubsect_simply}.
Put $V^-:=\psi^{-1}(W^-)$, and denote by $\G^-$ the full
subquiver of $\G$ with vertex set $V^-$. Thus $\G^-$
is the same graph as $G^-$, but with a change of labelling
of its vertices. (Compare for instance Figure~\ref{Fig2} and Figure~\ref{Fig6}.)

For every $i \not = j$ in $I$ with $c_{ij}\not = 0$, and every
$(i,m)$ in $V^-$,
we have in $\G^-$ an oriented cycle:
\begin{equation}\label{cycle}
\xymatrix@-1.0pc{
(i,m)\ar[rrdd] \\
(i,m-b_{ii})\ar[u]\\
{}\save[]{\vdots}\restore&&(j,m+b_{ij})\ar[lldd]\\
(i,m+2b_{ij}+b_{ii})\\
(i,m+2b_{ij})\ar[u]
}
\end{equation}
There are $2|b_{ij}|/b_{ii} = |c_{ij}|$
consecutive vertical up arrows, hence this cycle has length $2+|c_{ij}|$.
We define a \emph{potential} $S$ as the formal sum of all these
oriented cycles up to cyclic permutations, see \cite[\S3]{DWZ}. 
This is an infinite sum, but note that a given arrow
of $\G^-$ can only occur in a finite number of summands.
Hence all the cyclic derivatives of $S$, defined as in \cite[Definition 3.1]{DWZ}, 
are finite sums of paths in~$\G^-$. 
Let $R$ be the list of all cyclic derivatives of $S$.
Let $J$ denote the two-sided
ideal of the path algebra $\C\G^-$ generated by $R$.
Following \cite{DWZ}, we now introduce 
\begin{Def}
Let $A$ be the infinite-dimensional $\C$-algebra $\C\G^-\slash J$. 
\end{Def}
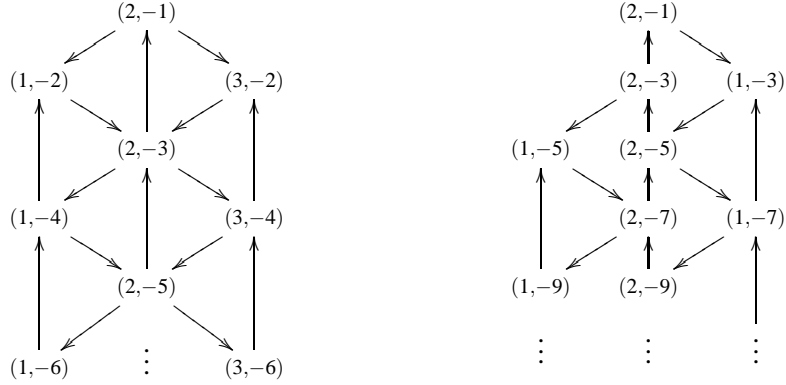
\begin{figure}[t]
\[
\def\objectstyle{\scriptstyle}
\def\lablestyle{\scriptstyle}
\xymatrix@-1.0pc{
&&&&\\
&&\ar[ld] (2,-1) \ar[rd]&&
\\
&(1,-2) \ar[rd] &&\ar[ld] (3,-2)
\\
&&\ar[ld] \ar[uu](2,-3) \ar[rd]&&
\\
&\ar[uu](1,-4) \ar[rd] &&\ar[ld] (3,-4)\ar[uu]
\\
&&\ar[ld] \ar[uu](2,-5) \ar[rd]&&
\\
&\ar[uu](1,-6) &{}\save[]{\vdots}\restore&\ar[uu] (3,-6) 
}
\qquad\qquad
\xymatrix@-1.0pc{
&&&&\\
&  & (\bdeux,-1) \ar[rd]&&
\\
&&\ar[u]\ar[ld] (\bdeux,-3) &\ar[ld] (\bun,-3) &
\\
&\ar[rd](\bun,-5) & \ar[u]\ar[rd](\bdeux,-5)&&
\\
&&\ar[u]\ar[ld] (\bdeux,-7) &\ar[ld] (\bun,-7) \ar[uu]  &
\\
&\ar[uu](\bun,-9)   & \ar[u](\bdeux,-9) &&
\\
&{}\save[]+<0cm,0ex>*{\vdots}\restore  &{}\save[]+<0cm,0ex>*{\vdots}\restore&{}\save[]+<0cm,0ex>*{\vdots}\ar[uu]\restore 
\\
}
\]
\caption{\label{Fig3} {\it The quivers $\G^-$ for $\g$ of type $A_3$ and $B_2$.}}
\end{figure}
\begin{example}
{\rm Let $\g$ be of type $A_3$.
Then $\G^-$ is the first graph in Figure~\ref{Fig3}.
The ideal $J$ is generated by the following 7 families of linear combinations of paths,
for every $m\in\Z_{<0}$,
\[
\begin{array}{l}
((1,2m), (2,2m-1), (1,2m-2)), \\[2mm]
((3,2m), (2,2m-1), (3,2m-2)), \\[2mm]
((1,2m), (1,2m+2), (2,2m+1)) + ((1,2m), (2,2m-1), (2,2m+1)),\\[2mm]
((3,2m), (3,2m+2), (2,2m+1)) + ((3,2m), (2,2m-1), (2,2m+1)),\\[2mm]
((2,2m-1), (1,2m-2), (1,2m)) + ((2,2m-1), (2,2m+1), (1,2m)),\\[2mm]
((2,2m-1), (3,2m-2), (3,2m)) + ((2,2m-1), (2,2m+1), (3,2m)),\\[2mm]
((2,2m+1), (1,2m), (2,2m-1)) + ((2,2m+1), (3,2m), (2,2m-1)).
\end{array}
\]
Here, using the fact that there is at most one arrow between two vertices of $\G^-$,
we have denoted unambiguously paths by sequences of vertices. Thus 
$((1,2m), (2,2m-1), (1,2m-2))$ denotes the path of length 2 starting at $(1,2m)$,
passing by $(2,2m-1)$ and ending in $(1,2m-2))$.
Also, for $m=-1$, the third and fourth linear combinations of paths reduce 
respectively to the single paths
\[
((1,-2), (2,-3), (2,-1)) \qquad \mbox{and}\qquad ((3,-2), (2,-3), (2,-1)). 
\]

}
\end{example}

\begin{example}\label{A_relations_B2}
{\rm 
Let $\g$ be of type $B_2$.
Then $\G^-$ is the second graph of Figure~\ref{Fig3}.
The ideal $J$ is generated by the following 4 families of linear combinations of paths,
for every $m\in\Z_{<0}$,
\[
\begin{array}{l}
((\bun,2m-1),(\bdeux,2m-3),(\bun,2m-5)),\\[2mm]
((\bun,2m-1),(\bun,2m+3),(\bdeux,2m+1)) + ((\bun,2m-1),(\bdeux,2m-3),(\bdeux,2m-1),(\bdeux,2m+1)),\\[2mm] 
((\bdeux,2m-3),(\bun,2m-5),(\bun,2m-1)) + ((\bdeux,2m-3),(\bdeux,2m-1),(\bdeux,2m+1),(\bun,2m-1)),\\[2mm]
((\bdeux,2m+1),(\bdeux,2m+3),(\bun,2m+1),(\bdeux,2m-1)) + ((\bdeux,2m+1),(\bun,2m-1),(\bdeux,2m-3),(\bdeux,2m-1)).
\end{array}
\]
For $m=-1$ and $m=-2$ the second linear combinations of paths reduce respectively to the single paths
\[
((\bun,-3), (\bdeux,-5),(\bdeux,-3),(\bdeux,-1))
\qquad\mbox{and}\qquad
((\bun,-5), (\bdeux,-7),(\bdeux,-5),(\bdeux,-3)).
\]
For $m=-1$ the fourth linear combination of paths reduces to the single path
\[
((\bdeux,-1),(\bun,-3),(\bdeux,-5),(\bdeux,-3)). 
\]
}
\end{example}

\subsection{$F$-polynomials of $A$-modules}\label{ssect-F-pol-A-mod}

Let $M$ be a finite-dimensional $A$-module, and let $e \in \N^{V^-}$ be a dimension vector.
Let $\Gr_e(M)$ be the variety of submodules of $M$ with dimension vector $e$.
This is a projective complex variety, and we denote by $\chi(\Gr_e(M))$ its Euler characteristic.  
Following \cite{DWZ2}, consider the polynomial
\begin{equation}\label{FpolAmod}
F_M=\sum_{e \in \N^{V^-}} \chi(\Gr_e(M)) \prod_{(i,r)\in V^-} v_{i,r}^{e_{i,r}} 
\end{equation}
in the indeterminates $v_{i,r}\ ((i,r)\in V^-)$, called the \emph{$F$-polynomial} of $M$.
Note that, for $\Gr_e(M)$ to be nonempty, one must take $e$ between $0$ and the dimension vector of $M$
(componentwise). Moreover, if $e=0$ or $e=\dimv(M)$, the variety $\Gr_e(M)$ is just a point,
so $F_M$ is a monic polynomial with constant term equal to 1.

In the sequel, we shall evaluate the variables of the $F$-polynomials at the inverses of the variables
$A_{i,r}$ introduced in (\ref{def_A_var}), namely:
\begin{equation}\label{eval_v}
v_{i,r} := A_{i,r}^{-1} =
 Y_{i,r-d_i}^{-1}Y_{i,r + d_i}^{-1} 
\prod_{j:\ c_{ji} = -1} Y_{j,r}
\prod_{j:\ c_{ji} = -2} Y_{j,r-1}Y_{j,r+1}
\prod_{j:\ c_{ji} = -3} Y_{j,r-2}Y_{j,r}Y_{j,r+2}.
\end{equation}

\subsection{Generic kernels}\label{ssect-A-mod}

Suppose that $X$ and $Y$ are $A$-modules such that $\Hom_A(X,Y)$ 
is finite-dimensional.
Assume also that there exists $f\in \Hom_A(X,Y)$ such that $\Ker(f)$ is finite-dimensional.
Then, there is an open dense subset $\widetilde{O}$ of 
$\Hom_A(X,Y)$ such that the kernels of all elements of $\widetilde{O}$ are finite-dimensional. 
Moreover, since the map sending a homomorphism $f$ to 
the $F$-polynomial of $\Ker(f)$ is constructible (see \cite[\S2]{Pa}), 
$\widetilde{O}$~contains an open dense subset $O$ of $\Hom_A(X,Y)$ such that
the $F$-polynomials of the kernels of all elements of $O$ coincide.
We shall say that an element of $O$ is a \emph{generic homomorphism}
from $X$ to $Y$.

Let us denote by $S_{i,m}$ the one-dimensional $A$-module
supported on $(i,m)\in V^-$.
Let $I_{i,m}$ be the (infinite-dimensional) injective 
$A$-module with socle isomorphic to $S_{i,m}$.
The $\C$-vector space $I_{i,m}$ has a basis indexed by classes modulo $J$
of paths in $\G^-$ with final vertex $(i,m)$. 
In particular, for every $k\ge 0$ we have in $\G^-$ a path 
\begin{equation}\label{path}
((i,m-kb_{ii}), (i,m-(k-1)b_{ii}), \ldots, (i,m))  
\end{equation}
of length $k$ from $(i,m-kb_{ii})$ to $(i,m)$,
whose class modulo $J$ is nonzero. Thus the $(i,m-kb_{ii})$-component
of the dimension vector of $I_{i,m}$ is nonzero, and it follows that
\begin{equation}
\Hom_A(I_{i,m},I_{i,m-kb_{ii}}) \not = 0, \qquad ((i,m)\in V^-,\ k\ge 0). 
\end{equation}
More precisely, $\Hom_A(I_{i,m},I_{i,m-kb_{ii}})$ has finite dimension
equal to the $(i,m-kb_{ii})$-component
of the dimension vector of $I_{i,m}$.
The next Lemma will be proven in \S\ref{proofs}.
\begin{Lem}\label{finite-dim}
There exists $f\in \Hom_A(I_{i,m},I_{i,m-kb_{ii}})$ with $\Ker(f)$ finite-dimensional. 
\end{Lem}

Because of this lemma, the following definition makes sense.

\begin{Def}
Let $K^{(i)}_{k,m}$ be the kernel of a generic $A$-module homomorphism
from $I_{i,m}$ to $I_{i,m-kb_{ii}}$.
\end{Def}

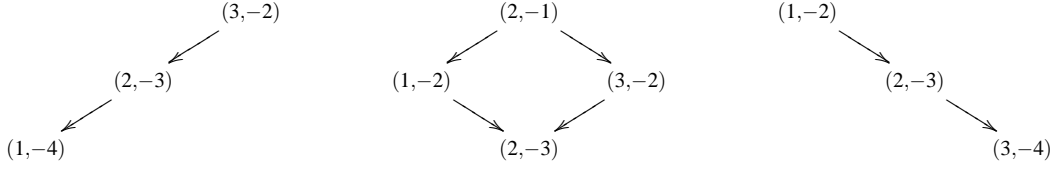
\begin{figure}[t]
\[
\def\objectstyle{\scriptstyle}
\def\lablestyle{\scriptstyle}
\xymatrix@-1.0pc{
&&
&\ar[ld] (3,-2) 
\\
&&\ar[ld] (2,-3) &&
\\
&(1,-4)  && 
\\
}
\def\objectstyle{\scriptstyle}
\def\lablestyle{\scriptstyle}
\xymatrix@-1.0pc{
&&\ar[ld] (2,-1) \ar[rd]&&
\\
&{(1,-2)}\ar[rd]&
&\ar[ld] (3,-2) 
\\
&& (2,-3) &&
\\
}
\def\objectstyle{\scriptstyle}
\def\lablestyle{\scriptstyle}
\xymatrix@-1.0pc{
&{(1,-2)}\ar[rd]&
& 
\\
&& (2,-3) \ar[rd]&&
\\
&  && (3,-4)
\\
}
\]
\caption{\label{Fig4} {\it The modules $K^{(1)}_{1,-4}$, $K^{(2)}_{1,-3}$, $K^{(3)}_{1,-4}$ for $\g$ of type $A_3$.}}
\end{figure}

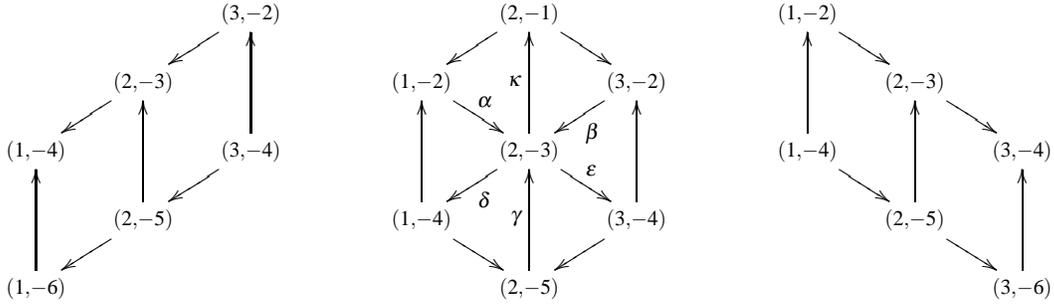
\begin{figure}[t]
\[
\def\objectstyle{\scriptstyle}
\def\lablestyle{\scriptstyle}
\xymatrix@-1.0pc{
&&
&\ar[ld] (3,-2) 
\\
&&\ar[ld] (2,-3) &&
\\
&(1,-4)&&\ar[ld] (3,-4)\ar[uu] 
\\
&&\ar[ld] (2,-5)\ar[uu] &&
\\
&(1,-6)\ar[uu]  && 
\\
}
\def\objectstyle{\scriptstyle}
\def\lablestyle{\scriptstyle}
\xymatrix@-1.0pc{
&&\ar[ld] (2,-1) \ar[rd]&&
\\
&{(1,-2)}\ar[rd]^{\alpha}&
&\ar[ld]^{\beta} (3,-2) 
\\
&& (2,-3)\ar[uu]^{\kappa}\ar[ld]^{\delta}\ar[rd]^{\epsilon} &&
\\
&{(1,-4)}\ar[rd]\ar[uu]&
&\ar[ld] (3,-4)\ar[uu] 
\\
&& (2,-5)\ar[uu]^{\gamma} &&
\\
}
\def\objectstyle{\scriptstyle}
\def\lablestyle{\scriptstyle}
\xymatrix@-1.0pc{
&{(1,-2)}\ar[rd]&
& 
\\
&& (2,-3) \ar[rd]&&
\\
&\ar[uu]{(1,-4)}\ar[rd]  && (3,-4)
\\
&& \ar[uu](2,-5) \ar[rd]&&
\\
&  && \ar[uu](3,-6)
\\
}
\]
\caption{\label{Fig5} {\it The modules $K^{(1)}_{2,-4}$, $K^{(2)}_{2,-3}$, and $K^{(3)}_{2,-4}$ for $\g$ of type $A_3$.}}
\end{figure}

\begin{example}\label{example_A3_modules}
{\rm 
Figure~\ref{Fig4} and Figure~\ref{Fig5} show the structure of some modules 
$K^{(i)}_{k,m}$ in type $A_3$. 
Our convention for displaying these quiver representations 
is the following. We only keep the vertices of $\G^-$ whose corresponding vector space is nonzero,
and the arrows whose corresponding linear map is nonzero.
Moreover, in these small examples, almost all vertices carry a vector space of dimension 1.
The only exception is the module $K^{(2)}_{2,-3}$ in Figure~\ref{Fig5}, whose vertex
$(2,-3)$ carries a vector space of dimension 2. The maps associated with the arrows incident
to this vertex have the following matrices  
\[
 \alpha = \beta = \gamma = \pmatrix{1\cr 0},
 \qquad
 \delta = \epsilon = \kappa = \pmatrix{0 & 1}.
\]
All other arrows carry linear maps with matrix $\pmatrix{\pm 1}$, whose sign is easily deduced
from the defining relations of $A$.

It is a nice exercise to check that the modules shown in Figure~\ref{Fig4} and Figure~\ref{Fig5}
are indeed the claimed modules $K^{(i)}_{k,m}$ (see also Example~\ref{K-ADE} below). 
For instance, one can easily see that the $(1,-6)$-component of the dimension 
vector of $I_{1,-4}$ is equal to 1.
Hence $\Hom_A(I_{1,-4},I_{1,-6})$ is of dimension 1, and 
$K^{(1)}_{1,-4}$ is the kernel of any nonzero homomorphism. 
It is also easy to see that the $(2,-5)$-component of the dimension 
vector of $I_{2,-3}$ is equal to 2. In this case 
we have a stratification of the 2-plane $\Hom_A(I_{2,-3},I_{2,-5})$
with three strata of dimension $0$, $1$, $2$. The module 
$K^{(2)}_{1,-3}$ is the kernel of any homomorphism in the open stratum,
that is, of any surjective homomorphism.
The image of any homomorphism in the one-dimensional stratum is
the unique submodule $X$ of $I_{2,-5}$ with dimension vector given by
\[
\dim(X_{i,m}) = 
\left\{
\begin{array}{cl}
1 & \mbox{if $i=2$ and $m = -5-2j$ for some $j\in\N$,}\\[2mm]
0 & \mbox{otherwise.}
\end{array}
\right.
\]
The kernel of such a homomorphism is infinite-dimensional.
}
\end{example}

\begin{example}\label{K-ADE}
{\rm
Let us assume that  
$\g$ is of type $A$, $D$, $E$.
In this case, the modules $K^{(i)}_{1,r}$
are closely related to the indecomposable injective modules over the 
preprojective algebra $\L$ of~$\de$. 

Consider the subalgebra $\widetilde{\L}$ of $A$
generated by the images modulo $J$ of the arrows of $\Gamma^-$
of the form $(i,m) \to (j,m-1)$, for every edge between $i$ and $j$ in $\de$,
and every $(i,m)\in V^-$.
In other words, if $\De^-_\de$ is the subquiver of
$\Gamma^-$ obtained by erasing all the vertical arrows $(i,m-2) \to (i,m)$,
then $\widetilde{\L}$ is isomorphic to the quotient of $\C\De^-_\de$
by the two-sided ideal generated by the relations
\[
\sum_{j:\ c_{ij}<0} ((i,m),(j,m-1),(i,m-2)) = 0,\qquad ((i,m)\in V^-).  
\]
Thus, $\widetilde{\L}$ is a $\Z_{<0}$-graded version of $\L$.
We can of course regard the simple $A$-module $S_{i,r}$ as a $\widetilde{\L}$-module.
Let $H_{i,r}$ be the injective $\widetilde{\L}$-module with socle $S_{i,r}$.
Then $H_{i,r}$ is finite-dimensional.
More precisely,
for $r\le 1-h$, $H_{i,r}$ is just a graded version of the 
indecomposable injective $\L$-module $I_i$ with socle the one-dimensional
$\L$-module $S_i$ supported on vertex $i$ of $\de$. For $r > 1-h$, 
$H_{i,r}$ is a graded version of a submodule of $I_i$. 

Any $\widetilde{\L}$-module $X$ can be given the structure of an $A$-module by letting the 
vertical arrows $(i,m-2) \to (i,m)$ act by 0 on $X$. In particular we can regard
$H_{i,r}$ as a finite-dimensional $A$-module.
Then one can check that $I_{i,r}$ has a unique submodule isomorphic to $H_{i,r}$,
giving rise to a non-split short exact sequence
\[
 0 \to H_{i,r} \to I_{i,r} \to I_{i,r-2} \to 0, \qquad ((i,r)\in V^-).
\]
It follows that 
the module $K^{(i)}_{1,m}$ is isomorphic to $H_{i,m}$.
In particular, when $m\le 1-h$, $K^{(i)}_{1,m}$ is a graded version of the injective
$\L$-module $I_i$.
}
\end{example}

\subsection{A geometric character formula}

Recall the $A$-module $K^{(i)}_{k,r}$ defined in \S\ref{ssect-A-mod}. 
We can now state our second main result.
\begin{Thm}\label{th_geom_form}
Let $(i,r)\in V^-$ and $k\in \N$. The 
$F$-polynomial of $K^{(i)}_{k,r}$ is equal to the normalized truncated $q$-character of the Kirillov-Reshetikhin
module $W^{(i)}_{k,\ r-(2k-1)d_i}$. More precisely, we have
\[
\chi_q^-\left(W^{(i)}_{k,\ r-(2k-1)d_i}\right) = \left(\prod_{s=1}^{k}  Y_{i,\ r-(2s-1)d_i}\right) F_{K^{(i)}_{k,r}}, 
\]
where the variables $v_{i,r}$ of the $F$-polynomial are evaluated as in {\rm(\ref{eval_v})}.
\end{Thm}

\begin{remark}\label{complete-character}
{\rm
If $r\le d_i-t\dh$, then the truncated $q$-character of $W^{(i)}_{k,\ r-(2k-1)d_i}$
is equal to the complete $q$-character.
Hence, Theorem~\ref{th_geom_form} gives a geometric formula for the $q$-character of any Kirillov-Reshetikhin module 
(up to a spectral shift).
}
\end{remark}
\begin{remark}
{\rm
If $M$ and $N$ are two finite-dimensional $A$-modules, then $F_{M\oplus N} = F_M F_N$
\cite[Proposition 3.2]{DWZ2}. It follows immediately that, replacing in Theorem~\ref{th_geom_form} the
module $K^{(i)}_{k,r}$ by a direct sum of such modules, we obtain a similar geometric
character formula for arbitrary tensor products of Kirillov-Reshetikhin modules.
In particular, we get a geometric formula for the standard modules, which are 
isomorphic to tensor products of fundamental modules.
}
\end{remark}

\begin{remark}
{\rm
Let $\g$ be of type $A$, $D$, $E$.
Let $\bV$ and $\bW$ be finite-dimensional vector spaces graded by~$V^-$.
In \cite{N1} (see also \cite{N3}), Nakajima has introduced a graded quiver variety 
$\mathfrak{L}^\bullet(\bV,\bW)$ and has endowed
the sum of cohomologies 
\[
\bigoplus_{\bV} H^*(\mathfrak{L}^\bullet(\bV,\bW)) 
\]
with the structure of a standard $U_q(\hg)$-module, 
with highest weight encoded by $\bW$.
It was proved by Lusztig (in the ungraded case), and by Savage and Tingley (in the graded case),
that $\mathfrak{L}^\bullet(\bV,\bW)$ is homeomorphic to a Grassmannian of submodules of an 
injective module over the graded preprojective algebra (see \cite[\S2.8]{L}).
Therefore, using the description of $K^{(i)}_{1,r}$ given in Example~\ref{K-ADE}, we see that
the varieties 
\[
\Gr_e\left(\bigoplus_{(i,r)}\left(K^{(i)}_{1,r}\right)^{\oplus a_{i,r}}\right) 
\]
involved in our geometric $q$-character formula for standard modules in the
simply laced case are homeomorphic to certain Nakajima varieties $\mathfrak{L}^\bullet(\bV,\bW)$.
Here, the multiplicities $a_{i,r}$ are the dimensions of the graded components of $\bW$,
and we assume that $a_{i,r}=0$ if $r>1-h$.
Similarly the graded dimension of $\bV$ is encoded by the dimension vector $e$.
}
\end{remark}

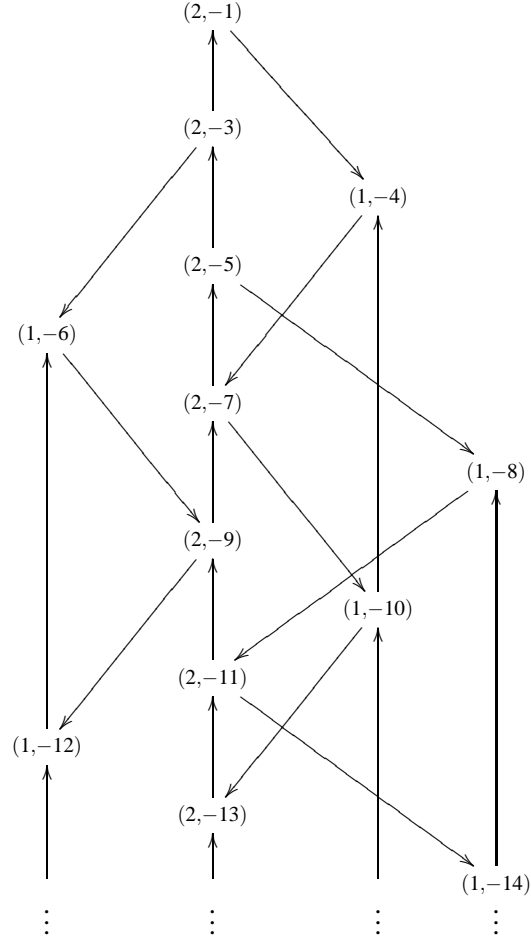
\begin{figure}[t]
\[
\def\objectstyle{\scriptstyle}
\def\labelstyle{\scriptstyle}
\xymatrix@-1.0pc{
&&&&\\
&&(\bdeux,-1)\ar[rrddd]& 
\\
&&&
\\
&&\ar[llddd]\ar[uu] (\bdeux,-3)& &
\\
&&&& (\bun,-4)\ar[llddd] \\
&&\ar[uu] (\bdeux,-5) \ar[rrrddd]&&
\\
(\bun,-6)\ar[rrddd] &&
\\
&&\ar[uu] \ar[rrddd](\bdeux,-7) &
\\
&&&&&(\bun,-8)\ar[lllddd]  
\\
&& \ar[uu](\bdeux,-9)\ar[llddd] 
\\
&&&& (\bun,-10)\ar[uuuuuu]\ar[llddd] &
\\
&&(\bdeux, -11)\ar[uu]\ar[rrrddd] & &
\\
(\bun,-12)\ar[uuuuuu]&  
\\
&&(\bdeux, -13)\ar[uu] & &
\\
{}\ar[uu]&&{}\ar[u]&&{}\ar[uuuu]&(\bun,-14)\ar[uuuuuu]
\\
{}\save[]+<0cm,2ex>*{\vdots}\restore&&{}\save[]+<0cm,2ex>*{\vdots}\restore
&&{}\save[]+<0cm,2ex>*{\vdots}\restore&{}\save[]+<0cm,2ex>*{\vdots}\restore
}
\]
\caption{\label{Fig6} {\it The quiver $\G^-$ for $\g$ of type $G_2$.}}
\end{figure}

\begin{example}
{\rm
Let $\g$ be of type $A_3$. We have
\[
v_{1,r} = Y_{1,r-1}^{-1} Y_{1,r+1}^{-1} Y_{2,r},\qquad
v_{2,r} = Y_{2,r-1}^{-1} Y_{2,r+1}^{-1} Y_{1,r} Y_{3,r},\qquad
v_{3,r} = Y_{3,r-1}^{-1} Y_{3,r+1}^{-1} Y_{2,r}.
\]
We continue Example~\ref{example_A3_modules}.
The submodule structure of the $A$-modules displayed in Figure~\ref{Fig4} is very simple.
Indeed, in this case, all the nonempty varieties $\Gr_e(K^{(i)}_{k,r})$ are reduced to
a single point, and their Euler characteristics are equal to 1.
Therefore the $F$-polynomial reduces to a generating polynomial for the dimension vectors of the 
(finitely many) submodules of $K^{(i)}_{k,r}$.
This yields the following well known formulas for the $q$-characters of the fundamental modules:
\[
\begin{array}{rcl}
\chi_q(L(Y_{1,-5})) &=& 
Y_{1,-5}(1 + v_{1,-4} + v_{1,-4}v_{2,-3}
+ v_{1,-4}v_{2,-3}v_{3,-2})\\[1mm]
&=& Y_{1,-5} + Y_{1,-3}^{-1}Y_{2,-4} + Y_{2,-2}^{-1}Y_{3,-3} + Y_{3,-1}^{-1},
\\[3mm]
\chi_q(L(Y_{2,-4})) &=& 
Y_{2,-4}(1 + v_{2,-3} + v_{1,-2}v_{2,-3}
+ v_{2,-3}v_{3,-2} + v_{1,-2}v_{2,-3}v_{3,-2}\\[1mm] 
&&+\ v_{1,-2}v_{2,-3}v_{3,-2}v_{2,-1})\\[1mm]
&=&Y_{2,-4} + Y_{1,-3}Y_{2,-2}^{-1}Y_{3,-3} + Y_{1,-1}^{-1}Y_{3,-3} + Y_{1,-3}Y_{3,-1}^{-1} 
+ Y_{1,-1}^{-1}Y_{2,-2}Y_{3,-1}^{-1}+Y_{2,0}^{-1},
\end{array}
\]
Similarly, the $A$-modules shown in Figure~\ref{Fig5}
give the following Kirillov-Reshetikhin $q$-characters:
\[
\begin{array}{rcl}
\chi_q(L(Y_{1,-7}Y_{1,-5})) &=& 
Y_{1,-7}Y_{1,-5}(1 + v_{1,-4}(1 + v_{1,-6} + v_{2,-3}
+ v_{1,-6}v_{2,-3} + v_{2,-3}v_{3,-2}\\[1mm]
&&
+\ v_{1,-6}v_{2,-3}v_{2,-5}
+ v_{1,-6}v_{2,-3}v_{3,-2}
+ v_{1,-6}v_{2,-3}v_{2,-5}v_{3,-2}\\[1mm]
&&+\ v_{1,-6}v_{2,-3}v_{2,-5}v_{3,-2}v_{3,-4})),\\[3mm]
\chi_q(L(Y_{2,-6}Y_{2,-4})) &=& 
Y_{2,-6}Y_{2,-4}(1 + v_{2,-3}(1 + v_{1,-2} + v_{2,-5} + v_{3,-2}
+ v_{1,-2}v_{2,-5} + v_{1,-2}v_{3,-2} 
\\[1mm]
&&
+\ v_{2,-5}v_{3,-2}
+ v_{1,-2}v_{2,-5}v_{3,-2}
+ v_{1,-2}v_{2,-5}v_{1,-4}
+ v_{1,-2}v_{3,-2}v_{2,-1}
\\[1mm]
&&
+\  v_{2,-5}v_{3,-2}v_{3,-4}
+ v_{1,-2}v_{2,-5}v_{3,-2}v_{1,-4}
+ v_{1,-2}v_{2,-5}v_{3,-2}v_{2,-1}
\\[1mm]
&&
+\ v_{1,-2}v_{2,-5}v_{3,-2}v_{3,-4}
+ v_{1,-2}v_{2,-5}v_{3,-2}v_{1,-4}v_{2,-1}
\\[1mm]
&&
+\ v_{1,-2}v_{2,-5}v_{3,-2}v_{1,-4}v_{3,-4}
+ v_{1,-2}v_{2,-5}v_{3,-2}v_{3,-4}v_{2,-1}
\\[1mm]
&&
+\ v_{1,-2}v_{2,-5}v_{3,-2}v_{1,-4}v_{2,-1}v_{3,-4}
+ v_{1,-2}v_{2,-5}v_{3,-2}v_{1,-4}v_{2,-1}v_{3,-4}v_{2,-3})),
\end{array}
\]
We omit the $q$-characters $\chi_q(L(Y_{3,-5}))$ and $\chi_q(L(Y_{3,-5}Y_{3,-7}))$, since they are
readily obtained from $\chi_q(L(Y_{1,-5}))$ and $\chi_q(L(Y_{1,-5}Y_{1,-7}))$ via the symmetry 
$1\leftrightarrow 3$.
}
\end{example}

\begin{example}\label{sect_ex_G2}
{\rm 
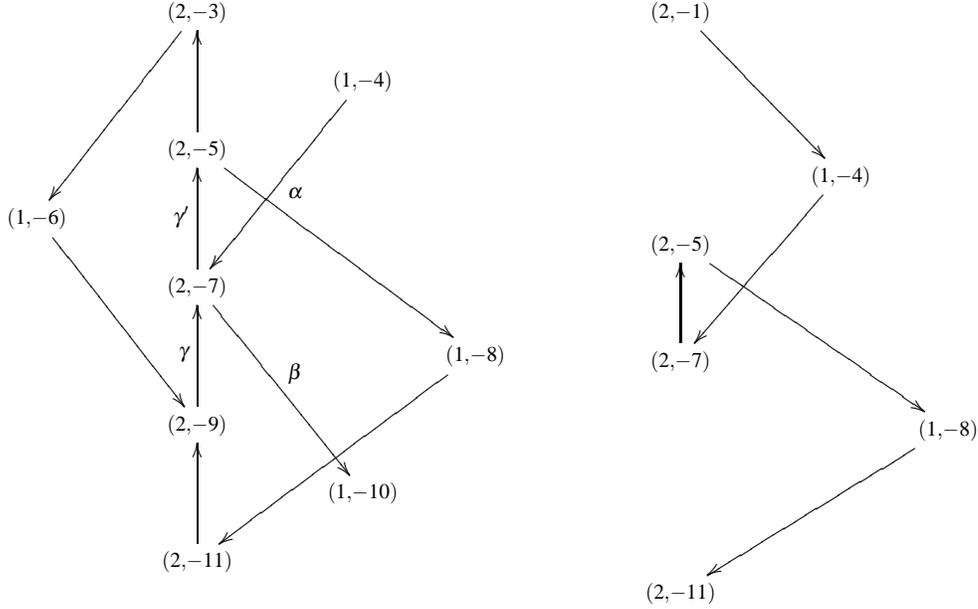
\begin{figure}[t]
\[
\def\objectstyle{\scriptstyle}
\xymatrix@-1.0pc{
&&&&\\
&&\ar[llddd](\bdeux,-3)& &
\\
&&&& (\bun,-4)\ar[llddd]^{\a} \\
&&\ar[uu] (\bdeux,-5) \ar[rrrddd]&&
\\
(\bun,-6)\ar[rrddd] &&
\\
&&\ar[uu]^{\ga'} \ar[rrddd]^{\b}(\bdeux,-7) &
\\
&&&&&(\bun,-8)\ar[lllddd]  
\\
&& \ar[uu]^{\ga}(\bdeux,-9) 
\\
&&&& (\bun,-10) &
\\
&&(\bdeux, -11)\ar[uu] & &
}
\quad
\xymatrix@-1.0pc{
&&&&\\
&&(\bdeux,-1)\ar[rrddd]& 
\\
&&&
\\
&&& &
\\
&&&& (\bun,-4)\ar[llddd] \\
&& (\bdeux,-5) \ar[rrrddd]&&
\\
 &&
\\
&&\ar[uu](\bdeux,-7) &
\\
&&&&&(\bun,-8)\ar[lllddd]  
\\
&&  
\\
&&&&&
\\
&&(\bdeux, -11)& &
}
\]
\caption{\label{Fig7} {\it The modules $K^{(\bun)}_{1,-10}$ and $K^{(\bdeux)}_{1,-11}$ for $\g$ of type $G_2$.}}
\end{figure}

Let $\g$ be of type $G_2$, with the long
root being $\a_{\bun}$.
The quiver $\G^-$ is shown in Figure~\ref{Fig6}.
The modules $K^{(\bun)}_{1,r}$ and $K^{(\bdeux)}_{1,s}$ with $r\le -10$ and 
$s\le -11$ have dimension 10 and 6, respectively.
For instance, $K^{(\bun)}_{1,-10}$ and $K^{(\bdeux)}_{1,-11}$ are represented
in Figure~\ref{Fig7}.
In the module $K^{(\bun)}_{1,-10}$ the vector space sitting at 
vertex $(\bdeux,-7)$ has dimension 2 (all other spaces have dimension 1).
The maps incident to this space are given by the following matrices
(see Figure~\ref{Fig7}):
\[
\a = \pmatrix{0 \cr 1},\quad   
\b = \pmatrix{1 & 0},\quad   
\ga = \pmatrix{1 \cr 0},\quad  \ga' = \pmatrix{0 & 1}.
\]
The corresponding fundamental modules have dimension 
\[
\dim L(Y_{\bun,-13}) = 15,\qquad
\dim L(Y_{\bdeux,-12}) = 7. 
\]
The Grassmannians of submodules of $K^{(\bun)}_{1,-10}$ and $K^{(\bdeux)}_{1,-11}$ are in this case again all reduced to points,
and the formula of Theorem~\ref{th_geom_form} amounts to an enumeration of the dimension
vectors of all submodules. This gives 
\[
\begin{array}{rcl}
\chi_q(L(Y_{\bun,-13}))&=&Y_{\bun,-13}(1+v_{\bun,-10}(1+v_{\bdeux,-7}(1+v_{\bdeux,-9}(1+v_{\bun,-6}+v_{\bdeux,-11}
+v_{\bun,-6}v_{\bdeux,-11}\\[1mm]
&&\ +\ v_{\bun,-6}v_{\bdeux,-3} +v_{\bdeux,-11}v_{\bun,-8} +  v_{\bun,-6}v_{\bdeux,-11}v_{\bdeux,-3} 
+ v_{\bun,-6}v_{\bdeux,-11}v_{\bun,-8}
\\[1mm]
&&\ +\  v_{\bun,-6}v_{\bdeux,-11}v_{\bdeux,-3}v_{\bun,-8}(1+v_{\bdeux,-5}(1+v_{\bdeux,-7}(1+v_{\bun,-4}))))))),
\\[2mm]
\chi_q(L(Y_{\bdeux,-12}))&=&Y_{\bdeux,-12}(1+v_{\bdeux,-11}(1+v_{\bun,-8}(1+v_{\bdeux,-5}(1+v_{\bdeux,-7}(1+ 
v_{\bun,-4}(1+v_{\bdeux,-1})))))),
\end{array} 
\]
where, following (\ref{eval_v}), we have
\[
v_{\bun,r} = Y_{\bun,r+3}^{-1} Y_{\bun,r-3}^{-1} Y_{\bdeux,r+2} Y_{\bdeux,r} Y_{\bdeux,r-2},
\qquad
v_{\bdeux,r} = Y_{\bdeux,r+1}^{-1}  Y_{\bdeux,r-1}^{-1} Y_{\bun,r}.
\]
}
\end{example}

\begin{remark}
{\rm
Assuming Theorem~\ref{th_geom_form}, we can easily calculate the dimension vectors of the $A$-modules
$K^{(i)}_{1,r}$ for $r\le d_i-t\dh$.
Indeed, by \cite[Lemma 6.8]{FM}, the lowest monomial of $\chi_q(Y_{i,r-d_i})$
is equal to $Y_{\nu(i),r-d_i+t\dh}^{-1}$, where $\nu$ is the involution of
$I$ defined by $w_0(\a_{i})=-\a_{\nu(i)}$. Denote by $\left(d_{j,s}(K^{(i)}_{1,r})\right)$
the dimension vector of $K^{(i)}_{1,r}$. Then, we have
\[
Y_{\nu(i),r-d_i+t\dh}^{-1} =  Y_{i,r-d_i} \prod_{(j,s)\in V^-} v_{j,s}^{d_{j,s}(K^{(i)}_{1,r})},
\]
and using (\ref{eval_v}), this equation determines the numbers 
$d_{j,s}(K^{(i)}_{1,r})$.
In particular, if we introduce the \emph{ungraded} dimension vector $(d_j(i))$ of $K^{(i)}_{1,r}$ by
\[
   d_j(i) := \sum_s d_{j,s}(K^{(i)}_{1,r}),\qquad (r\le d_i-t\dh), 
\]
we can deduce from this the nice formula
\begin{equation}\label{dimensionK}
 \sum_{i,j \in I} d_j(i)\a_j = \sum_{\b\in \Phi_{>0}} \b, 
\end{equation}
where $\Phi_{>0}$ is the set of positive roots of $\g$.
This can be observed in Figure~\ref{Fig4} and Figure~\ref{Fig7}
(see also \S\ref{append-B2}, \S\ref{append-B3}, \S\ref{append-C3}, \S\ref{append-F4} below).
When $\g$ is of type $A$, $D$, $E$, as explained in Remark~\ref{K-ADE} the 
modules $K^{(i)}_{1,r}$ are graded versions of the indecomposable injective
modules over the preprojective algebra~$\L$, and formula~(\ref{dimensionK})
recovers a well known property of $\L$.
}
\end{remark}

\subsection{Proof of the theorem} 
The proof relies on Theorem~\ref{thm1}, and on the categorification of cluster algebras by means of quivers with potentials,
developed by Derksen, Weyman and Zelevinsky \cite{DWZ,DWZ2}.
This categorification provides (among other things) a description of cluster variables in terms 
of Grassmannians of submodules, which will be our key ingredient.
An important additional result will be borrowed from Plamondon~\cite{P}.

\subsubsection{$F$-polynomials and $g$-vectors of cluster variables}

Recall the cluster algebra $\AA$ of \S\ref{subsect_ca}, with initial
seed $(\bz^-,G^-)$. Following \cite[(3.7)]{FZ}, define
\begin{equation}\label{def_yhat}
\widehat{y}_{i,r} := \prod_{(i,r)\to (j,s)} z_{j,s} \prod_{(j,s)\to (i,r)} z_{j,s}^{-1},
\qquad ((i,r)\in W^-).
\end{equation}
Here the first (\resp second) product is over all outgoing (\resp incoming) arrows at the vertex $(i,r)$
of the graph $G^-$.
The following result is similar to \cite[Lemma 7.2]{HL}.
\begin{Lem}\label{lem_hat_y}
After performing in {\rm (\ref{def_yhat})} the change of variables {\rm (\ref{chvar})}, there holds
\[
\widehat{y}_{i,r} = A_{i,\, r-d_i}^{-1}, \qquad ((i,r)\in W^-),
\]
where the Laurent monomials $A_{i,r}$ are given by {\rm (\ref{def_A_var})}.
\end{Lem}
\proof
Using the definition of the quiver $G^-$, we can rewrite (\ref{def_yhat}) as 
\[
\widehat{y}_{i,r} = \frac{z_{i,\,r+b_{ii}}}{z_{i,\,r-b_{ii}}} 
\prod_{j\not = i} \frac{z_{j,\,r+b_{ij}+d_j-d_i}}{z_{j,\,r-b_{ij}+d_j-d_i}},
\]
where the product is over all $j$'s such that $c_{ij}\not = 0$.
Here we use the convention that $z_{i,s}=1$ for every $(i,s)$ with $s>0$.
Using the change of variables (\ref{chvar}), we obtain
\[
\widehat{y}_{i,r} = 
Y_{i,\,r-b_{ii}}^{-1}Y_{i,r}^{-1}
\prod_{j\not = i;\ c_{ij}\not = 0} Y_{r,\,r-d_i+b_{ij}+d_j} Y_{r,\,r-d_i+b_{ij}+3d_j}\cdots  Y_{r,\,r-d_i-b_{ij}-d_j}.
\]
The result then follows by comparison with (\ref{def_A_var}), if we notice again that
$b_{ij}+d_j = c_{ji}+1$ because of (\ref{easy-prop}).
\cqfd

In \cite{FZ} Fomin and Zelevinsky attach to every cluster variable $x$ of $\AA$
a polynomial $F_x$ with integer coefficients in the set of variables 
$\widehat{\by}=\{\widehat{y}_{i,r}\mid (i,r)\in W^-\}$,
and a vector $\bg_x \in \Z^{(W^-)}$, such that \cite[Corollary 6.3]{FZ}
\begin{equation}\label{form-Fg}
 x = \bz^{\bg_x}\, F_x(\widehat{\by}).
\end{equation}
Note that $\AA$ has no frozen cluster variables, so there is no denominator
in (\ref{form-Fg}).
The polynomial $F_x$ and the integer vector $g_x$ are called the \emph{$F$-polynomial} 
and \emph{$g$-vector} of the cluster variable $x$, respectively. We refer the
reader to \cite{FZ} for their definition.

On the other hand, it follows from the theory of $q$-characters that for every 
simple $U_q(\hg)$-module $L(m)$ in the category $\CC^-$, the truncated $q$-character
$\chi_q^-(L(m))$ can be written as
\begin{equation}\label{normalized}
\chi_q^-(L(m)) = m P_m, 
\end{equation}
where $P_m$ is a polynomial with integer coefficients in the variables 
$\{A_{i,\,r-d_i}^{-1}\mid (i,r)\in W^-\}$.
Moreover, $P_m$ has constant term $1$.

Now, by the proof of Theorem~\ref{thm1}, among the cluster variables of $\AA$, we find
all the truncated $q$-characters of the Kirillov-Reshetikhin modules of $\CC^-$.
These are of the form $L(m)$ with  
\begin{equation}
m = m_{k,r}^{(i)}:=\prod_{j=0}^{k-1} Y_{i,\,r+jb_{ii}},\qquad ((i,r)\in W^-,\ r+(k-1)b_{ii}\le 0).  
\end{equation}

\begin{Prop}\label{gvector}
The $g$-vector of the truncated $q$-character of the Kirillov-Reshetikhin module 
$W^{(i)}_{k,r}=L\left(m_{k,r}^{(i)}\right)$,
considered as a cluster variable of $\AA$, is given by
\[
g_{j,s} =
\left\{
\begin{array}{cl}
1 & \mbox{if $(j,s) = (i, r)$,}\\[1mm]
-1 & \mbox{if $(j,s) = (i, r+kb_{ii})$ and $r+kb_{ii}\le 0$,}\\[1mm]
0 & \mbox{otherwise.}
\end{array}
\right.
\]
\end{Prop}
\proof
Write for short $m = m_{k,r}^{(i)}$, and denote by $x$ the cluster variable 
$\chi_q^-(L(m))$.
Then, comparing (\ref{form-Fg}) with (\ref{normalized}), we have
\[
P_m = m^{-1} \bz^{\bg_x} F_x, 
\]
where, by Lemma~\ref{lem_hat_y}, $P_m$ and $F_x$ are polynomials in the same 
variables 
\[
\widehat{y}_{i,r} = A^{-1}_{i, r-d_i}.
\]
Since $P_m$ has constant term 1, it follows that $ m \bz^{-\bg_x}$
is a monomial in the variables $\widehat{y}_{i,r}$ which divides
the $F$-polynomial $F_x$. But, by \cite[Proposition 5.2]{FZ}, $F_x$
is not divisible by any $\widehat{y}_{i,r}$. So, using (\ref{chvar}), 
\[
\bz^{\bg_x} = m = \frac{z_{i,r}}{z_{i,\,r+kb_{ii}}}, 
\]
where as above, we set $z_{i,s} = 1$ if $s>0$. 
\cqfd

\subsubsection{Truncated algebras}

Let $\ell \in \Z_{<0}$. Let $\G_\ell^{-}$ be the full subquiver of
$\G^{-}$ with set of vertices  
\[
V_\ell^{-}:=\{(i,m)\in V^- \mid m\ge \ell\}.
\]
Let $S_\ell$ be the corresponding truncation of the potential $S$, that is, $S_\ell$
is defined as the sum of all cycles in $S$ which only involve vertices of $V_\ell^{-}$.
Let $J_\ell$ denote the two-sided ideal of $\C\G_\ell^{-}$ generated
by all cyclic derivatives of~$S_\ell$. Finally, define the \emph{truncated algebra
at height $\ell$} as
\[
A_\ell := \C\G_\ell^{-} / J_\ell. 
\]

\begin{Prop} \label{finite+rigid}
For every $\ell$ we have:
\begin{itemize}
 \item[{\rm (i)}] the algebra $A_\ell$ is finite-dimensional;
 \item[{\rm (ii)}] the quiver with potential $(\G_\ell^{-},J_\ell)$ is rigid.
\end{itemize}
\end{Prop}

\proof
The proof is similar to \cite[Example 8.7]{DWZ}.
Let $\pi : \C\G_\ell^{-} \to A_\ell$ be the natural projection.
To prove (i), we show that $A_\ell$ is spanned by the images under $\pi$
of a finite number of paths.
The arrows of $\G_\ell^{-}$ are of two types:
\begin{itemize}
 \item[(a)] the \emph{vertical} arrows of the form $(i,m) \to (i,m+b_{ii})$;
 \item[(b)] the \emph{oblique} arrows of the form $(i,m) \to (j,m+b_{ij})$
 provided $c_{ij} < 0$.
\end{itemize}
Let us say that a path from $(i,m)$ to $(j,s)$ in $\G_\ell^{-}$
is \emph{going up} (\resp \emph{down}) if $m<s$ (\resp $m>s$).
Note that all vertical arrows go up
and all oblique arrows go down.
Each oblique arrow of the boundary of $\G_\ell^{-}$ belongs to a single cycle of the potential $S_\ell$,
and each interior oblique arrow belongs to exactly two cycles.
Therefore each interior oblique arrow gives rise to a ``commutativity relation'' in $A_\ell$:
\[
\begin{array}{ll}
&\pi((j,m+b_{ji}),(i,m+2b_{ji}),(i,m+2b_{ji}+b_{ii}),\ldots,(i,m-b_{ii}),(i,m)) \\[2mm]
& = - \pi((j,m+b_{ji}),(j,m+b_{ji}+b_{jj}),\ldots,(j,m-b_{ji}-b_{jj}),(j,m-b_{ji}),(i,m)) 
\end{array}
\]
The path in the left-hand side consists of an oblique arrow followed by $|c_{ij}|$ vertical arrows,
while the right-hand side has $|c_{ji}|$ vertical arrows followed by an oblique arrow.
Let $p$ be a path in $\G_\ell^{-}$ with origin $(i,m)$.
Using only the above type of commutativity relations, we can bring a number
of vertical arrows to the front of $p$ and write 
\[
\pi(p) = \pi(p_2) \pi(p_1),
\]
where $p_1$ is a path with origin $(i,m)$ consisting only of vertical arrows, and
$p_2$ is a path satisfying the following property: if $q$ is a maximal factor of $p_2$
containing only vertical arrows, then $q$ is preceded by at least one oblique arrow,
say $(j,s)\to (k,s+b_{jk})$, and $q$ contains \emph{less} than $|c_{kj}|$ arrows.
Hence $q$ can be non trivial only if $|c_{kj}|>1$. 

In particular in the simply laced case, then $p_2$ contains only oblique arrows.
In that case, we can immediately conclude that all arrows of $p_1$ go up and
all arrows of $p_2$ go down, so the lengths of $p_1$ and $p_2$ are both 
bounded by $\ell$, and therefore $A_\ell$ is finite-dimensional. 

Otherwise, 
if $q$ is non trivial and $p_2$ contains other vertical arrows after $q$, 
then $q$ needs to be followed by at least \emph{two} oblique arrows. 
Indeed, using the same notation as above, $q$ consists of $N$ vertical arrows
of the form $(k, r) \to (k, r + b_{kk})$ with $1\le N<|c_{kj}|$.
Now, by (\ref{easy-prop}), 
the inequality
$|c_{kj}| > 1$ implies $d_k=1$ and $d_j=|b_{kj}|$.
Let $(k,t) \to (l, t + b_{kl})$ be the first arrow coming after~$q$.
Then, since $d_k=1$ we have $|c_{lk}|=1$. 
If this oblique arrow 
is followed by a vertical one
$(l, t + b_{kl}) \to (l, t +b_{kl} + b_{ll})$, then 
we can use the commutativity relation and bring it, together with all the vertical arrows possibly following it, on top
of $q$.
In this way, we replace $q$ by a vertical path $q'$ followed by two consecutive 
oblique arrows.

One then easily checks by inspection that the subpath of $p_2$ containing $q$
together with the oblique arrow preceding it and the oblique arrow following it, is going down. 
Therefore, by induction, $p_2$ can be factored into a product of paths, each of them of length 
less than $t+2$, and all these paths go down (except possibly the last one, which might end with
less than $t$ vertical arrows).
So again, the length of $p_2$ is bounded above, and this proves (i) in all cases.

To prove (ii), it is enough to show that every cycle of the form (\ref{cycle})
is cyclically equivalent to an element of $J_\ell$. Up to cyclic equivalence,
this cycle $\gamma$ can be written with origin in $(i,m)$. Then, we 
have:
\[
\begin{array}{ll}
\pi(\gamma)&=\ \pi((i,m), (j,m+b_{ij}), (i,m+2b_{ij}),  (i,m+2b_{ij}+b_{ii}),
\ldots, (i,m-b_{ii}),(i,m)) \\[2mm]
& = \ 
\pi((i,m), (j,m+b_{ij}), (j,m+b_{ij}+b_{jj}), \ldots, (j,m-b_{ij}-b_{jj}),
(j,m-b_{ij}),(i,m)) \\[2mm]
& = \ 
\pi((i,m), (i,m+b_{ii}), \ldots, (i,m-2b_{ij}-b_{ii}), (i,m-2b_{ij}),
(j,m-b_{ij}),(i,m)), 
\end{array}
\]
and the last path is cyclically equivalent to
\[
((i,m-2b_{ij}), (j,m-b_{ij}), (i,m),  (i,m+b_{ii}),
\ldots, (i,m-2b_{ij}-b_{ii}),(i,m-2b_{ij})). 
\]
This cycle is nothing else than $\gamma$ shifted vertically up by $-2b_{ij}$.
Hence, iterating this process, we can replace, modulo $J_\ell$ and cyclic equivalence, 
any cycle $\gamma$ of the form (\ref{cycle}) by 
a similar cycle $\gamma'$ sitting at the top boundary of $\G_\ell^{-}$.
Now the upper oblique arrow of $\gamma'$ does not belong to any other cycle,
so it gives rise to a zero relation in~$A_\ell$. In other words, 
$\gamma'$ is cyclically equivalent to an element of~$J_\ell$.
This proves (ii). 
\cqfd

\begin{remark}
{\rm
In the simply laced case and when
$|\ell|$ is less than the Coxeter number, the algebra $A_\ell$ arises as the
endomorphism algebra of a (finite-dimensional) rigid module over the preprojective algebra $\L$
associated with $\de$, and appears in the works of Geiss, Schr\"oer and the 
second author (see \cite{GLS1,GLS2}). 
This gives another proof of Proposition~\ref{finite+rigid}~(i)
in this case.
}
\end{remark}

\subsubsection{Proof of Lemma~\ref{finite-dim} and Theorem~\ref{th_geom_form}}\label{proofs}

Let $(i,r)\in V^-$ and $k\in \N$.
By Theorem~\ref{thm1}, the truncated $q$-character $\chi_q^-\left(W^{(i)}_{k,\ r-(2k-1)d_i}\right)$
is a cluster variable $x$ of $\AA$.
By Proposition~\ref{gvector}, the $g$-vector of $x$ is given by 
\begin{equation}\label{gvector2}
g_{j,s} =
\left\{
\begin{array}{cl}
1 & \mbox{if $(j,s) = (i, r-2kd_i+d_i)$,}\\[1mm]
-1 & \mbox{if $(j,s) = (i, r+d_i)$,}\\[1mm]
0 & \mbox{otherwise.}
\end{array}
\right.
\end{equation}
Note that, since $(i,r)\in V^-$, we have $(i, r+d_i)\in W^-$.
For $\ell < 0$, let $W_\ell^{-} := \psi(V_\ell^{-})$, and put
$\bz_\ell^{-} = \{z_{i,r}\mid (i,r)\in W_\ell^{-}\}$. 
We denote by $G_\ell^{-}$ the same quiver as $\G_\ell^{-}$,
but with vertices labelled by $W_\ell^{-}$.
Clearly, the cluster variable $x$ is a Laurent polynomial in the 
variables of $\bz_\ell^{-}$ for some $\ell \ll 0$, and can be 
regarded as a cluster variable of the cluster algebra $\AA_\ell$
defined by the initial seed $\left(\bz_\ell^{-}, G_\ell^{-}\right)$.
By Proposition~\ref{finite+rigid}~(ii), we can apply the theory of \cite{DWZ,DWZ2} 
and deduce that the $F$-polynomial of $x$ coincides with the polynomial
$F_M$ associated with a certain $A_\ell$-module~$M$.
In order to identify this module, we apply \cite[Remark 4.1]{P}, which
states in our setting that $M$ is the kernel of a generic element of
the homomorphism space
between two injective $A_\ell$-modules corresponding to the negative and
positive components of the $g$-vector of $x$. More precisely,  
let us denote by $S_{i,m}^{\ell}$ the one-dimensional $A_\ell$-module
supported on $(i,m)\in V_\ell^{-}$.
Let $I^\ell_{i,m}$ be the injective 
$A_\ell$-module with socle isomorphic to $S^\ell_{i,m}$.
Then, using (\ref{gvector2}) and taking into account the change of labelling 
$\psi\colon V_\ell^{-} \to W_\ell^{-}$ given by (\ref{def_psi}), we get that $M$ is the 
kernel of a generic element of 
$\Hom_{A_\ell}(I^\ell_{i,r}, I^\ell_{i,r-kb_{ii}})$.

Finally we can identify $M$ with the kernel of a generic homomorphism
between injective $A$-modules. Indeed, for $m < \ell <0$ we have a natural
projection $A_m \to A_\ell$ whose kernel is generated by all arrows of $\G_m^-$
starting or ending at a vertex $v \in V_m^-\setminus V_\ell^-$.
This induces for every $(i,r)\in V^-_\ell$ an embedding $I_{(i,r)}^\ell \to I_{(i,r)}^m$,
and we can regard the $A$-module $I_{(i,r)}$ as the direct limit of $I_{(i,r)}^\ell$
along these maps.
Since $F_M$ is independent of $\ell \ll 0$, we see that
$M$ is also the kernel of a generic element of $\Hom_{A}(I_{i,r}, I_{i,r-kb_{ii}})$,
that is, 
$M = K^{(i)}_{k,r}$. In particular $K^{(i)}_{k,r}$ is finite-dimensional.
This proves Lemma~\ref{finite-dim} and
finishes the proof of Theorem~\ref{th_geom_form}.

\begin{remark}
{\rm
Using the same formula as (\ref{FpolAmod}), we can attach
to the infinite-dimensional $A$-module $I_{i,m}$ a formal 
power series $F_{I_{i,m}}$ in the variables $v_{j,r}$. 
This series also has an interpretation in terms of quantum affine algebras. 
Indeed, by \cite{HJ}, the category of finite-dimensional 
$U_q(\hg)$-modules can be seen as a subcategory of a category $\mathcal{O}$ 
of (possibly infinite-dimensional) representations of a Borel subalgebra of $U_q(\hg)$. 
The $q$-character morphism can be extended to the Grothendieck ring of  
$\mathcal{O}$ (the target ring is also completed). 
This category contains distinguished simple representations called negative 
fundamental representations $L_{{i},a}^-$ (${i}\in {I}$, $a\in\mathbb{C}^*$)
\cite[Definition 3.7]{HJ}.
Denote by $\widetilde{\chi}_q(L_{i,a}^-)$ the normalized $q$-character of $L_{i,a}^-$, 
that is, its $q$-character divided by its highest weight monomial.
This normalized $q$-character is a formal power series in the variables $A_{{j},b}^{-1}$ 
\cite[Theorem 6.1]{HJ}, and it
is obtained as a limit of normalized $q$-characters of Kirillov-Reshetikhin modules.
It is not difficult to deduce from Theorem~\ref{th_geom_form} and Remark~\ref{complete-character} 
that, for $m\le d_i-t\dh$,
\[
\widetilde{\chi}_q\left(L_{i,\,q^{m-d_i}}^-\right) = F_{I_{i,m}}.
\]
This is the first geometric description of the $q$-character
of these negative fundamental representations. 
}
\end{remark}


\section{Beyond Kirillov-Reshetikhin modules}\label{section4}

\subsection{Grothendieck rings}

Let us consider again the cluster algebra $\AA$, with initial seed
$\Sigma = (\bz^-,G^-)$ whose 
cluster variables $z_{i,r}$ are given by (\ref{chvar}).
The Laurent phenomenon for cluster algebras implies that $\AA$ is 
a subring of $\Z[Y_{i,r}^{\pm1}\mid Y_{i,r}\in \bY^-]$.
The following theorem gives the precise relationship between  
$\AA$ and the Grothendieck ring of the category $\CC^-$.

\begin{Thm}\label{thm_Grothendieck_ring}
The cluster algebra $\AA$ is equal to
the image of the injective ring homomorphism
from $K_0(\CC^-)$ to $\Z[Y_{i,r}^{\pm1}\mid Y_{i,r}\in \bY^-]$
given by $[L(m)]\mapsto \chi_q^-(m)$ (see  Proposition~\ref{prop_truncated}).
Hence $\AA$ is isomorphic to the Grothendieck ring of $\CC^-$.
\end{Thm}
\proof
Let $R^-$ denote the image of the homomorphism $[L(m)]\mapsto \chi_q^-(m)$.
By \cite{FR}, $K_0(\CC^-)$ is the polynomial ring in the classes of the fundamental
modules of $\CC^-$, hence $R^-$ is the polynomial ring in the truncated $q$-characters
$\chi_q^-(Y_{i,r})\ (Y_{i,r}\in \bY^-)$.
By Theorem~\ref{thm1}, $\AA$ contains all these fundamental truncated 
$q$-characters, hence $\AA$ contains $R^-$. 

To prove the reverse inclusion, we will use a description of the image of
the $q$-character homomorphism as an intersection of kernels of screening
operators \cite{FR,FM}. To do this, we need to work with complete
(\ie untruncated) $q$-characters. So let us consider as in \S\ref{ssect-truncated}
the larger set of variables $\bY$. Following \cite[\S7.1]{FR}, for every $i\in I$,
we have a linear operator $S_i$ from the ring
$\Z[Y_{i,r}^{\pm1}\mid Y_{i,r} \in \bY]$ to a certain free module
$\Y_i$ over this ring, which satisfies the Leibniz rule
\[
S_i(xy) =  x\,S_i(y) + y\,S_i(x),\qquad (x,y \in \Z[Y_{i,r}^{\pm1}\mid Y_{i,r} \in \bY]). 
\]
It was conjectured in \cite{FR} and proved in \cite{FM} that 
an element of $\Z[Y_{i,r}^{\pm1}\mid Y_{i,r} \in \bY]$ is a polynomial
in the $q$-characters $\chi_q(Y_{i,r})\ (Y_{i,r}\in \bY)$ if and only if
it belongs to 
\[
\bigcap_{i\in I} \Ker S_i. 
\]

Let us now introduce an auxilliary cluster algebra $\AA'$. 
It is defined using the same initial seed $(\bz^-, G^-)$ as $\AA$, 
but the initial variables of $\AA'$ are given by the following modification 
of (\ref{chvar})
\[
z'_{i,r} := \prod_{k\ge 0,\ r+kb_{ii}\le 0} Y_{i,\,r+kb_{ii}+2t\dh}, 
\]
in which the spectral parameters are all shifted upwards by $2t\dh$.
By Theorem~\ref{thm1}, if we apply to this initial seed of $\AA'$ the 
sequence of mutations $\mu_\SS$ repeated $\dh$ times, we will obtain 
a new seed $\Sigma'$ with the same quiver $ G^-$. 
Moreover, the cluster variable of $\Si'$ sitting at vertex $(i,r)\in W^-$
is nothing else than the \emph{complete} $q$-character
$\chi_q(W^{(i)}_{k_{i,r},r})$.

Consider a cluster variable $x$ of $\AA$. By definition, $x$ is obtained 
from $\Sigma$ by a finite sequence of mutations $\mu_x$. We want to show that 
$x$ belongs to $R^-$. By Theorem~\ref{thm1}, all cluster variables
of $\Si$ belong to $R^-$, so by induction on the length,
we may assume that the last exchange relation of $\mu_x$ is of the form
\[
xy = M_1 + M_2, 
\]
where $y$ is a cluster variable of $\AA$, $M_1$ and $M_2$ are cluster
monomials of $\AA$, and $y$, $M_1$, $M_2$ belong to $R^-$.
Let us apply the same sequence of mutations $\mu_x$ in the cluster
algebra $\AA'$ to the seed $\Si'$. The last exchange relation will
be of the form 
\[
x'y' = M_1' + M_2', 
\]
where $y'$, $M_1'$, $M_2'$ are polynomials in the complete fundamental
$q$-characters  $\chi_q(Y_{i,r})\ (Y_{i,r}\in \bY^-)$.
Moreover, $x'$, $y'$, $M_1'$, $M_2'$ give back $x$, $y$, $M_1$, $M_2$
by application of the truncation ring homomorphism.
By the Laurent phenomenon \cite{FZ1} in the cluster algebra $\AA'$, we know
that $x'$, $y'$, $M_1'$, $M_2'$ are Laurent polynomials in the variables of $\bY$.
Since $S_i$ is a derivation, we have
\[
S_i(x'y') = x'S_i(y') + y'S_i(x')= S_i(M_1') + S_i(M_2'),
\]
hence $S_i(x') = 0$ because $S_i(y') = S_i(M_1') = S_i(M_2') =0$.
It follows that $x'$ is annihilated by all the screening operators, 
so $x'$ is a polynomial in the
$q$-characters  $\chi_q(Y_{i,r})\ (Y_{i,r}\in \bY^-)$.
This implies that $x$ is  a polynomial in the truncated
$q$-characters  $\chi_q^-(Y_{i,r})\ (Y_{i,r}\in \bY^-)$,
that is, $x\in R^-$. \cqfd

\subsection{Conjectures}

\subsubsection{Cluster monomials}
In view of Theorem~\ref{thm_Grothendieck_ring}, it is natural to formulate
some conjectures. Following \cite{Le-imag}, let us say that a simple $U_q(\hg)$-module $S$ 
is \emph{real} if $S\otimes S$ is simple.

\begin{Conj}\label{conj1}
In the above identification of the cluster algebra $\AA$ with 
the ring of truncated $q$-characters of $\CC^-$, the cluster
monomials get identified with the truncated $q$-characters of
the real simple modules of $\CC^-$.
\end{Conj}

When $\g$ is of type $A$, $D$, $E$, Conjecture~\ref{conj1} is essentially
equivalent to \cite[Conjecture 13.2]{HL}. But the initial seed used here
is different and allows a direct connection between cluster expansions
and (truncated) $q$-characters. 

\subsubsection{Geometric $q$-character formulas}
Using the methods and tools of \S\ref{section3}, we can translate 
Conjecture~\ref{conj1} into a new conjectural geometric formula
for the (truncated) $q$-character of a real simple module of $\CC^-$.

Let $m$ be a dominant monomial in the variables  $Y_{i,r}\in\bY^-$.
Using the change of variables (\ref{chvar}), which we can express as
\[
Y_{i,r} = \frac{z_{i,r}}{z_{i,r+b_{ii}}}, \qquad ((i,r)\in W^-),
\]
(where we understand $z_{i,s} = 1$ if $s>0$), we can rewrite
\[
m = \bz^{g(m)} := \prod_{(i,r)\in W^-}  z_{i,r}^{g_{i,r}(m)}.
\]
Let us call the integer vector $g(m)\in \Z^{(W^-)}$ the \emph{$g$-vector
of $L(m)$}. 
Following \S\ref{ssect-A-mod}, let us attach to $m$ the $A$-module
$K(m)$ defined as the kernel of a generic $A$-module homomorphism
from the injective $A$-module $I(m)^-$ to the injective $A$-module $I(m)^+$,
where 
\[
I(m)^+ = \bigoplus_{g_{i,r}(m)>0} I_{i,r-d_i}^{\oplus g_{i,r}(m)},
\qquad
I(m)^- = \bigoplus_{g_{i,r}(m)<0} I_{i,r-d_i}^{\oplus |g_{i,r}(m)|}.
\]
Finally define the $F$-polynomial $F_{K(m)}$ of $K(m)$ as in \S\ref{ssect-F-pol-A-mod}.
We can now state the following conjectural generalization of Theorem~\ref{th_geom_form}.
\begin{Conj}\label{Conj2}
Suppose that $L(m)$ is an irreducible real $U_q(\hg)$-module in $\CC^-$.
Then the truncated $q$-character of $L(m)$ is equal to
\[
\chi_q^-(L(m)) = m F_{K(m)}, 
\]
where the variables $v_{i,r}$ of the $F$-polynomial are evaluated as in {\rm(\ref{eval_v})}.
\end{Conj}

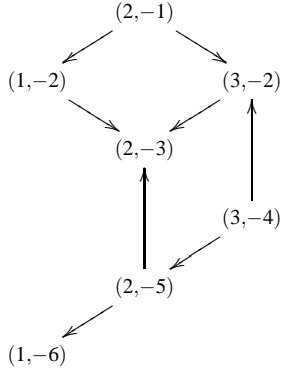
\begin{figure}[t]
\[
\def\objectstyle{\scriptstyle}
\def\lablestyle{\scriptstyle}
\xymatrix@-1.0pc{
&& \ar[ld](2,-1)\ar[rd]&&
\\
&(1,-2)\ar[rd] &&\ar[ld] (3,-2) 
\\
&& (2,-3) &&
\\
&&&\ar[ld] (3,-4)\ar[uu] 
\\
&&\ar[ld] (2,-5)\ar[uu] &&
\\
&(1,-6)  && 
\\
}
\]
\caption{\label{Fig8} {\it The $A$-module $K(m)$ for $m=Y_{1,-7}Y_{2,-4}$ in type $A_3$.}}
\end{figure}

\begin{example}
{\rm
Let $\g$ be of type $A_3$. Take $m=Y_{1,-7}Y_{2,-4}$.
We have
\[
I(m)^+ = I_{1,-8}\oplus I_{2,-5},\qquad I(m)^- =  I_{1,-6}\oplus I_{2,-3}. 
\]
The module $K(m)$ has dimension 7 and is displayed in Figure~\ref{Fig8}.
Using for instance the fact that $L(m)$ is a minimal affinization (in the 
sense of \cite{C}), we can compute its $q$-character. We find:
\[
\begin{array}{rcl}
\chi_q(L(Y_{1,-7}Y_{2,-4})) &=& 
Y_{1,-7}Y_{2,-4}\left(1 + v_{1,-6}+ v_{2,-3}+ v_{1,-6}v_{2,-3}
+ v_{1,-2}v_{2,-3} + v_{2,-3}v_{3,-2}
\right.\\[2mm]
&&\left.
+\ v_{1,-6}v_{1,-2}v_{2,-3} 
+ v_{1,-6}v_{2,-3}v_{3,-2}
+ v_{1,-6}v_{2,-3}v_{2,-5}
+ v_{1,-2}v_{2,-3}v_{3,-2}
\right.\\[2mm]
&&\left. 
+\ v_{1,-6}v_{1,-2}v_{2,-5}v_{2,-3}
+ v_{1,-6}v_{1,-2}v_{2,-3}v_{3,-2}
+ v_{1,-6}v_{2,-5}v_{2,-3}v_{3,-2}
\right.\\[2mm]
&&\left. 
+\ v_{1,-2}v_{2,-5}v_{2,-3}v_{3,-2} 
+ v_{1,-6}v_{1,-2}v_{2,-5}v_{2,-3}v_{3,-2} 
\right.\\[2mm]
&&\left. 
+\ v_{1,-6}v_{1,-2}v_{2,-3}v_{2,-1}v_{3,-2}
+ v_{1,-6}v_{2,-5}v_{2,-3}v_{3,-4}v_{3,-2}
\right.\\[2mm]
&&\left. 
+\ v_{1,-6}v_{1,-2}v_{2,-5}v_{2,-3}v_{3,-4}v_{3,-2}
+ v_{1,-6}v_{1,-2}v_{2,-5}v_{2,-3}v_{2,-1}v_{3,-2}  
\right.\\[2mm]
&&\left. 
+\ v_{1,-6}v_{1,-2}v_{2,-5}v_{2,-3}v_{2,-1}v_{3,-4}v_{3,-2} 
\right),
\end{array}
\]
in agreement with Conjecture~\ref{Conj2}.
}
\end{example}


\section{Appendix}

\subsection{Mutation sequence in type $A_2$}\label{appendixA3}

We display the sequence of mutated quivers obtained from $G^-$ 
at each step of the mutation sequence $\mu_\SS$. The first quiver
is $G^-$, and in the next quivers the box indicates  
at which vertex a mutation has been performed.

\[
\def\objectstyle{\scriptscriptstyle}
\xymatrix@-1.0pc{
&\ar[ld](2,0)\\
{(1,-1)}\ar[rd]& 
\\
&\ar[ld] (2,-2)\ar[uu] 
\\
{(1,-3)}\ar[rd]\ar[uu]& 
\\
&\ar[ld] (2,-4) \ar[uu]
\\
\ar[uu]{(1,-5)}\ar[rd]&
\\
&\ar[ld]\ar[uu] (2,-6) 
\\
\ar[uu]{(1,-7)}
&{}\save[]+<0cm,2ex>*{\vdots}\restore
}
\hskip1.2cm
\xymatrix@-1.0pc{
&\ar[dd]\fbox{$\scriptscriptstyle(2,0)$}
\\
{(1,-1)}\ar[ru]& 
\\
&\ar[ld] (2,-2) 
\\
{(1,-3)}\ar[rd]\ar[uu]& 
\\
&\ar[ld] (2,-4) \ar[uu]
\\
\ar[uu]{(1,-5)}\ar[rd]&
\\
&\ar[ld]\ar[uu] (2,-6) 
\\
\ar[uu]{(1,-7)}
&{}\save[]+<0cm,2ex>*{\vdots}\restore&
}
\hskip1cm
\xymatrix@-1.0pc{
&\ar[lddd](2,0)
\\
{(1,-1)}\ar[ru]& 
\\
&\ar[uu]\ar[dd] \fbox{$\scriptscriptstyle(2,-2)$}
\\
{(1,-3)}\ar[ru]\ar[uu]& 
\\
&\ar[ld] (2,-4)
\\
\ar[uu]{(1,-5)}\ar[rd]&
\\
&\ar[ld]\ar[uu] (2,-6) 
\\
\ar[uu]{(1,-7)}
&{}\save[]+<0cm,2ex>*{\vdots}\restore&
}
\hskip0.8cm
\xymatrix@-1.0pc{
&\ar[lddd](2,0)
\\
{(1,-1)}\ar[ru]& 
\\
&\ar[uu]\ar[lddd] (2,-2)
\\
{(1,-3)}\ar[ru]\ar[uu]& 
\\
&\ar[uu] \fbox{$\scriptscriptstyle(2,-4)$} \ar[dd]
\\
\ar[uu]{(1,-5)}\ar[ru]&
\\
&\ar[ld] (2,-6) 
\\
\ar[uu]{(1,-7)}
&{}\save[]+<0cm,2ex>*{\vdots}\restore&
}
\]
\[
\def\objectstyle{\scriptscriptstyle}
\xymatrix@-1.0pc{
&\ar[lddd](2,0)
\\
{(1,-1)}\ar[ru]&
\\
&\ar[uu]\ar[lddd] (2,-2)
\\
{(1,-3)}\ar[ru]\ar[uu]& 
\\
&\ar[uu] (2,-4)\ar[lddd]
\\
\ar[uu]{(1,-5)}\ar[ru]&
\\
&\ar[uu] \fbox{$\scriptscriptstyle(2,-6)$} 
\\
\ar[uu]{(1,-7)}\ar[ru]
&{}\save[]+<0cm,2ex>*{\vdots}\restore&
}
\xymatrix@-1.0pc{
&\ar[ld](2,0)
\\
\fbox{$\scriptscriptstyle(1,-1)$}\ar[dd] 
\\
&\ar[uu]\ar[lddd] (2,-2)&
\\
{(1,-3)}\ar[ru]& 
\\
&\ar[uu] (2,-4)\ar[lddd]
\\
\ar[uu]{(1,-5)}\ar[ru]&
\\
&\ar[uu] (2,-6)  
\\
\ar[uu]{(1,-7)}\ar[ru]
&{}\save[]+<0cm,2ex>*{\vdots}\restore&
}
\def\objectstyle{\scriptscriptstyle}
\xymatrix@-1.0pc{
&\ar[ld](2,0)
\\
(1,-1)\ar[rd]& 
\\
&\ar[uu]\ar[ld] (2,-2)
\\
\fbox{$\scriptscriptstyle(1,-3)$}\ar[uu]\ar[dd]& 
\\
&\ar[uu] (2,-4)\ar[lddd]
\\
{(1,-5)}\ar[ru]&
\\
&\ar[uu] (2,-6) 
\\
\ar[uu]{(1,-7)}\ar[ru]
&{}\save[]+<0cm,2ex>*{\vdots}\restore&
}
\xymatrix@-1.0pc{
&\ar[ld](2,0)&
\\
(1,-1)\ar[rd]&& 
\\
&\ar[uu]\ar[ld] (2,-2)&&
\\
(1,-3)\ar[uu]\ar[rd]&&
\\
&\ar[uu] (2,-4)\ar[ld]&&
\\
\fbox{$\scriptscriptstyle(1,-5)$}\ar[uu]\ar[dd]&
\\
&\ar[uu] (2,-6)  &&
\\
{(1,-7)}\ar[ru]
&{}\save[]+<0cm,2ex>*{\vdots}\restore&
}
\xymatrix@-1.0pc{
&\ar[ld](2,0)
\\
(1,-1)\ar[rd]& 
\\
&\ar[uu]\ar[ld] (2,-2)
\\
(1,-3)\ar[uu]\ar[rd]& 
\\
&\ar[uu] (2,-4)\ar[ld]
\\
(1,-5)\ar[uu]\ar[rd]&
\\
&\ar[uu] (2,-6)\ar[ld] 
\\
\fbox{$\scriptscriptstyle(1,-7)$}\ar[uu]
&{}\save[]+<0cm,2ex>*{\vdots}\restore&
} 
\]

\subsection{Mutation sequence in type $B_2$}\label{appendixB2}

We display the sequence of mutated quivers obtained from $G^-$ 
at each step of the mutation sequence $\mu_\SS$. 

\[
\def\objectstyle{\scriptscriptstyle}
\xymatrix@-1.0pc{
&&&&\\
&&(\bdeux,0)\ar[rd]
& 
\\
&&\ar[ld]\ar[u] (\bdeux,-2)& \ar[ld](\bun,-1)&
\\
&{(\bun,-3)}\ar[rd] &\ar[u] (\bdeux,-4) \ar[rd]&&
\\
&&\ar[u] \ar[ld](\bdeux,-6) &\ar[ld] (\bun,-5) \ar[uu]  &
\\
&(\bun,-7)\ar[uu]\ar[rd]    & \ar[u](\bdeux,-8)\ar[rd] &&
\\
&&(\bdeux, -10)\ar[u]\ar[ld] & (\bun,-9)\ar[uu] &
\\
&(\bun,-11)\ar[uu]& {}\save[]+<0cm,0ex>*{\vdots}\restore& {}\save[]+<0cm,0ex>*{\vdots}\restore&
}
\xymatrix@-1.0pc{
&&&&\\
&&\fbox{$\scriptscriptstyle(\bdeux,0)$}\ar[d]
& 
\\
&&\ar[ld]\ar[r] (\bdeux,-2)& \ar[ld]\ar[lu](\bun,-1)&
\\
&{(\bun,-3)}\ar[rd] &\ar[u] (\bdeux,-4) \ar[rd]&&
\\
&&\ar[u] \ar[ld](\bdeux,-6) &\ar[ld] (\bun,-5) \ar[uu]  &
\\
&(\bun,-7)\ar[uu]\ar[rd]    & \ar[u](\bdeux,-8)\ar[rd] &&
\\
&&(\bdeux, -10)\ar[u]\ar[ld] & (\bun,-9)\ar[uu] &
\\
&(\bun,-11)\ar[uu]& {}\save[]+<0cm,0ex>*{\vdots}\restore& {}\save[]+<0cm,0ex>*{\vdots}\restore&
}
\xymatrix@-1.0pc{
&&&&\\
&&\ar[ldd](\bdeux,0)
& 
\\
&&\ar[u]\ar[d] \fbox{$\scriptscriptstyle(\bdeux,-2)$}& \ar[l](\bun,-1)&
\\
&{(\bun,-3)}\ar[rd]\ar[ru] & \ar[l](\bdeux,-4) \ar[rd]&&
\\
&&\ar[u] \ar[ld](\bdeux,-6) &\ar[ld] (\bun,-5) \ar[uu]  &
\\
&(\bun,-7)\ar[uu]\ar[rd]    & \ar[u](\bdeux,-8)\ar[rd] &&
\\
&&(\bdeux, -10)\ar[u]\ar[ld] & (\bun,-9)\ar[uu] &
\\
&(\bun,-11)\ar[uu]& {}\save[]+<0cm,0ex>*{\vdots}\restore& {}\save[]+<0cm,0ex>*{\vdots}\restore&
}
\]
\[
\def\objectstyle{\scriptscriptstyle}
\xymatrix@-1.0pc{
&&&&\\
&&\ar[ldd](\bdeux,0)
& 
\\
&&\ar[u] (\bdeux,-2)\ar[rdd]& \ar[l](\bun,-1)&
\\
&{(\bun,-3)}\ar[r] & \ar[u]\ar[d]\fbox{$\scriptscriptstyle(\bdeux,-4)$} &&
\\
&&\ar[ld](\bdeux,-6)\ar[r] &\ar[ld]\ar[lu] (\bun,-5) \ar[uu]  &
\\
&(\bun,-7)\ar[uu]\ar[rd]    & \ar[u](\bdeux,-8)\ar[rd] &&
\\
&&(\bdeux, -10)\ar[u]\ar[ld] & (\bun,-9)\ar[uu] &
\\
&(\bun,-11)\ar[uu]& {}\save[]+<0cm,0ex>*{\vdots}\restore& {}\save[]+<0cm,0ex>*{\vdots}\restore&
}
\xymatrix@-1.0pc{
&&&&\\
&&\ar[ldd](\bdeux,0)
& 
\\
&&\ar[u] (\bdeux,-2)\ar[rdd]& \ar[l](\bun,-1)&
\\
&{(\bun,-3)}\ar[r] & \ar[ldd]\ar[u](\bdeux,-4) &&
\\
&&\ar[u]\ar[d]\fbox{$\scriptscriptstyle(\bdeux,-6)$} &\ar[l] (\bun,-5) \ar[uu]  &
\\
&(\bun,-7)\ar[ru]\ar[uu]\ar[rd]    & \ar[l](\bdeux,-8)\ar[rd] &&
\\
&&(\bdeux, -10)\ar[u]\ar[ld] & (\bun,-9)\ar[uu] &
\\
&(\bun,-11)\ar[uu]& {}\save[]+<0cm,0ex>*{\vdots}\restore& {}\save[]+<0cm,0ex>*{\vdots}\restore&
}
\xymatrix@-1.0pc{
&&&&\\
&&\ar[ldd](\bdeux,0)
& 
\\
&&\ar[u] (\bdeux,-2)\ar[rdd]& \ar[l](\bun,-1)&
\\
&{(\bun,-3)}\ar[r] & \ar[ldd]\ar[u](\bdeux,-4) &&
\\
&&\ar[u](\bdeux,-6)\ar[rdd] &\ar[l] (\bun,-5) \ar[uu]  &
\\
&(\bun,-7)\ar[uu]\ar[r]    & \ar[u]\ar[d]\fbox{$\scriptscriptstyle(\bdeux,-8)$} &&
\\
&&(\bdeux, -10)\ar[ld]\ar[r] & \ar[lu](\bun,-9)\ar[uu] &
\\
&(\bun,-11)\ar[uu]& {}\save[]+<0cm,0ex>*{\vdots}\restore& {}\save[]+<0cm,0ex>*{\vdots}\restore&
}
\]
\[
\def\objectstyle{\scriptscriptstyle}
\xymatrix@-1.0pc{
&&&&\\
&&\ar[ldd](\bdeux,0)
& 
\\
&&\ar[u] (\bdeux,-2)\ar[rdd]& \ar[l](\bun,-1)&
\\
&{(\bun,-3)}\ar[r] & \ar[ldd]\ar[u](\bdeux,-4) &&
\\
&&\ar[u](\bdeux,-6)\ar[rdd] &\ar[l] (\bun,-5) \ar[uu]  &
\\
&(\bun,-7)\ar[uu]\ar[r]    &\ar[ldd] \ar[u](\bdeux,-8) &&
\\
&&\fbox{$\scriptscriptstyle(\bdeux, -10)$}\ar[u]& \ar[l](\bun,-9)\ar[uu] &
\\
&(\bun,-11)\ar[uu]\ar[ru]& {}\save[]+<0cm,0ex>*{\vdots}\restore& {}\save[]+<0cm,0ex>*{\vdots}\restore&
}
\xymatrix@-1.0pc{
&&&&\\
&&\ar[ldd](\bdeux,0)
& 
\\
&&\ar[u] (\bdeux,-2)\ar[r]& \ar[dd]\fbox{$\scriptscriptstyle(\bun,-1)$}&
\\
&{(\bun,-3)}\ar[r] & \ar[ldd]\ar[u](\bdeux,-4) &&
\\
&&\ar[u](\bdeux,-6)\ar[rdd] &\ar[l] (\bun,-5)   &
\\
&(\bun,-7)\ar[uu]\ar[r]    &\ar[ldd] \ar[u](\bdeux,-8) &&
\\
&&(\bdeux, -10)\ar[u]& \ar[l](\bun,-9)\ar[uu] &
\\
&(\bun,-11)\ar[uu]& {}\save[]+<0cm,0ex>*{\vdots}\restore& {}\save[]+<0cm,0ex>*{\vdots}\restore&
}
\xymatrix@-1.0pc{
&&&&\\
&&\ar[ldd](\bdeux,0)
& 
\\
&&\ar[u] (\bdeux,-2)\ar[r]&\ar[ldd] (\bun,-1)&
\\
&{(\bun,-3)}\ar[r] & \ar[ldd]\ar[u](\bdeux,-4) &&
\\
&&\ar[u](\bdeux,-6)\ar[r] &\ar[uu]\ar[dd] \fbox{$\scriptscriptstyle(\bun,-5)$}   &
\\
&(\bun,-7)\ar[uu]\ar[r]    &\ar[ldd] \ar[u](\bdeux,-8) &&
\\
&&(\bdeux, -10)\ar[u]& \ar[l](\bun,-9) &
\\
&(\bun,-11)\ar[uu]& {}\save[]+<0cm,0ex>*{\vdots}\restore& {}\save[]+<0cm,0ex>*{\vdots}\restore&
}
\]

\[
\def\objectstyle{\scriptscriptstyle}
\xymatrix@-1.0pc{
&&&&\\
&&\ar[ldd](\bdeux,0)
& 
\\
&&\ar[u] (\bdeux,-2)\ar[r]&\ar[ldd] (\bun,-1)&
\\
&{(\bun,-3)}\ar[r] & \ar[ldd]\ar[u](\bdeux,-4) &&
\\
&&\ar[u](\bdeux,-6)\ar[r] &\ar[uu]\ar[ldd] (\bun,-5)   &
\\
&(\bun,-7)\ar[uu]\ar[r]    &\ar[ldd] \ar[u](\bdeux,-8) &&
\\
&&(\bdeux, -10)\ar[u] \ar[r]&\fbox{$\scriptscriptstyle(\bun,-9)$}\ar[uu] &
\\
&(\bun,-11)\ar[uu]& {}\save[]+<0cm,0ex>*{\vdots}\restore& {}\save[]+<0cm,0ex>*{\vdots}\restore&
}
\xymatrix@-1.0pc{
&&&&\\
&&\fbox{$\scriptscriptstyle(\bdeux,0)$}\ar[d]
& 
\\
&&\ar[ld] (\bdeux,-2)\ar[r]&\ar[ldd] (\bun,-1)&
\\
&{(\bun,-3)}\ar[r]\ar[ruu] & \ar[ldd]\ar[u](\bdeux,-4) &&
\\
&&\ar[u](\bdeux,-6)\ar[r] &\ar[uu]\ar[ldd] (\bun,-5)   &
\\
&(\bun,-7)\ar[uu]\ar[r]    &\ar[ldd] \ar[u](\bdeux,-8) &&
\\
&&(\bdeux, -10)\ar[u] \ar[r]&(\bun,-9)\ar[uu] &
\\
&(\bun,-11)\ar[uu]& {}\save[]+<0cm,0ex>*{\vdots}\restore& {}\save[]+<0cm,0ex>*{\vdots}\restore&
}
\xymatrix@-1.0pc{
&&&&\\
&&(\bdeux,0)\ar[rd]
& 
\\
&& \ar[u]\ar[d]\fbox{$\scriptscriptstyle(\bdeux,-2)$}&\ar[l]\ar[ldd] (\bun,-1)&
\\
&{(\bun,-3)}\ar[ru] & \ar[ldd](\bdeux,-4)\ar[ru] &&
\\
&&\ar[u](\bdeux,-6)\ar[r] &\ar[uu]\ar[ldd] (\bun,-5)   &
\\
&(\bun,-7)\ar[uu]\ar[r]    &\ar[ldd] \ar[u](\bdeux,-8) &&
\\
&&(\bdeux, -10)\ar[u] \ar[r]&(\bun,-9)\ar[uu] &
\\
&(\bun,-11)\ar[uu]& {}\save[]+<0cm,0ex>*{\vdots}\restore& {}\save[]+<0cm,0ex>*{\vdots}\restore&
}
\]
\[
\def\objectstyle{\scriptscriptstyle}
\xymatrix@-1.0pc{
&&&&\\
&&(\bdeux,0)\ar[rd]
& 
\\
&& \ar[u](\bdeux,-2)\ar[lddd]&\ar[ld] (\bun,-1)&
\\
&{(\bun,-3)}\ar[ru] & \ar[u]\ar[d]\fbox{$\scriptscriptstyle(\bdeux,-4)$} &&
\\
&&\ar[ld](\bdeux,-6)\ar[r] &\ar[uu]\ar[ldd] (\bun,-5)   &
\\
&(\bun,-7)\ar[uu]\ar[r]\ar[ruu]    &\ar[ldd] \ar[u](\bdeux,-8) &&
\\
&&(\bdeux, -10)\ar[u] \ar[r]&(\bun,-9)\ar[uu] &
\\
&(\bun,-11)\ar[uu]& {}\save[]+<0cm,0ex>*{\vdots}\restore& {}\save[]+<0cm,0ex>*{\vdots}\restore&
}
\xymatrix@-1.0pc{
&&&&\\
&&(\bdeux,0)\ar[rd]
& 
\\
&& \ar[u](\bdeux,-2)\ar[lddd]&\ar[ld] (\bun,-1)&
\\
&{(\bun,-3)}\ar[ru] & \ar[u](\bdeux,-4)\ar[rd] &&
\\
&&\ar[u]\ar[d]\fbox{$\scriptscriptstyle(\bdeux,-6)$}&\ar[l]\ar[uu]\ar[ldd] (\bun,-5)   &
\\
&(\bun,-7)\ar[uu]\ar[ru]    &\ar[ldd] (\bdeux,-8)\ar[ru] &&
\\
&&(\bdeux, -10)\ar[u] \ar[r]&(\bun,-9)\ar[uu] &
\\
&(\bun,-11)\ar[uu]& {}\save[]+<0cm,0ex>*{\vdots}\restore& {}\save[]+<0cm,0ex>*{\vdots}\restore&
}
\xymatrix@-1.0pc{
&&&&\\
&&(\bdeux,0)\ar[rd]
& 
\\
&& \ar[u](\bdeux,-2)\ar[lddd]&\ar[ld] (\bun,-1)&
\\
&{(\bun,-3)}\ar[ru] & \ar[u](\bdeux,-4)\ar[rd] &&
\\
&&\ar[u](\bdeux,-6)\ar[lddd]&\ar[ld]\ar[uu] (\bun,-5)   &
\\
&(\bun,-7)\ar[uu]\ar[ru]    &\ar[u]\ar[d] \fbox{$\scriptscriptstyle(\bdeux,-8)$} &&
\\
&&(\bdeux, -10)\ar[ld] \ar[r]&(\bun,-9)\ar[uu] &
\\
&(\bun,-11)\ar[uu]\ar[ruu]& {}\save[]+<0cm,0ex>*{\vdots}\restore& {}\save[]+<0cm,0ex>*{\vdots}\restore&
}
\]
\[
\def\objectstyle{\scriptscriptstyle}
\xymatrix@-1.0pc{
&&&&\\
&&(\bdeux,0)\ar[rd]
& 
\\
&& \ar[u](\bdeux,-2)\ar[lddd]&\ar[ld] (\bun,-1)&
\\
&{(\bun,-3)}\ar[ru] & \ar[u](\bdeux,-4)\ar[rd] &&
\\
&&\ar[u](\bdeux,-6)\ar[lddd]&\ar[ld]\ar[uu] (\bun,-5)   &
\\
&(\bun,-7)\ar[uu]\ar[ru]    &\ar[u] (\bdeux,-8)\ar[rd] &&
\\
&&\fbox{$\scriptscriptstyle(\bdeux, -10)$}\ar[u] &\ar[l](\bun,-9)\ar[uu] &
\\
&(\bun,-11)\ar[uu]\ar[ru]& {}\save[]+<0cm,0ex>*{\vdots}\restore& {}\save[]+<0cm,0ex>*{\vdots}\restore&
}
\xymatrix@-1.0pc{
&&&&\\
&&(\bdeux,0)\ar[rd]
& 
\\
&& \ar[u](\bdeux,-2)\ar[ld]&\ar[ld] (\bun,-1)&
\\
&\fbox{$\scriptscriptstyle(\bun,-3)$}\ar[dd] & \ar[u](\bdeux,-4)\ar[rd] &&
\\
&&\ar[u](\bdeux,-6)\ar[lddd]&\ar[ld]\ar[uu] (\bun,-5)   &
\\
&(\bun,-7)\ar[ru]    &\ar[u] (\bdeux,-8)\ar[rd] &&
\\
&&(\bdeux, -10)\ar[u] &\ar[l](\bun,-9)\ar[uu] &
\\
&(\bun,-11)\ar[uu]\ar[ru]& {}\save[]+<0cm,0ex>*{\vdots}\restore& {}\save[]+<0cm,0ex>*{\vdots}\restore&
}
\xymatrix@-1.0pc{
&&&&\\
&&(\bdeux,0)\ar[rd]
& 
\\
&& \ar[u](\bdeux,-2)\ar[ld]&\ar[ld] (\bun,-1)&
\\
&(\bun,-3)\ar[rd] & \ar[u](\bdeux,-4)\ar[rd] &&
\\
&&\ar[u](\bdeux,-6)\ar[ld]&\ar[ld]\ar[uu] (\bun,-5)   &
\\
&\fbox{$\scriptscriptstyle(\bun,-7)$}\ar[uu]\ar[dd]    &\ar[u] (\bdeux,-8)\ar[rd] &&
\\
&&(\bdeux, -10)\ar[u] &\ar[l](\bun,-9)\ar[uu] &
\\
&(\bun,-11)\ar[ru]& {}\save[]+<0cm,0ex>*{\vdots}\restore& {}\save[]+<0cm,0ex>*{\vdots}\restore&
}
\]

\subsection{Mutation sequence in type $G_2$}\label{appendixG2}

We display the sequence of mutated quivers obtained from $G^-$ 
at each step of the mutation sequence $\mu_\SS$. 

\[
\def\objectstyle{\scriptscriptstyle}
\xymatrix@-1.0pc{
&(\bdeux,0)\ar[rd]
&
\\
&\ar[ld]\ar[u] (\bdeux,-2)& \ar[ldd](\bun,-1)&
\\
{(\bun,-3)}\ar[ddr] &\ar[u] (\bdeux,-4) \ar[rrd]&&
\\
&\ar[u] \ar[rd](\bdeux,-6) && (\bun,-5)\ar[lldd]   
\\
  &(\bdeux,-8)\ar[u]\ar[dl]& \ar[uuu](\bun,-7)\ar[ldd] &
\\
(\bun,-9)\ar[uuu]&(\bdeux,-10)\ar[u]\ar[rrd]&&
\\ 
{}\save[]+<0cm,0ex>*{\vdots}\restore \ar[u]& (\bdeux,-12)\ar[u] &{}\save[]+<0cm,0ex>*{\vdots}\restore\ar[uu]&(\bun,-11)\ar[uuu]
}
\xymatrix@-1.0pc{
&\fbox{$\scriptscriptstyle(\bdeux,0)$}\ar[d]
&
\\
&\ar[ld]\ar[r] (\bdeux,-2)& \ar[ldd](\bun,-1)\ar[lu]&
\\
{(\bun,-3)}\ar[ddr] &\ar[u] (\bdeux,-4) \ar[rrd]&&
\\
&\ar[u] \ar[rd](\bdeux,-6) && (\bun,-5)\ar[lldd]   
\\
  &(\bdeux,-8)\ar[u]\ar[dl]& \ar[uuu](\bun,-7)\ar[ldd] &
\\
(\bun,-9)\ar[uuu]&(\bdeux,-10)\ar[u]\ar[rrd]&&
\\ 
{}\save[]+<0cm,0ex>*{\vdots}\restore \ar[u]& (\bdeux,-12)\ar[u] &{}\save[]+<0cm,0ex>*{\vdots}\restore\ar[uu]&(\bun,-11)\ar[uuu]
}
\xymatrix@-1.0pc{
&(\bdeux,0)\ar[ddl]
&
\\
&\ar[d]\ar[u] \fbox{$\scriptscriptstyle(\bdeux,-2)$}& \ar[ldd]\ar[l](\bun,-1)&
\\
{(\bun,-3)}\ar[ru]\ar[ddr] & (\bdeux,-4) \ar[rrd]\ar[ur]\ar[l]&&
\\
&\ar[u] \ar[rd](\bdeux,-6) && (\bun,-5)\ar[lldd]   
\\
  &(\bdeux,-8)\ar[u]\ar[dl]& \ar[uuu](\bun,-7)\ar[ldd] &
\\
(\bun,-9)\ar[uuu]&(\bdeux,-10)\ar[u]\ar[rrd]&&
\\ 
{}\save[]+<0cm,0ex>*{\vdots}\restore \ar[u]& (\bdeux,-12)\ar[u] &{}\save[]+<0cm,0ex>*{\vdots}\restore\ar[uu]&(\bun,-11)\ar[uuu]
}
\]
\[
\def\objectstyle{\scriptscriptstyle}
\xymatrix@-1.0pc{
&(\bdeux,0)\ar[ddl]
&
\\
&\ar[u](\bdeux,-2)\ar[ddrr]& (\bun,-1)\ar[ld]&
\\
{(\bun,-3)}\ar[r]\ar[ddr] & \fbox{$\scriptscriptstyle(\bdeux,-4)$}\ar[u]\ar[d] &&
\\
& \ar[rd]\ar[ul]\ar[rr](\bdeux,-6) && (\bun,-5)  \ar[llu]\ar[lldd]
\\
  &(\bdeux,-8)\ar[u]\ar[dl]& \ar[uuu](\bun,-7)\ar[ldd] &
\\
(\bun,-9)\ar[uuu]&(\bdeux,-10)\ar[u]\ar[rrd]&&
\\ 
{}\save[]+<0cm,0ex>*{\vdots}\restore \ar[u]& (\bdeux,-12)\ar[u] &{}\save[]+<0cm,0ex>*{\vdots}\restore\ar[uu]&(\bun,-11)\ar[uuu]
}
\xymatrix@-1.0pc{
&(\bdeux,0)\ar[ddl]
&
\\
&\ar[u](\bdeux,-2)\ar[ddrr]& (\bun,-1)\ar[ld]&
\\
{(\bun,-3)}\ar[dr] & (\bdeux,-4)\ar[u]\ar[ddr] &&
\\
& 
\fbox{$\scriptscriptstyle(\bdeux,-6)$}\ar[u]\ar[d]
 && (\bun,-5) \ar[ll]\ar[ddll] 
\\
  &(\bdeux,-8)\ar[dl]\ar[r]\ar[urr]& \ar[uuu](\bun,-7)\ar[ul]\ar[ldd] &
\\
(\bun,-9)\ar[uuu]&(\bdeux,-10)\ar[u]\ar[rrd]&&
\\ 
{}\save[]+<0cm,0ex>*{\vdots}\restore \ar[u]& (\bdeux,-12)\ar[u] &{}\save[]+<0cm,0ex>*{\vdots}\restore\ar[uu]&(\bun,-11)\ar[uuu]
}
\xymatrix@-1.0pc{
&(\bdeux,0)\ar[ddl]
&
\\
&\ar[u](\bdeux,-2)\ar[ddrr]& (\bun,-1)\ar[ld]&
\\
{(\bun,-3)}\ar[dr] & (\bdeux,-4)\ar[u]\ar[ddr] &&
\\
& 
(\bdeux,-6)\ar[u]\ar[ddl]
 && (\bun,-5) \ar[dll]
\\
  & \fbox{$\scriptscriptstyle(\bdeux,-8)$}\ar[u]\ar[d] & \ar[uuu](\bun,-7)\ar[l]\ar[ldd] &
\\
(\bun,-9)\ar[uuu]\ar[ur]&(\bdeux,-10)\ar[rrd]\ar[l]\ar[ur]&&
\\ 
{}\save[]+<0cm,0ex>*{\vdots}\restore \ar[u]& (\bdeux,-12)\ar[u] &{}\save[]+<0cm,0ex>*{\vdots}\restore\ar[uu]&(\bun,-11)\ar[uuu]
}
\]
\[
\def\objectstyle{\scriptscriptstyle}
\xymatrix@-1.0pc{
&(\bdeux,0)\ar[ddl]
&
\\
&\ar[u](\bdeux,-2)\ar[ddrr]& (\bun,-1)\ar[ld]&
\\
{(\bun,-3)}\ar[dr] & (\bdeux,-4)\ar[u]\ar[ddr] &&
\\
& 
(\bdeux,-6)\ar[u]\ar[ddl]
 && (\bun,-5) \ar[dll]
\\
  & (\bdeux,-8)\ar[u]\ar[ddrr] & \ar[uuu](\bun,-7)\ar[ld] &
\\
(\bun,-9)\ar[uuu]\ar[r]&\fbox{$\scriptscriptstyle(\bdeux,-10)$}\ar[u]\ar[d]&&
\\ 
{}\save[]+<0cm,0ex>*{\vdots}\restore \ar[u]& (\bdeux,-12)\ar[ul]\ar[rr] &{}\save[]+<0cm,0ex>*{\vdots}\restore\ar[uu]&(\bun,-11)\ar[uuu]\ar[llu]
}
\xymatrix@-1.0pc{
&(\bdeux,0)\ar[ddl]
&
\\
&\ar[u](\bdeux,-2)\ar[ddrr]& (\bun,-1)\ar[ld]&
\\
{(\bun,-3)}\ar[dr] & (\bdeux,-4)\ar[u]\ar[ddr] &&
\\
& 
(\bdeux,-6)\ar[u]\ar[ddl]
 && (\bun,-5) \ar[dll]
\\
  & (\bdeux,-8)\ar[u]\ar[ddrr] & \ar[uuu](\bun,-7)\ar[ld] &
\\
(\bun,-9)\ar[uuu]\ar[rd]&(\bdeux,-10)\ar[u]&&
\\ 
{}\save[]+<0cm,0ex>*{\vdots}\restore \ar[u]& \fbox{$\scriptscriptstyle(\bdeux,-12)$}\ar[u] &{}\save[]+<0cm,0ex>*{\vdots}\restore\ar[uu]&(\bun,-11)\ar[uuu]\ar[ll]
}
\xymatrix@-1.0pc{
&(\bdeux,0)\ar[ddl]
&
\\
&\ar[u](\bdeux,-2)\ar[ddrr]& 
\fbox{$\scriptscriptstyle(\bun,-1)$}\ar[ddd]&
\\
{(\bun,-3)}\ar[dr] & (\bdeux,-4)\ar[ru]\ar[u] &&
\\
& 
(\bdeux,-6)\ar[u]\ar[ldd]
 && (\bun,-5) \ar[lld] 
\\
 {} & (\bdeux,-8)\ar[u]\ar[rrdd] & (\bun,-7)  \ar[ld]&{}
\\
(\bun,-9)\ar[uuu]\ar[rd]&(\bdeux,-10)\ar[u]&&
\\ 
{}\save[]+<0cm,0ex>*{\vdots}\restore \ar[u]& (\bdeux,-12)\ar[u] &{}\save[]+<0cm,0ex>*{\vdots}\restore\ar[uu]&(\bun,-11)\ar[uuu]
}
\]

\[
\def\objectstyle{\scriptscriptstyle}
\xymatrix@-1.0pc{
&(\bdeux,0)\ar[ddl]
&
\\
&\ar[u](\bdeux,-2)\ar[ddrr]& 
(\bun,-1)\ar[ldddd]&
\\
{(\bun,-3)}\ar[dr] & (\bdeux,-4)\ar[ru]\ar[u] &&
\\
& 
(\bdeux,-6)\ar[u]\ar[ldd]
 && (\bun,-5) \ar[lld] 
\\
 {} & (\bdeux,-8)\ar[u]\ar[rrdd] & \fbox{$\scriptscriptstyle(\bun,-7)$}\ar[uuu]\ar[dd] &{}
\\
(\bun,-9)\ar[uuu]\ar[rd]&(\bdeux,-10)\ar[u]\ar[ur]&&
\\ 
{}\save[]+<0cm,0ex>*{\vdots}\restore \ar[u]& (\bdeux,-12)\ar[u] &{}\save[]+<0cm,0ex>*{\vdots}\restore&(\bun,-11)\ar[uuu]
}
\xymatrix@-1.0pc{
&\fbox{$\scriptscriptstyle(\bdeux,0)$}\ar[d]
&
\\
&(\bdeux,-2)\ar[ddrr]\ar[dl]& 
(\bun,-1)\ar[ldddd]&
\\
{(\bun,-3)}\ar[uur]\ar[dr] & (\bdeux,-4)\ar[ru]\ar[u] &&
\\
& 
(\bdeux,-6)\ar[u]\ar[ldd]
 && (\bun,-5) \ar[lld] 
\\
 {} & (\bdeux,-8)\ar[u]\ar[rrdd] & (\bun,-7)\ar[uuu] &{}
\\
(\bun,-9)\ar[uuu]\ar[rd]&(\bdeux,-10)\ar[u]\ar[ur]&&
\\ 
{}\save[]+<0cm,0ex>*{\vdots}\restore \ar[u]& (\bdeux,-12)\ar[u] &{}\save[]+<0cm,0ex>*{\vdots}\restore\ar[uu]&(\bun,-11)\ar[uuu]
}
\xymatrix@-1.0pc{
&(\bdeux,0)\ar[dddrr]
&
\\
&\fbox{$\scriptscriptstyle(\bdeux,-2)$}\ar[u]\ar[d]
& 
(\bun,-1)\ar[ldddd]&
\\
{(\bun,-3)}\ar[ur]\ar[dr]& (\bdeux,-4)\ar[ru]\ar[l]\ar[drr] &&
\\
{}& 
(\bdeux,-6)\ar[u]\ar[ldd]
 && (\bun,-5)\ar[uull] \ar[dll] 
\\
 {} & (\bdeux,-8)\ar[u]\ar[rrdd] & (\bun,-7)\ar[uuu] &{}
\\
(\bun,-9)\ar[uuu]\ar[rd]&(\bdeux,-10)\ar[u]\ar[ur]&&
\\ 
{}\save[]+<0cm,0ex>*{\vdots}\restore \ar[u]& (\bdeux,-12)\ar[u] &{}\save[]+<0cm,0ex>*{\vdots}\restore\ar[uu]&(\bun,-11)\ar[uuu]
} 
\]

\[
\def\objectstyle{\scriptscriptstyle}
\xymatrix@-1.0pc{
&(\bdeux,0)\ar[dddrr]
&
\\
&(\bdeux,-2)\ar[u]\ar[r]
& 
(\bun,-1)\ar[dl]\ar[ldddd]&
\\
{(\bun,-3)}\ar[r]& \fbox{$\scriptscriptstyle(\bdeux,-4)$}\ar[d]\ar[u] &&
\\
& 
(\bdeux,-6)\ar[rr]\ar[uur]\ar[ddl]
 && (\bun,-5)\ar[ull]\ar[dll] 
\\
 {} & (\bdeux,-8)\ar[u]\ar[rrdd] & (\bun,-7)\ar[uuu] &{}
\\
(\bun,-9)\ar[uuu]\ar[rd]&(\bdeux,-10)\ar[u]\ar[ur]&&
\\ 
{}\save[]+<0cm,0ex>*{\vdots}\restore \ar[u]& (\bdeux,-12)\ar[u] &{}\save[]+<0cm,0ex>*{\vdots}\restore\ar[uu]&(\bun,-11)\ar[uuu]
}
\xymatrix@-1.0pc{
&(\bdeux,0)\ar[dddrr]
&
\\
&(\bdeux,-2)\ar[u]\ar[r]
& 
(\bun,-1)\ar[ddl]\ar[ldddd]&
\\
{(\bun,-3)}\ar[r]& (\bdeux,-4)\ar[u]\ar[lddd] \ar[dddl]&&
\\
& 
\fbox{$\scriptscriptstyle(\bdeux,-6)$}\ar[d]\ar[u]
 && (\bun,-5)\ar[ll] 
\\
 {} & (\bdeux,-8)\ar[dl]\ar[ruuu]\ar[rrdd] & (\bun,-7)\ar[uuu] &{}
\\
(\bun,-9)\ar[uuu]\ar[rd]\ar[uur]&(\bdeux,-10)\ar[ur]\ar[u]&&
\\ 
{}\save[]+<0cm,0ex>*{\vdots}\restore \ar[u]& (\bdeux,-12)\ar[u] &{}\save[]+<0cm,0ex>*{\vdots}\restore\ar[uu]&(\bun,-11)\ar[uuu]
}
\xymatrix@-1.0pc{
&(\bdeux,0)\ar[dddrr]
&
\\
&(\bdeux,-2)\ar[u]\ar[r]
& 
(\bun,-1)\ar[dddl]&
\\
{(\bun,-3)}\ar[r]& (\bdeux,-4)\ar[u] \ar[dddl]&&
\\
& 
(\bdeux,-6)\ar[u]\ar[rrddd]
 && (\bun,-5)\ar[ll] 
\\
  &\fbox{$\scriptscriptstyle(\bdeux,-8)$}\ar[u]\ar[d]& 
   (\bun,-7)\ar[uuu]
     &
  \\
(\bun,-9)\ar[uuu]\ar[rd]\ar[ur]&(\bdeux,-10)\ar[l]\ar[ur]\ar[drr]&&
\\ 
{}\save[]+<0cm,0ex>*{\vdots}\restore \ar[u]& (\bdeux,-12)\ar[u] &{}\save[]+<0cm,0ex>*{\vdots}\restore\ar[uu]&(\bun,-11)\ar[uull]\ar[uuu] }
\]

\[
\def\objectstyle{\scriptscriptstyle}
\xymatrix@-1.0pc{
&(\bdeux,0)\ar[dddrr]
&
\\
&(\bdeux,-2)\ar[u]\ar[r]
& 
(\bun,-1)\ar[dddl]&
\\
{(\bun,-3)}\ar[r]& (\bdeux,-4)\ar[u] \ar[dddl]&&
\\
& 
(\bdeux,-6)\ar[u]\ar[rrddd]
 && (\bun,-5)\ar[ll] 
\\
  &(\bdeux,-8)\ar[u]\ar[r]& 
   (\bun,-7)\ar[uuu]\ar[dl]
     &
  \\
(\bun,-9)\ar[uuu]\ar[r]&\fbox{$\scriptscriptstyle(\bdeux,-10)$}\ar[u]\ar[d]&&
\\ 
{}\save[]+<0cm,0ex>*{\vdots}\restore \ar[u]& (\bdeux,-12)\ar[rr]\ar[uur] &{}\save[]+<0cm,0ex>*{\vdots}\restore\ar[uu]&(\bun,-11)\ar[uuu]\ar[llu] }
\xymatrix@-1.0pc{
&(\bdeux,0)\ar[dddrr]
&
\\
&(\bdeux,-2)\ar[u]\ar[r]
& 
(\bun,-1)\ar[dddl]&
\\
{(\bun,-3)}\ar[r]& (\bdeux,-4)\ar[u] \ar[dddl]&&
\\
& 
(\bdeux,-6)\ar[u]\ar[rrddd]
 && (\bun,-5)\ar[ll] 
\\
  &(\bdeux,-8)\ar[u]\ar[r]& 
   (\bun,-7)\ar[uuu]\ar[ddl]
     &
  \\
(\bun,-9)\ar[uuu]\ar[r]\ar[rd]&(\bdeux,-10)\ar[u]&&
\\ 
{}\save[]+<0cm,0ex>*{\vdots}\restore \ar[u]& \fbox{$\scriptscriptstyle(\bdeux,-12)$}\ar[u] &{}\save[]+<0cm,0ex>*{\vdots}\restore\ar[uu]&(\bun,-11)\ar[uuu]\ar[ll] }
\xymatrix@-1.0pc{
&(\bdeux,0)\ar[dddrr]
&
\\
&(\bdeux,-2)\ar[u]\ar[r]
& 
(\bun,-1)\ar[dddl]&
\\
\fbox{$\scriptscriptstyle(\bun,-3)$}\ar[ddd]& (\bdeux,-4)\ar[u]\ar[l] &&
\\
& 
(\bdeux,-6)\ar[u]\ar[rrddd]
 && (\bun,-5)\ar[ll] 
 \\
  &(\bdeux,-8)\ar[u]\ar[r]& 
   (\bun,-7)\ar[uuu]
     &
  \\
(\bun,-9)\ar[r]&(\bdeux,-10)\ar[u]&&
\\ 
{}\save[]+<0cm,0ex>*{\vdots}\restore \ar[u]& (\bdeux,-12)\ar[u] &{}\save[]+<0cm,0ex>*{\vdots}\restore\ar[uu]&(\bun,-11)\ar[uuu]\ar[ll] }
\]

\[
\def\objectstyle{\scriptscriptstyle}
\xymatrix@-1.0pc{
&(\bdeux,0)\ar[dddrr]
&
\\
&(\bdeux,-2)\ar[u]\ar[r]
& 
(\bun,-1)\ar[dddl]&
\\
(\bun,-3)\ar[dddr]& (\bdeux,-4)\ar[u]\ar[l] &&
\\
& 
(\bdeux,-6)\ar[u]\ar[rrddd]
 && (\bun,-5)\ar[ll] 
 \\
  &(\bdeux,-8)\ar[u]\ar[r]& 
   (\bun,-7)\ar[uuu]
     &
  \\
\fbox{$\scriptscriptstyle(\bun,-9)$}\ar[d]\ar[uuu]&(\bdeux,-10)\ar[l]\ar[u]&&
\\ 
{}\save[]+<0cm,0ex>*{\vdots}\restore & (\bdeux,-12)\ar[u] &{}\save[]+<0cm,0ex>*{\vdots}\restore\ar[uu]&(\bun,-11)\ar[uuu]\ar[ll] }
\xymatrix@-1.0pc{
&\fbox{$\scriptscriptstyle(\bdeux,0)$}\ar[d]
&
\\
&(\bdeux,-2)\ar[r]\ar[ddrr]
& 
(\bun,-1)\ar[dddl]&
\\
(\bun,-3)\ar[dddr]& (\bdeux,-4)\ar[u]\ar[l] &&
\\
& 
(\bdeux,-6)\ar[u]\ar[rrddd]
 && (\bun,-5)\ar[ll]\ar[uuull] 
 \\
  &(\bdeux,-8)\ar[u]\ar[r]& 
   (\bun,-7)\ar[uuu]&
  \\
(\bun,-9)\ar[uuu]&(\bdeux,-10)\ar[l]\ar[u]&&
\\ 
{}\save[]+<0cm,0ex>*{\vdots}\restore \ar[u]& (\bdeux,-12)\ar[u] &{}\save[]+<0cm,0ex>*{\vdots}\restore\ar[uu]&(\bun,-11)\ar[uuu]\ar[ll] }
\xymatrix@-1.0pc{
&(\bdeux,0)\ar[dr]
&
\\
&\fbox{$\scriptscriptstyle(\bdeux,-2)$}\ar[u]\ar[d]
& 
(\bun,-1)\ar[dddl]\ar[l]&
\\
(\bun,-3)\ar[rddd]& (\bdeux,-4)\ar[ur]\ar[l]\ar[rrd] &&
\\
& 
(\bdeux,-6)\ar[u]\ar[rrddd]
 && (\bun,-5)\ar[ll]\ar[uull] 
 \\
  &(\bdeux,-8)\ar[u]\ar[r]& 
   (\bun,-7)\ar[uuu]&
  \\
(\bun,-9)\ar[uuu]&(\bdeux,-10)\ar[l]\ar[u]&&
\\ 
{}\save[]+<0cm,0ex>*{\vdots}\restore \ar[u]& (\bdeux,-12)\ar[u] &{}\save[]+<0cm,0ex>*{\vdots}\restore\ar[uu]&(\bun,-11)\ar[uuu]\ar[ll] }
\]

\[
\def\objectstyle{\scriptscriptstyle}
\xymatrix@-1.0pc{
&(\bdeux,0)\ar[dr]
&
\\
&(\bdeux,-2)\ar[u]\ar[ld]
& 
(\bun,-1)\ar[dl]\ar[dddl]&
\\
(\bun,-3)\ar[r]\ar[rddd]& \fbox{$\scriptscriptstyle(\bdeux,-4)$}\ar[d]\ar[u] &&
\\
& 
(\bdeux,-6)\ar[ul]\ar[uur]\ar[rrddd]
 && (\bun,-5)\ar[ull] 
 \\
  &(\bdeux,-8)\ar[u]\ar[r]& 
   (\bun,-7)\ar[uuu]&
  \\
(\bun,-9)\ar[uuu]&(\bdeux,-10)\ar[l]\ar[u]&&
\\ 
{}\save[]+<0cm,0ex>*{\vdots}\restore \ar[u]& (\bdeux,-12)\ar[u] &{}\save[]+<0cm,0ex>*{\vdots}\restore\ar[uu]&(\bun,-11)\ar[uuu]\ar[ll] }
\xymatrix@-1.0pc{
&(\bdeux,0)\ar[dr]
&
\\
&(\bdeux,-2)\ar[u]\ar[ld]
& 
(\bun,-1)\ar[ddl]&
\\
(\bun,-3)\ar[rd]\ar[rddd]& (\bdeux,-4)\ar[u]\ar[rrdddd] &&
\\
& 
\fbox{$\scriptscriptstyle(\bdeux,-6)$}\ar[u]\ar[d]
 && (\bun,-5)\ar[ull] 
  \\
  &(\bdeux,-8)\ar[r]\ar[uul]\ar[rrdd]& 
   (\bun,-7)\ar[uuu]&
  \\
(\bun,-9)\ar[uuu]&(\bdeux,-10)\ar[l]\ar[u]&&
\\ 
{}\save[]+<0cm,0ex>*{\vdots}\restore \ar[u]& (\bdeux,-12)\ar[u] &{}\save[]+<0cm,0ex>*{\vdots}\restore\ar[uu]&(\bun,-11)\ar[uuu]\ar[ll]\ar[uuull] }
\xymatrix@-1.0pc{
&(\bdeux,0)\ar[dr]
&
\\
&(\bdeux,-2)\ar[u]\ar[ld]
& 
(\bun,-1)\ar[ddl]&
\\
(\bun,-3)\ar[rdd]& (\bdeux,-4)\ar[u]\ar[rrdddd] &&
\\
& 
(\bdeux,-6)\ar[u]\ar[dr]
 && (\bun,-5)\ar[ull] 
  \\
  &\fbox{$\scriptscriptstyle(\bdeux,-8)$}\ar[u]\ar[d]& 
   (\bun,-7)\ar[uuu]\ar[l]&
  \\
(\bun,-9)\ar[uuu]&(\bdeux,-10)\ar[l]\ar[ur]\ar[drr]&&
\\ 
{}\save[]+<0cm,0ex>*{\vdots}\restore \ar[u]& (\bdeux,-12)\ar[u] &{}\save[]+<0cm,0ex>*{\vdots}\restore\ar[uu]&(\bun,-11)\ar[uull]\ar[uuu]\ar[ll] }
\]

\[
\def\objectstyle{\scriptscriptstyle}
\xymatrix@-1.0pc{
&(\bdeux,0)\ar[dr]
&
\\
&(\bdeux,-2)\ar[u]\ar[ld]
& 
(\bun,-1)\ar[ddl]&
\\
(\bun,-3)\ar[rdd]& (\bdeux,-4)\ar[u]\ar[rrdddd] &&
\\
& 
(\bdeux,-6)\ar[u]\ar[dr]
 && (\bun,-5)\ar[ull] 
  \\
  &(\bdeux,-8)\ar[u]\ar[dl]& 
   (\bun,-7)\ar[uuu]\ar[dl]&
  \\
(\bun,-9)\ar[uuu]\ar[r]&\fbox{$\scriptscriptstyle(\bdeux,-10)$}\ar[u]\ar[d]&&
\\ 
{}\save[]+<0cm,0ex>*{\vdots}\restore \ar[u]& (\bdeux,-12)\ar[ul]\ar[uur] &{}\save[]+<0cm,0ex>*{\vdots}\restore\ar[uu]&(\bun,-11)\ar[ull]\ar[uuu]}
\xymatrix@-1.0pc{
&(\bdeux,0)\ar[dr]
&
\\
&(\bdeux,-2)\ar[u]\ar[ld]
& 
(\bun,-1)\ar[ddl]&
\\
(\bun,-3)\ar[rdd]& (\bdeux,-4)\ar[u]\ar[rrdddd] &&
\\
& 
(\bdeux,-6)\ar[u]\ar[dr]
 && (\bun,-5)\ar[ull] 
  \\
  &(\bdeux,-8)\ar[u]\ar[dl]& 
   (\bun,-7)\ar[uuu]\ar[ddl]&
  \\
(\bun,-9)\ar[uuu]\ar[dr]&(\bdeux,-10)\ar[u]&&
\\ 
{}\save[]+<0cm,0ex>*{\vdots}\restore \ar[u]&\fbox{$\scriptscriptstyle(\bdeux,-12)$}\ar[u] &{}\save[]+<0cm,0ex>*{\vdots}\restore\ar[uu]&(\bun,-11)\ar[ull]\ar[uuu]}
\xymatrix@-1.0pc{
&(\bdeux,0)\ar[dr]
&
\\
&(\bdeux,-2)\ar[u]\ar[ld]
& 
(\bun,-1)\ar[ddl]&
\\
(\bun,-3)\ar[rdd]& (\bdeux,-4)\ar[u]\ar[rrd] &&
\\
& 
(\bdeux,-6)\ar[u]\ar[dr]
 && \fbox{$\scriptscriptstyle(\bun,-5)$}\ar[ddd]
  \\
  &(\bdeux,-8)\ar[u]\ar[dl]& 
   (\bun,-7)\ar[uuu]\ar[ddl]&
  \\
(\bun,-9)\ar[uuu]&(\bdeux,-10)\ar[u]&&
\\ 
{}\save[]+<0cm,0ex>*{\vdots}\restore \ar[u]&(\bdeux,-12)\ar[u] &{}\save[]+<0cm,0ex>*{\vdots}\restore\ar[uu]&(\bun,-11)\ar[ull]}
\]

\subsection{Examples of $A$-modules for $\g$ of type $B_2$}\label{append-B2}

We describe some $A$-modules $K^{(i)}_{k,m}$ for $\g$ of type $B_2$.
The quiver $\G^-$ is: 
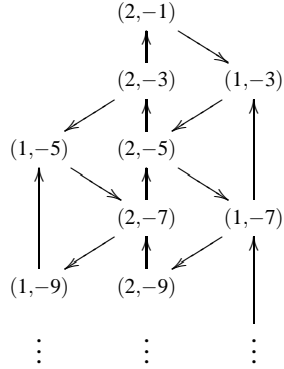
\[
\def\objectstyle{\scriptstyle}
\xymatrix@-1.0pc{
&  & (\bdeux,-1) \ar[rd]&&
\\
&&\ar[u]\ar[ld] (\bdeux,-3) &\ar[ld] (\bun,-3) &
\\
&\ar[rd](\bun,-5) & \ar[u]\ar[rd](\bdeux,-5)&&
\\
&&\ar[u]\ar[ld] (\bdeux,-7) &\ar[ld] (\bun,-7) \ar[uu]  &
\\
&\ar[uu](\bun,-9)   & \ar[u](\bdeux,-9) &&
\\
&{}\save[]+<0cm,0ex>*{\vdots}\restore  &{}\save[]+<0cm,0ex>*{\vdots}\restore&{}\save[]+<0cm,0ex>*{\vdots}\ar[uu]\restore 
\\
}
\]
Following the convention of Example~\ref{example_A3_modules}, unless otherwise
specified, in the following figures the vertices carry one-dimensional spaces, and the arrows carry
linear maps with matrix $\pmatrix{\pm1}$.

The modules $K^{(\bun)}_{1,-5}$ and $K^{(\bun)}_{1,-7}$ are:
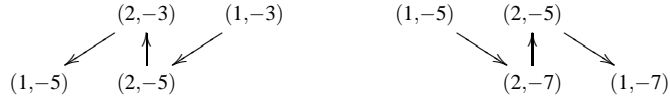
\[
\def\objectstyle{\scriptstyle}
\xymatrix@-1.0pc{
&&\ar[ld] (\bdeux,-3)& \ar[ld](\bun,-3)&
\\
&{(\bun,-5)} &\ar[u] (\bdeux,-5)&&
}
\xymatrix@-1.0pc{
&{(\bun,-5)}\ar[rd] & (\bdeux,-5) \ar[rd]&&
\\
&&\ar[u] (\bdeux,-7) &(\bun,-7)  &
}
\]
The modules $K^{(\bdeux)}_{1,-5}$ and $K^{(\bdeux)}_{1,-7}$ are:
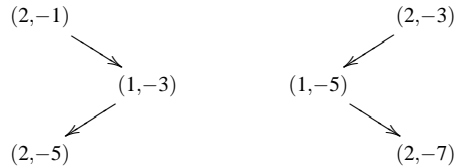
\[
\def\objectstyle{\scriptstyle}
\xymatrix@-1.0pc{
(\bdeux,-1) \ar[rd]&&
\\
 &\ar[ld] (\bun,-3)   &
\\
(\bdeux,-5) &&
}
\xymatrix@-1.0pc{
&&\ar[ld] (\bdeux,-3)& 
\\
&{(\bun,-5)}\ar[rd] &
\\
&&(\bdeux,-7) 
}
\]
Applying Theorem~\ref{th_geom_form}, we recover the following well known formulas for the $q$-characters 
of the fundamental $U_q(\hg)$-modules:  
\[
\begin{array}{rcl}
\chi_q(L(Y_{\bun,-7}))&=& Y_{\bun,-7}(1+v_{\bun,-5}(1 + v_{\bdeux,-3}
(1+ v_{\bdeux,-5}(1 + v_{\bun,-3})))),\\[3mm]
\chi_q(L(Y_{\bdeux,-6}))&=& Y_{\bdeux,-6}(1+v_{\bdeux,-5}(1 + v_{\bun,-3}
(1+ v_{\bdeux,-1}))).  
\end{array}
\]

The modules $K^{(\bun)}_{2,-5}$ and $K^{(\bdeux)}_{2,-7}$ are:
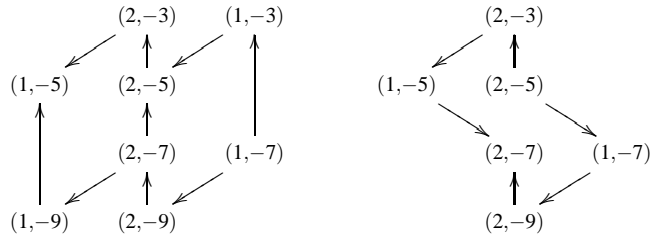
\[
\def\objectstyle{\scriptstyle}
\xymatrix@-1.0pc{
&\ar[ld] (\bdeux,-3)& \ar[ld](\bun,-3)&
\\
{(\bun,-5)}&\ar[u] (\bdeux,-5) &&
\\
&\ar[u]\ar[ld] (\bdeux,-7) &\ar[ld] (\bun,-7) \ar[uu]  &
\\
\ar[uu](\bun,-9)   & \ar[u](\bdeux,-9)&&
}
\quad
\xymatrix@-1.0pc{
&\ar[ld] (\bdeux,-3)&&
\\
{(\bun,-5)}\ar[rd] &\ar[u] (\bdeux,-5) \ar[rd]&&
\\
& (\bdeux,-7) &\ar[ld] (\bun,-7) &
\\
& \ar[u](\bdeux,-9)&&
}
\]
They correspond under Theorem~\ref{th_geom_form} to the Kirillov-Reshetikhin modules 
\[
W^{(\bun)}_{2,-11}= L(Y_{\bun,-11}Y_{\bun,-7})
\quad \mbox{and}\quad 
W^{(\bdeux)}_{2,-10}= L(Y_{\bdeux,-10}Y_{\bdeux,-8}). 
\]

The modules $K^{(\bun)}_{3,-5}$ and $K^{(\bdeux)}_{3,-7}$ are:
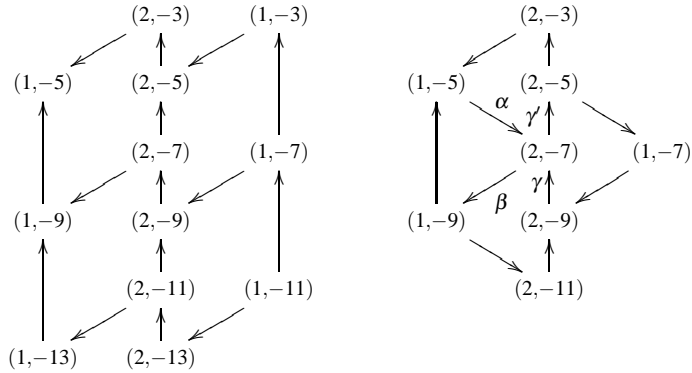
\[
\def\objectstyle{\scriptstyle}
\xymatrix@-1.0pc{
&\ar[ld] (\bdeux,-3)& \ar[ld](\bun,-3)&
\\
{(\bun,-5)}&\ar[u] (\bdeux,-5) &&
\\
&\ar[u]\ar[ld] (\bdeux,-7) &\ar[ld] (\bun,-7) \ar[uu]  &
\\
\ar[uu](\bun,-9)   & \ar[u](\bdeux,-9)&&
\\
&\ar[u]\ar[ld] (\bdeux,-11) &\ar[ld] (\bun,-11) \ar[uu]  &
\\
\ar[uu](\bun,-13)   & \ar[u](\bdeux,-13)&&
}
\quad
\xymatrix@-1.0pc{
&\ar[ld] (\bdeux,-3)&&
\\
{(\bun,-5)}\ar[rd]^{\a} &\ar[u] (\bdeux,-5) \ar[rd]&&
\\
& (\bdeux,-7)\ar[u]^{\ga'}\ar[ld]^{\b} &\ar[ld] (\bun,-7) &
\\
(\bun,-9)\ar[uu]\ar[rd]& \ar[u]^{\ga}(\bdeux,-9)&&
\\
& \ar[u](\bdeux,-11)&&
}
\]
In $K^{(\bdeux)}_{3,-7}$, the vertex $(\bdeux,-7)$ carries a two-dimensional
vector space. The linear maps carried by the adjacent arrows have the following
matrices:
\[
\a = \ga = \pmatrix{1\cr 0},\quad \b = \ga' = \pmatrix{0&1}. 
\]
They correspond under Theorem~\ref{th_geom_form} to the Kirillov-Reshetikhin modules: 
\[
W^{(\bun)}_{3,-15}= L(Y_{\bun,-15}Y_{\bun,-11}Y_{\bun,-7})
\quad \mbox{and}\quad 
W^{(\bdeux)}_{3,-12}= L(Y_{\bdeux,-12}Y_{\bdeux,-10}Y_{\bdeux,-8}). 
\]

\subsection{Examples of $A$-modules for $\g$ of type $B_3$}\label{append-B3}

Let $\g$ be of type $B_3$, with the short
root being $\a_{\btrois}$.
The quiver $\G^-$ is: 
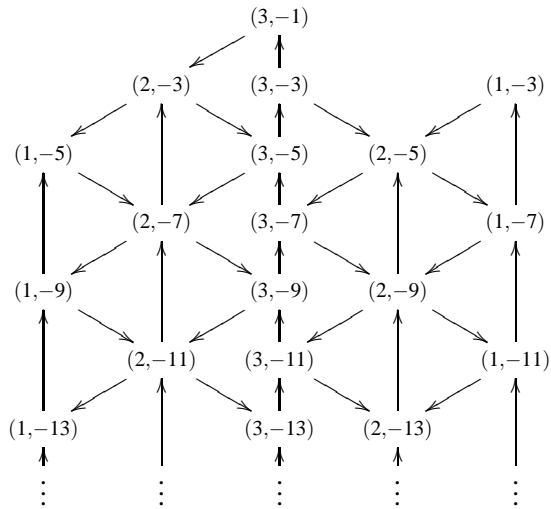
\[
\def\objectstyle{\scriptstyle}
\xymatrix@-1.0pc{
&&&\ar[ld] (\btrois,-1)& &&
\\
&&{(\bdeux,-3)}\ar[rd]\ar[ld] &\ar[u] (\btrois,-3) \ar[rd]&&(\bun,-3)\ar[ld]&
\\
&(\bun,-5)\ar[rd]&&\ar[u]\ar[ld] (\btrois,-5) &\ar[ld] (\bdeux,-5) \ar[rd]  &&
\\
&&\ar[uu](\bdeux,-7)\ar[rd]\ar[ld]   & \ar[u](\btrois,-7) \ar[rd]&&(\bun,-7)\ar[ld]\ar[uu]&
\\
&(\bun,-9)\ar[rd]\ar[uu]&&\ar[u]\ar[ld] (\btrois,-9) &\ar[ld] (\bdeux,-9) \ar[uu]\ar[rd] &&
\\
&&\ar[uu](\bdeux,-11)\ar[ld]\ar[rd] & \ar[u](\btrois,-11)\ar[rd]&&\ar[ld](\bun,-11)\ar[uu]&
\\
&(\bun,-13)\ar[uu] && (\btrois,-13)\ar[u] & (\bdeux,-13) \ar[uu]
\\
&{}\save[]+<0cm,0ex>*{\vdots}\ar[u]\restore&{}\save[]+<0cm,0ex>*{\vdots}\ar[uu]\restore  
&{}\save[]+<0cm,0ex>*{\vdots}\ar[u]\restore&{}\save[]+<0cm,0ex>*{\vdots}\ar[u]\restore 
&{}\save[]+<0cm,0ex>*{\vdots}\ar[uu]\restore 
\\
}
\]
The module $K^{(\bun)}_{1,-9}$ is:
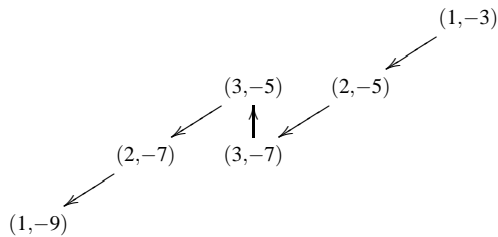
\[
\def\objectstyle{\scriptstyle}
\xymatrix@-1.0pc{
&&&&&(\bun,-3)\ar[ld]
\\
&&&\ar[ld] (\btrois,-5)& \ar[ld](\bdeux,-5)&
\\
&&{(\bdeux,-7)}\ar[ld] &\ar[u] (\btrois,-7)&&
\\
&(\bun,-9)
}
\]
The modules $K^{(\bdeux)}_{1,-11}$ and $K^{(\btrois)}_{1,-11}$ are:
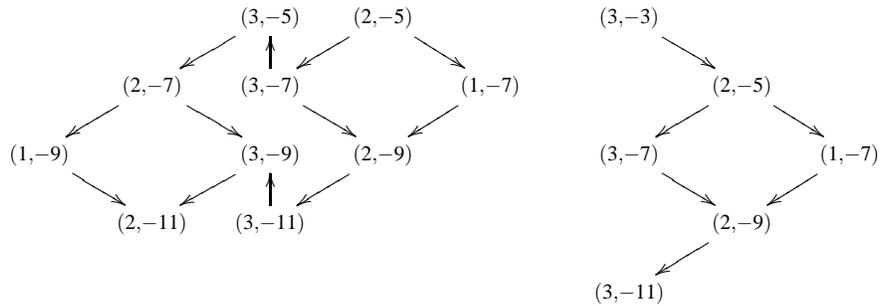
\[
\def\objectstyle{\scriptstyle}
\def\lablestyle{\scriptstyle}
\xymatrix@-1.0pc{
&&&(\btrois,-5)\ar[ld]&(\bdeux,-5)\ar[ld]\ar[rd]
\\
&&(\bdeux,-7)\ar[rd]\ar[ld]&(\btrois,-7)\ar[u]\ar[rd]&&(\bun,-7)\ar[ld]
\\
&(\bun,-9)\ar[rd]&&\ar[ld] (\btrois,-9)& \ar[ld](\bdeux,-9)&
\\
&&{(\bdeux,-11)} &\ar[u] (\btrois,-11)&&
}
\qquad
\xymatrix@-1.0pc{
(\btrois,-3)\ar[rd]
\\
&(\bdeux,-5)\ar[ld]\ar[rd]
\\
(\btrois,-7)\ar[rd]&&(\bun,-7)\ar[ld]
\\
& \ar[ld](\bdeux,-9)&
\\
(\btrois,-11)&&
}
\]
The corresponding fundamental $U_q(\hg)$-modules are $L(Y_{\bun,-11})$,
$L(Y_{\bdeux,-13})$, and $L(Y_{\btrois,-12})$, of respective dimensions
7, 22, and 8.

\subsection{Examples of $A$-modules for $\g$ of type $C_3$}\label{append-C3}

Let $\g$ is of type $C_3$, with the long
root being $\a_{\btrois}$.
The quiver $\G^-$ is: 
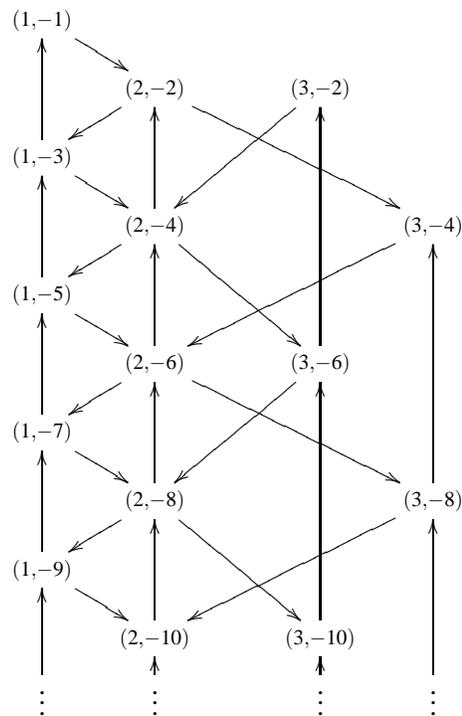
\[
\def\objectstyle{\scriptstyle}
\xymatrix@-1.0pc{
(\bun,-1)\ar[rd]&
\\
&\ar[ld](\bdeux,-2)\ar[rrrdd]&&(\btrois,-2)\ar[lldd]& 
\\
\ar[uu](\bun,-3)\ar[rd]&
\\
&\ar[uu](\bdeux,-4)\ar[ld]\ar[rrdd]&&&\ar[llldd](\btrois,-4) 
\\
\ar[uu](\bun,-5)\ar[rd]&
\\
&\ar[uu](\bdeux,-6)\ar[rrrdd]\ar[ld]&&\ar[uuuu](\btrois,-6)\ar[lldd]& 
\\
\ar[uu](\bun,-7)\ar[rd]&
\\
&\ar[uu](\bdeux,-8)\ar[ld]\ar[rrdd]&&&\ar[uuuu](\btrois,-8)\ar[llldd]
\\
\ar[uu](\bun,-9)\ar[rd]&
\\
&\ar[uu](\bdeux,-10)&&\ar[uuuu](\btrois,-10)&  
\\
{}\save[]+<0cm,0ex>*{\vdots}\ar[uu]\restore&{}\save[]+<0cm,0ex>*{\vdots}\ar[u]\restore&&{}\save[]+<0cm,0ex>*{\vdots}\ar[u]\restore
&{}\save[]+<0cm,0ex>*{\vdots}\ar[uuu]\restore
}
\]
The modules $K^{(\bun)}_{1,-7}$ and $K^{(\bdeux)}_{1,-8}$ are:
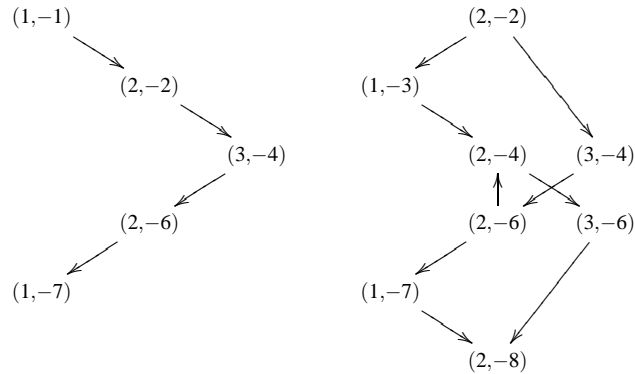
\[
\def\objectstyle{\scriptstyle}
\def\lablestyle{\scriptstyle}
\xymatrix@-1.0pc{
(\bun,-1)\ar[rd]
\\
&(\bdeux,-2)\ar[rd]
\\
&&(\btrois,-4)\ar[ld]
\\
& \ar[ld](\bdeux,-6)&
\\
(\bun,-7)&&
}
\qquad
\xymatrix@-1.0pc{
&\ar[ld](\bdeux,-2)\ar[rdd]
\\
(\bun,-3)\ar[rd]&
\\
&(\bdeux,-4)\ar[rd]&(\btrois,-4)\ar[ld]
\\
&\ar[u](\bdeux,-6)\ar[ld]&\ar[ldd](\btrois,-6)
\\
(\bun,-7)\ar[rd]&
\\
&(\bdeux,-8)
}
\]
The module $K^{(\btrois)}_{1,-8}$ is:
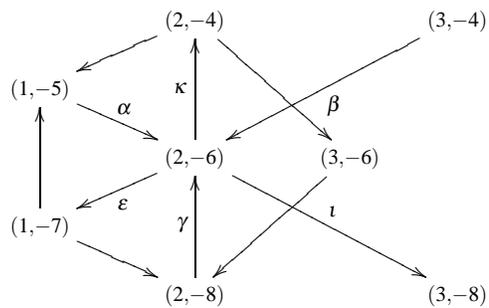
\[
\def\objectstyle{\scriptstyle}
\xymatrix@-1.0pc{
&&(\bdeux,-4)\ar[lld]\ar[rrdd]&&&\ar[llldd]^{\b}(\btrois,-4)
\\
(\bun,-5)\ar[rrd]^{\a}&
\\
&&\ar[uu]^{\kappa}(\bdeux,-6)\ar[rrrdd]^{\iota}\ar[lld]^{\epsilon}&&(\btrois,-6)\ar[lldd]& 
\\
\ar[uu](\bun,-7)\ar[rrd]&
\\
&&\ar[uu]^{\ga}(\bdeux,-8)&&&(\btrois,-8)
}
\]
Here, the vector space sitting at vertex $(\bdeux,-6)$ has dimension 2.
The maps incident to this space are given by the following matrices:
\[
\a = \pmatrix{1 \cr 0},\  
\b = \pmatrix{0 \cr 1},\  
\ga = \pmatrix{1 \cr 0},\ 
\epsilon = \pmatrix{0 & 1},\ 
\kappa = \pmatrix{0 & 1},\
\iota = \pmatrix{1 & 0}.
\]
The corresponding fundamental $U_q(\hg)$-modules are $L(Y_{\bun,-8})$,
$L(Y_{\bdeux,-10})$, and $L(Y_{\btrois,-10})$, of respective dimensions
6, 14, and 14.

\subsection{Examples of $A$-modules for $\g$ of type $F_4$} \label{append-F4}

Let $\g$ be of type $F_4$. 
We label the simple roots $\a_\bun$, $\a_\bdeux$, $\a_\btrois$, $\a_\bquat$,
so that the short simple roots are $\a_{\bun}$ and $\a_\bdeux$.
The quiver $\G^-$ is: 
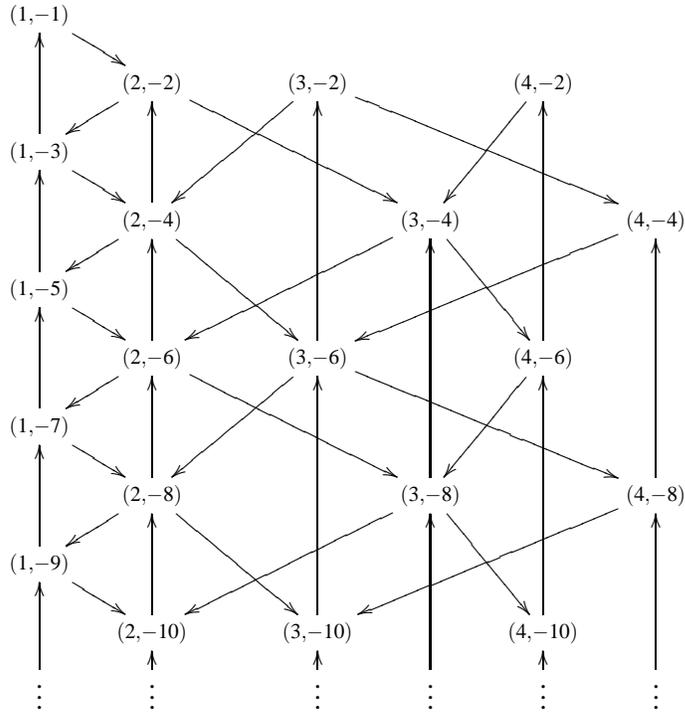
\[
\def\objectstyle{\scriptstyle}
\xymatrix@-1.0pc{
(\bun,-1)\ar[rd]&
\\
&\ar[ld](\bdeux,-2)\ar[rrrdd]&&(\btrois,-2)\ar[lldd]\ar[rrrdd]&&(\bquat,-2)\ar[ldd] 
\\
\ar[uu](\bun,-3)\ar[rd]&
\\
&\ar[uu](\bdeux,-4)\ar[ld]\ar[rrdd]&&&\ar[llldd](\btrois,-4)\ar[rdd]&&(\bquat,-4)\ar[llldd] 
\\
\ar[uu](\bun,-5)\ar[rd]&
\\
&\ar[uu](\bdeux,-6)\ar[rrrdd]\ar[ld]&&\ar[uuuu](\btrois,-6)\ar[lldd]\ar[rrrdd]&&\ar[uuuu](\bquat,-6)\ar[ldd] 
\\
\ar[uu](\bun,-7)\ar[rd]&
\\
&\ar[uu](\bdeux,-8)\ar[ld]\ar[rrdd]&&&\ar[uuuu](\btrois,-8)\ar[llldd]\ar[rdd]&&\ar[uuuu](\bquat,-8)\ar[llldd]
\\
\ar[uu](\bun,-9)\ar[rd]&
\\
&\ar[uu](\bdeux,-10)&&\ar[uuuu](\btrois,-10)&  
&\ar[uuuu](\bquat,-10)&
\\
{}\save[]+<0cm,0ex>*{\vdots}\ar[uu]\restore&{}\save[]+<0cm,0ex>*{\vdots}\ar[u]\restore&
&{}\save[]+<0cm,0ex>*{\vdots}\ar[u]\restore
&{}\save[]+<0cm,0ex>*{\vdots}\ar[uuu]\restore&{}\save[]+<0cm,0ex>*{\vdots}\ar[u]\restore
&{}\save[]+<0cm,0ex>*{\vdots}\ar[uuu]\restore
}
\]
The module $K^{(\bun)}_{1,-17}$ is:
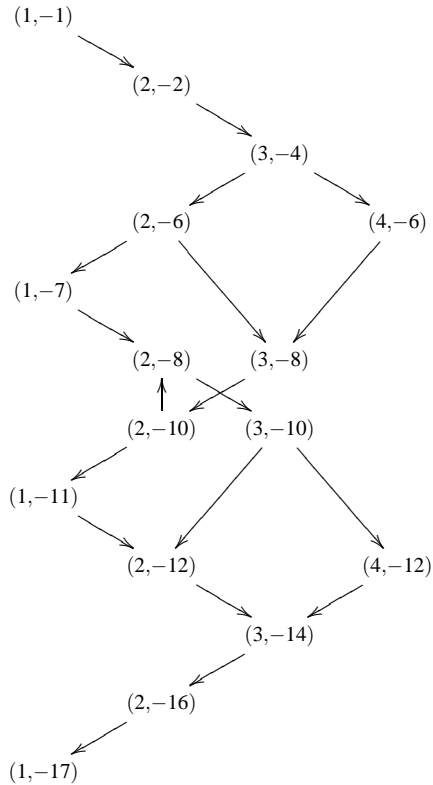
\[
\def\objectstyle{\scriptstyle}
\def\lablestyle{\scriptstyle}
\xymatrix@-1.0pc{
(\bun,-1)\ar[rd]
\\
&(\bdeux,-2)\ar[rd]
\\
&&(\btrois,-4)\ar[ld]\ar[rd]
\\
&\ar[ld](\bdeux,-6)\ar[rdd]&&(\bquat,-6)\ar[ldd]
\\
(\bun,-7)\ar[rd]&
\\
&(\bdeux,-8)\ar[rd]&(\btrois,-8)\ar[ld]
\\
&\ar[u](\bdeux,-10)\ar[ld]&\ar[ldd](\btrois,-10)\ar[rdd]
\\
(\bun,-11)\ar[rd]&
\\
&(\bdeux,-12)\ar[rd]&&(\bquat,-12)\ar[ld]
\\
&&(\btrois,-14)\ar[ld]
\\
&(\bdeux,-16)\ar[ld]
\\
(\bun,-17)
}
\]
The module $K^{(\bquat)}_{1,-16}$ is:
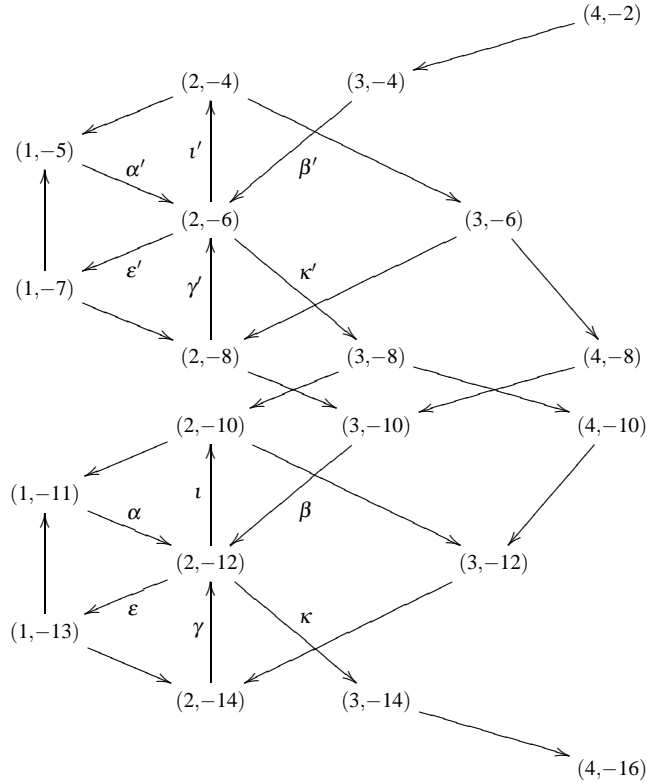
\[
\def\objectstyle{\scriptstyle}
\def\lablestyle{\scriptstyle}
\xymatrix@-1.0pc{
&&&&&&(\bquat,-2)\ar[lld]
\\
&&(\bdeux,-4)\ar[lld]\ar[rrrdd]&&\ar[lldd]^{\b'}(\btrois,-4)
\\
(\bun,-5)\ar[rrd]^{\a'}&
\\
&&\ar[uu]^{\iota'}(\bdeux,-6)\ar[rrdd]^{\kappa'}\ar[lld]^{\epsilon'}&&&(\btrois,-6)\ar[llldd]\ar[rdd]& 
\\
\ar[uu](\bun,-7)\ar[rrd]&
\\
&&\ar[uu]^{\ga'}(\bdeux,-8)\ar[rrd]&&\ar[lld](\btrois,-8)\ar[rrd]&&\ar[lld](\bquat,-8)
\\
&&(\bdeux,-10)\ar[lld]\ar[rrrdd]&&\ar[lldd]^{\b}(\btrois,-10)&& (\bquat,-10)\ar[ldd]
\\
(\bun,-11)\ar[rrd]^{\a}&
\\
&&\ar[uu]^{\iota}(\bdeux,-12)\ar[rrdd]^{\kappa}\ar[lld]^{\epsilon}&&&(\btrois,-12)\ar[llldd]& 
\\
\ar[uu](\bun,-13)\ar[rrd]&
\\
&&\ar[uu]^{\ga}(\bdeux,-14)&&(\btrois,-14)\ar[rrd]
\\
&&&&&&(\bquat,-16)
}
\]
Here, the vector spaces sitting at vertex $(\bdeux,-6)$ and $(\bdeux,-12)$ have dimension 2.
The maps incident to these spaces are given by the following matrices:
\[
\a = \pmatrix{1 \cr 0},\  \kappa = \pmatrix{1 & 0},\ 
\b = \pmatrix{0 \cr 1},\  \epsilon = \pmatrix{0 & 1},\ 
\ga = \pmatrix{1 \cr 0},\  \iota = \pmatrix{0 & 1},
\]
\[
\a' = \pmatrix{1 \cr 0},\  \kappa' = \pmatrix{1 & 0},\ 
\b' = \pmatrix{0 \cr 1},\  \epsilon' = \pmatrix{0 & 1},\ 
\ga' = \pmatrix{1 \cr 0},\  \iota' = \pmatrix{0 & 1}.
\]
The corresponding fundamental $U_q(\hg)$-modules are $L(Y_{\bun,-18})$
and $L(Y_{\bquat,-18})$, of respective dimensions
26, and 53.

\subsection*{Acknowledgements}
We are grateful to P.-G. Plamondon for very helpful comments and discussions.


\smallskip
\small
\noindent
\begin{tabular}{ll}
David {\sc Hernandez} : &
Universit\'e Paris Diderot-Paris 7,\\
& Institut de Math\'ematiques de Jussieu - Paris Rive Gauche CNRS UMR 7586,\\
& B\^atiment Sophie Germain, Case 7012,\\ 
& 75205 Paris Cedex 13, France. \\
&email : {\tt hernandez@math.jussieu.fr}\\[5mm]
Bernard {\sc Leclerc} : & Normandie Univ, France\\ 
&UNICAEN, LMNO F-14032 Caen, France\\
&CNRS UMR 6139, F-14032 Caen, France\\
&Institut Universitaire de France.\\
&email : {\tt bernard.leclerc@unicaen.fr}
\end{tabular}


\begin{thebibliography}{ABCD}


\footnotesize


\bibitem[\bf C]{C}
{\sc V. Chari}, {\em Minimal affinizations of representations of quantum groups: the rank $2$ case},
Publ. Res. Inst. Math. Sci.  {\bf 31}  (1995), 873--911.

\bibitem[\bf CH]{CH} {\sc V. Chari, D. Hernandez},
\emph{Beyond Kirillov-Reshetikhin modules}, 
In Quantum affine algebras, extended affine Lie algebras, 
and their applications,  
Contemp. Math., {\bf 506}, AMS Providence, 2010, 49--81.

\bibitem[\bf CP1]{CP} {\sc V. Chari, A. Pressley},
\emph{A guide to quantum groups}.
Cambridge University Press 1994.

\bibitem[\bf CP2]{CP2}{\sc V. Chari, A. Pressley},
{\it  Minimal affinizations of representations of quantum groups: the
simply laced case},
J. Algebra {\bf 184} (1996), 1--30. 

\bibitem[\bf DWZ1]{DWZ}
{\sc H. Derksen, J. Weyman, A. Zelevinsky},
{\em Quivers with potential and their representations I: Mutations}, 
Selecta Math., {\bf 14} (2008), 59--119. 

\bibitem[\bf DWZ2]{DWZ2}
{\sc H. Derksen, J. Weyman, A. Zelevinsky},
{\em Quivers with potential and their representations II: 
Applications to cluster algebras}, 
J. Amer. Math. Soc. {\bf 23} (2010), 749--790.

\bibitem[\bf FM]{FM} 
{\sc E. Frenkel, E. Mukhin}, 
\emph{Combinatorics of $q$-characters of finite-dimensional
representations of quantum affine algebras},  
Comm. Math. Phys. \textbf{216}  (2001),  23--57.

\bibitem[\bf FR]{FR}{\sc E. Frenkel, N. Reshetikhin},
{\it The $q$-characters of representations of quantum affine algebras},
Recent developments in quantum affine algebras and related topics,
Contemp. Math. {\bf 248} (1999), 163--205.

\bibitem[\bf FZ1]{FZ1}{\sc S. Fomin, A. Zelevinsky}, 
\emph{Cluster algebras I: Foundations},
J. Amer. Math. Soc. {\bf 15} (2002), 497--529.

\bibitem[\bf FZ2]{FZsurvey}{\sc S. Fomin, A. Zelevinsky}, 
\emph{Cluster algebras: notes for the CDM-03 conference}, in  
Current developments in mathematics, 2003,  1--34, Int. Press,
Somerville, MA, 2003.

\bibitem[\bf FZ3]{FZ}{\sc S. Fomin, A. Zelevinsky},
{\it Cluster algebras IV: coefficients},
Compos. Math.  {\bf 143}  (2007),  112--164.

\bibitem[\bf GG]{GG}{\sc J. Grabowski, S. Gratz}
\emph{Cluster algebras of infinite rank},
(2012) arxiv 1212.3528.

\bibitem[\bf GLS1]{GLS1}
{\sc C. Geiss, B. Leclerc, J. Schr\"oer},
\emph{Auslander algebras and initial seeds for cluster algebras},
J. London Math. Soc. {\bf 75} (2007), 718--740.

\bibitem[\bf GLS2]{GLS2}
{\sc C. Geiss, B. Leclerc, J. Schr\"oer},
\emph{Kac-Moody groups and cluster algebras},
Advances Math. {\bf 228} (2011), 329--443.

\bibitem[\bf GSV]{GSV}
{\sc M. Gekhtman, M. Shapiro, A. Vainshtein}, 
{\it Cluster algebras and Poisson geometry}, 
AMS Math. Survey and Monographs {\bf 167}, AMS 2010.

\bibitem[\bf H]{H} {\sc D. Hernandez},
\emph{The Kirillov-Reshetikhin conjecture and solutions of $T$-systems},
J. Reine Angew. Math. \textbf{596} (2006), 63--87.

\bibitem[\bf HJ]{HJ}{\sc D. Hernandez, M. Jimbo},
\emph{Asymptotic representations and Drinfeld rational fractions},
Compositio Math. {\bf 148} (2012), 1593--1623.

\bibitem[\bf HL1]{HL}
{\sc D. Hernandez, B. Leclerc},
{\em Cluster algebras and quantum affine algebras}, 
Duke Math. J. {\bf 154} (2010), 265--341.

\bibitem[\bf HL2]{HL2}
{\sc D. Hernandez, B. Leclerc},
{\em Monoidal categorifications of cluster algebras of type $A$ and $D$}, 
in Symmetries, integrable systems and representations, (K. Iohara, S. Morier-Genoud, B. R\'emy, eds.),
Springer proceedings in mathematics and statistics {\bf 40}, 2013, 175--193.

\bibitem[\bf HoJo]{HoJo}
{\sc T. Holm, P. Jorgensen},
\emph{On a cluster category of infinite Dynkin type, and the relation to triangulations of the infinitygon},
Math. Z. {\bf 270} (2012), 277--295.

\bibitem[\bf IIKKN1]{IIKKN1}
{\sc R. Inoue, O. Iyama, B. Keller, A. Kuniba, T. Nakanishi},
\emph{Periodicities of $T$ and $Y$-systems, dilogarithm identities,
and cluster algebras I, Type $B_r$}, 
Publ. RIMS, (to appear).

\bibitem[\bf IIKKN2]{IIKKN2}
{\sc R. Inoue, O. Iyama, B. Keller, A. Kuniba, T. Nakanishi},
\emph{Periodicities of $T$ and $Y$-systems, dilogarithm identities,
and cluster algebras II, Type $C_r, F_4$ and $G_2$}, 
Publ. RIMS, (to appear).

\bibitem[\bf IIKNS]{IIKNS}
{\sc R. Inoue, O. Iyama, A. Kuniba, T. Nakanishi, J. Suzuki},
\emph{Periodicities of $T$ and $Y$-systems},
Nagoya Math. J. {\bf 197} (2010), 59--174.


\bibitem[\bf Ka]{K}
{\sc V. G. Kac},
{\em Infinite dimensional Lie algebras}, Cambridge University Press 1990.


\bibitem[\bf KNS1]{KNS1}
{\sc A. Kuniba, T. Nakanishi, J. Suzuki},
\emph{Functional relations in solvable lattice models: I. Functional relations and representation
theory}, Int. J. Mod. Phys. A{\bf9} (1994), 5215--5266.

\bibitem[\bf KNS2]{KNS2}
{\sc A. Kuniba, T. Nakanishi, J. Suzuki},
\emph{$T$-systems and $Y$-systems in integrable systems},
J. Phys. A {\bf 44} (2011), no. 10, 103001, 146 pp.

\bibitem[\bf Le1]{Le-imag}{\sc B. Leclerc},
{\it Imaginary vectors in the dual canonical basis of
  $U_q(\mathfrak{n})$},
Transform. Groups  {\bf 8}  (2003),  95--104. 

\bibitem[\bf Le2]{L}{\sc B. Leclerc},
{\it Quantum loop algebras, quiver varieties, and cluster algebras},
in Representations of Algebras and Related Topics, (A. Skowro\'nski and K. Yamagata, eds.),
European Math. Soc. Series of Congress Reports, 2011, 117--152.

\bibitem[\bf Lu]{Lu2}{\sc G. Lusztig},
{\it On quiver varieties}, Adv. Math. {\bf 136} (1998), 141--182.

\bibitem[\bf N1]{N1}{\sc H. Nakajima}, 
{\em Quiver varieties and finite-dimensional representations of 
quantum affine algebras}.  
J. Amer. Math. Soc.  {\bf 14}  (2001),  145--238.

\bibitem[\bf N2]{N2}{\sc H. Nakajima}, 
\emph{$t$-analogs of $q$-characters of Kirillov-Reshetikhin modules 
of quantum affine algebras},  
Represent. Theory  \textbf{7}  (2003), 259--274. 

\bibitem[\bf N3]{NAnnals}{\sc H. Nakajima}, 
{\em Quiver varieties and $t$-analogs of $q$-characters of 
quantum affine algebras}.  
Annals of Math.  {\bf 160}  (2004),  1057--1097.

\bibitem[\bf N4]{N3}{\sc H. Nakajima}, 
{\em Quiver varieties and cluster algebras},
Kyoto J. Math. {\bf 51} (2011), 71--126. 

\bibitem[\bf Pa]{Pa}{\sc Y. Palu},
\emph{Cluster characters II: A multiplication formula},
Proceedings London Math. Soc. (to appear).

\bibitem[\bf Pl1]{P1}
{\sc P.-G. Plamondon},
\emph{Cluster characters for cluster categories with infinite-dimensional morphism spaces},
Adv. Math. {\bf 227} (2011), 1--39.

\bibitem[\bf Pl2]{P}
{\sc P.-G. Plamondon},
{\it Generic bases for cluster algebras from the cluster category},
Int. Math. Research. Notices, (to appear).   
\end{thebibliography}
\end{document}